\documentclass[12pt]{article}
\usepackage{amsmath,amsthm,amsfonts,amssymb}

\usepackage[usenames]{color} 

\newtheorem{thm}{Theorem}[section]

\newtheorem{lemma}[thm]{Lemma}

\newtheorem{cor}[thm]{Corollary}
\newtheorem{corollary}[thm]{Corollary}
\newtheorem{proposition}[thm]{Proposition}
\newtheorem{question}[thm]{Question}

\newtheorem{prop}[thm]{Proposition}

\newtheorem{sublem}[thm]{Sublemma}
\theoremstyle{definition}

\newtheorem{definition}[thm]{Definition} 

\newtheorem{notn}[thm]{Notation}
\newtheorem{com}[thm]{Remark}
\newtheorem{ex}[thm]{Example}

\newtheorem{remark}[thm]{Remark}
\theoremstyle{remark}

\newcounter{remarks}
{\paragraph*{Remarks}\smallskip
     \begin{list}{\arabic{remarks}. }{\usecounter{remarks}%
          \setlength{\leftmargin}{0in}%
          \setlength{\rightmargin}{0in}%
          \setlength{\labelsep}{0pt}%
          \setlength{\labelwidth}{0pt}%
          \setlength{\listparindent}{0pt}%
     }
}
{
\end{list}
}

\DeclareMathOperator{\Fix}{Fix}
\DeclareMathOperator{\Per}{Per}
\DeclareMathOperator{\PA}{P}

\DeclareMathOperator{\PF}{PF}
\DeclareMathOperator{\EL}{EL}
\DeclareMathOperator{\IA}{IA_n}

\newcommand{\va}{{\mathbf a}}
\newcommand{\vb}{{\mathbf b}}

\newcommand{\R}{{\mathbb R}}
\newcommand{\Z}{{\mathbb Z}}

\newcommand{\T}{{\mathbb T}}

\newcommand{\f}{F_n}
\newcommand{\D}{{\mathcal D}}
\newcommand{\E}{{\mathcal E}}
\newcommand{\V}{\mathcal V}
\newcommand{\oone}{\phi} 
\newcommand{\otwo}{\psi}
\newcommand{\othree}{\theta}
\newcommand{\aone}{\Phi}
\newcommand{\atwo}{\Psi}
\newcommand{\athree}{\Theta}
\newcommand{\Out}{\mathsf{Out}}
\newcommand{\Aut}{\mathsf{Aut}}
\newcommand{\Inn}{\mathsf{Inn}}
\newcommand{\mcg}{\mathsf{MCG}}
\newcommand{\Stab}{\mathsf{Stab}}
\newcommand{\ffs}{free factor system}

\newcommand{\EG}{EG}

\newcommand{\F}{\mathcal F}
\newcommand{\rtt}{relative train track map}

\def\L{\mathcal L}

\newcommand{\A}{\mathcal A}

\newcommand{\fG} {f : G \to G}
\newcommand{\ti} {\tilde}
\newcommand{\iNp} {indivisible Nielsen path}

\newcommand{\filt}{\emptyset = G_0 \subset G_1 \subset \ldots  \subset G_N = G}

\newcommand{\ofn}{\Out(F_n)}
\newcommand{\eg}{EG}
\newcommand{\noneg}{NEG}

\newcommand{\cs}{completely split}

\newcommand{\ds}{$QE$-splitting}
\newcommand{\qe}{quasi-exceptional}

\newcommand{\fpss}{FPS subgraph}

\newcommand{\wt}{\widetilde}

\newcommand{\ias}{almost invariant subgraph}

\newcommand{\comment}[1]{}

\title{ Abelian subgroups of $\Out(F_n)$}
\author{Mark Feighn\thanks{This material is based upon work supported by the
National Science Foundation under Grant No.~DMS0504917.}\and
Michael Handel\thanks{This material is based upon work supported by
the National Science Foundation under Grant No.~DMS0405814.}}

\begin{document}
\maketitle
\begin{abstract}
We classify abelian subgroups of $\Out(F_n)$ up to finite index in an algorithmic and computationally friendly way.  A process called disintegration is used to canonically decompose a single rotationless element $\phi$ into  a composition of finitely many elements and then use these elements to generate an abelian subgroup $\A(\phi)$ that contains $\phi$.   The main theorem is that up to finite index every abelian subgroup is realized by this construction.  As an application  we classify, up to
finite index, abelian subgroups of $\Out(F_n)$ and of $\IA$ with
maximal rank. 
\end{abstract}

\tableofcontents \label{disint}
\section{Introduction}

In this paper we classify abelian subgroups of $\Out(F_n)$ up to
finite index in an algorithmic and computationally friendly way.
There are two steps.  The first is to construct an abelian subgroup
$\D(\phi)$ from a given $\phi \in \Out(F_n)$ by a process that we call
{\em disintegration}.  The subgroup $\D(\phi)$ is very well understood
in terms of relative train track maps and has natural coordinates that
embed it into some $\Z^M$.  The second step is to prove the following
theorem.

\vspace{.1in}

\noindent {\bf Theorem~\ref{disintegration groups are enough}.}\ \ {\em  
For every abelian subgroup $A$ of $\Out(F_n)$ there exists $\phi \in
A$ such that $A \cap D(\oone)$ has finite index in $A$.}

\vspace{.1in}

To motivate the disintegration process, consider a pure element $\mu$
of the mapping class group $\mcg(S)$ of a compact oriented surface
$S$.  By the Thurston classification theorem
\cite{wpt:pseudo},\cite{flp:thurston}, there is a decomposition of $S$
into subsurfaces $S_l$, some of which are annuli and the rest of which
have negative Euler characteristic, and there is a homeomorphism $h :
S \to S$ representing $\mu$, called a {\em normal form} for $\mu$,
that preserves each $S_l$.  If $S_l$ is an annulus then $h|S_l$ is a
non-trivial Dehn twist.  If $S_l$ has negative Euler characteristic
then $h|S_l$ is either the identity or is pseudo-Anosov.  In all
cases, $h|\partial S_l$ is the identity.

We may assume that the $S_l$'s are numbered so that $h|S_l$ is the
identity if and only if $l > M$ for some $M$.  For each $M$-tuple of
integers $\va = (a_1,\ldots,a_M)$ let $h_{\va}: S \to S$ be the
homeomorphism that agrees with $h^{a_l}$ on $S_l$ for $1 \le l \le M$
and is the identity on the remaining $S_l$'s.  Then $h_{\va}$ is a
normal form for an element $\mu_{\va} \in \mcg(S)$ and we define
$\D(\mu)$ to be the subgroup consisting of all such $\mu_{\va}$.  It
is easy to check that $ \mu_{\va} \to \va$ defines an isomorphism
between $\D(\mu)$ and $\Z^M$.

An element $\phi$ of $\Out(F_n)$ has finite sets of
natural invariants on which it acts by permutation.  If these actions
are trivial then we say that $\phi$ is {\it rotationless}; complete
details can be found in section~\ref{s:lifts}.  This property is
similar to being pure, which is defined as acting trivially on
$H_1(F_n,\Z_3)$.  Abelian subgroups are both virtually rotationless
and virtually pure.  The latter is obvious and the former follows from
Corollary~\ref{is birot}.  We work in the rotationless category since
it is more natural for our constructions.

Suppose that $\phi$ is a rotationless element of $\Out(F_n)$.  The
analog of a normal form $h :S \to S$ is a relative train track map
$\fG$, which is a particularly nice homotopy equivalence of a marked
graph that represents $\phi$ in the sense that the outer automorphism
of $\pi_1(G)$ that it induces is identified with $\phi$ by
the marking.  There is an associated maximal filtration $\filt$ by
$f$-invariant subgraphs.  The $i^{th}$ stratum $H_i$ is the closure of
$G_i \setminus G_{i-1}$.  The exact properties satisfied by $\fG$ and
$\filt$ are detailed in section~\ref{s:background}.

As a first attempt to mimic the construction of $\D(\mu)$, let
$X_1,\ldots, X_M$ be the strata that are not fixed by $f$, let $\va =
(a_1,\ldots,a_M)$ be an $M$-tuple of non-negative integers and define
$f_{\va}$ to agree with $f^{a_l}$ on $X_l$ and to be the identity on
the subgraph of edges fixed by $f$.  Although it is not obvious,
$f_{\va} : G \to G$ is a homotopy equivalence (see Lemma~\ref{homeq})
and so defines an element $\phi_{\va} \in \Out(F_n)$.

Without some restrictions on $\va$ however, the subgroup generated by
the $\phi_{\va}$'s need not be abelian.  In the following examples, we
do not distinguish between a homotopy equivalence of the rose and the
outer automorphism that it represents.

\begin{ex}\label{ex:not invariant}  Let $G$ be the graph with one vertex and with edges labelled $A,B$ and $C$.   Define $\fG$ by  
$$ A\mapsto A \qquad B\mapsto BA \qquad C\mapsto CB.
$$
Let $X_1=\{B\}$ and $X_2=\{C\}$ and $\va=(m,n)$.   Then 
$$
f_{(m,n)}\circ f(C)=f_{(m,n)}(CB)=f^n(C)f^m(B)
$$
and 
$$
f \circ f_{(m,n)}(C)=f(f^n(C))=f^n(f(C))=f^n(CB)=f^n(C)f^n(B).
$$ This shows that $f_{(m,n)}$ commutes with $f = f_{(1,1)}$ if and
only if $m=n$.

The underlying problem is that strata are not invariant.  It does not
matter that the path $f(B)$ crosses $A$ since $A$ is fixed by $f$.
The lack of commutativity stems from the fact that $f(C)$ crosses $B$.
\end{ex}

To address this problem we enlarge the $X_i$'s to be unions of
strata. It is not necessary to choose the $X_i$'s to be fully
invariant (i.e. to satisfy $f(X_i) \subset X_i$) but they must be {\em
almost invariant} as made precise in Definition~\ref{almost
invariant}.

The next example illustrates a more complicated relation on the
coordinates of $\va$ that is needed to insure that the $f_{\va}$'s
commute.

\begin{ex} \label{quasiexceptional} Let $G$ be the graph with one vertex and with edges labelled $A,B,C, D$ and $E$.   Define $\fG$ by  
$$ A\mapsto A \qquad B\mapsto BA^2 \qquad C \mapsto CB \qquad D\mapsto
DA^5 \qquad E \mapsto DC\bar B
$$
where $\bar B$ is $B$ with its orientation reversed. Let $X_1 =\{B,C\}, X_2 = \{D\}$ and $X_3=\{E\}$ and let  $\va = (m,n,p)$. 
Then 
$$
f\circ f_{\va}(D)=f(f^p(D))=f^p(f(D))=f^p(DC\bar B) =f^p(D)f^p(C \bar B)
$$
and
$$
f_{\va}\circ f(D)=f_{\va}(DC\bar B)=f^p(D)f_{\va}(C\bar B). 
$$
If $f$ commutes with $f_{\va}$ then 
$$
f^p(C \bar B) = f_{\va}(C\bar B). 
$$ Thus $CA^{3p}\bar B = CA^{5n-2m} \bar B$ and $3p=5n-2m$.  One can
check that the converse holds as well.  Namely if we require that
$\va$ be an element of the linear subspace of $\Z^3 = \{(m,n,p)\}$
defined by $3p=5n-2m$ then the $\phi_{\va}$'s commute.
\end{ex}

The path $C \bar B$ of Example~\ref{quasiexceptional} is {\em
quasi-exceptional} as defined in section~\ref{s:disintegration}.  When
the image of an edge in $X_k$ contains a quasi-exceptional path with
initial edge in $X_i$ and terminal edge in $X_j$ then there is an
induced relation between the $i^{th}, j^{th}$ and $k^{th}$ coefficient
of $\va$.  These define a linear subspace of $\Z^M$.  The non-negative
$M$-tuples that lie in this subspace are said to be {\em admissible}.
The map $\va \to \phi_{\va}$ on admissible $M$-tuples extends to an
injective homomorphism of the full linear subspace and we define the
image of this subspace to be $\D(\phi)$.

The mapping class group version of Theorem~\ref{disintegration groups
are enough} is a straightforward consequence of two easily proved,
well known facts.  The first (see for example Corollary~5.2 of
\cite{fhp:sabel}) is that the subsurfaces $S_l$ can be chosen
independently of $\mu \in A$.  The second (see for example Lemma~2.10
of \cite{fhp:abel}) is that an abelian subgroup containing a
pseudo-Anosov element is virtually cyclic.

The proof for $\Out(F_n)$ is considerably harder.  This is due, in
part, to the fact that disintegration in $\Out(F_n)$ is a more
complicated operation, as illustrated by the examples, than it is
$\mcg(S)$.  Another factor is that, unlike normal forms in the mapping
class group, relative train track maps representing an element $\phi
\in \Out(F_n)$ are not unique.  No matter how canonical a construction
is with respect to a particular $\fG$, one must still check the extent
to which it is independent of the choice of $\fG$.  The most
technically difficult argument in this paper (section~\ref{s:finite
index}) is a proof that the rank of the admissible linear subspace of
$\Z^M$ described above depends only on $\phi$ and not the choice of
$\fG$.

Recall that $\IA$ is the subgroup of $\Out(F_n)$ consisting of
elements that act as the identity on $H_1(F_n)$.  As an application of
Theorem~\ref{disintegration groups are enough} we classify, up to
finite index, abelian subgroups of $\Out(F_n)$ and of $\IA$ with
maximal rank.  The exact statements appear as Proposition~\ref{max
abel} and Proposition~\ref{IA max abel}.  Roughly speaking, we prove
that if $\D(\phi)$ has maximal rank in $\Out(F_n)$ then $\fG$ has
$2n-3$ strata, each of which is either a single linear edge or is
exponentially growing and is closely related to a pseudo-Anosov
homeomorphism of a four times punctured sphere.  If $\D(\phi)$ has
maximal rank in $\IA$ then $\fG$ has $2n-4$ such strata and pointwise
fixes a rank two subgraph.

From an algebraic point of view, the natural abelian subgroup
associated to an element $\phi \in \Out(F_n)$ is the center
$Z(C(\phi))$ of the centralizer $C(\phi)$ of $\phi$ which can also be
described as the intersection of all maximal (with respect to
inclusion) abelian subgroups that contain $\phi$.  In our context it
is natural to look at the weak center $WZ(C(\phi))$ of $C(\phi)$
defined as the subgroup of elements that commute with an iterate of
each element of $C(\phi)$.  The following result is a step toward an
algorithmic construction of $Z(C(\phi))$.  Every abelian subgroup $A$
has a finite index subgroup $A_R$, all of whose elements are
rotationless.

\vspace{.1in}

\noindent {\bf Theorem~\ref{WZC}.}\ {\em  
$\D_R(\phi) \subset WZ(C(\phi))$ for all rotationless $\phi$.}

\vspace{.1in}

In section~\ref{s:for commensurator} we apply this theorem to give
algebraic characterizations of certain maximal rank abelian subgroups
of $\Out(F_n)$ and $\IA$.  This characterization is needed in the
calculation of the commensurator group of $\Out(F_n)$
\cite{FarbH:commensurator}.

In section~\ref{s:lifts} we define what it means for $\oone \in
\Out(F_n)$ to be {\em rotationless}, prove that the rotationless
elements of any abelian subgroup $A$ form a finite index subgroup
$A_R$ and consider lifts of $A_R$ from $\Out(F_n)$ to $\Aut(F_n)$.
These lifts are essential to our approach and are similar to ones used
in \cite{bfh:tits3}.

In section~\ref{s:generic} we define a natural embedding of $A_R$ into a
lattice in Euclidean space and say what it means for an element of $A_R$
to be {\em generic} with respect to this embedding.  

In section~\ref{s:Aone} we associate an abelian subgroup $\A(\oone)$ to each rotationless
$\oone$ and prove that if $\oone$ is generic in $A_R$ then $A_R \subset
\A(\oone)$.    We also prove (Corollary~\ref{in the center}) that $\A(\phi) \subset WZ(C(\phi))$.

In section~\ref{s:disintegration}, we define $\D(\phi)$ and prove (Corollary~\ref{D is contained in A}) that
$\D_R(\phi) \subset \A(\oone)$, thereby completing the proof of Theorem~\ref{WZC}.    

In section~\ref{s:finite index} we prove  (Theorem~\ref{D has finite index}) that $\D_R(\phi)$ has finite index  in $\A(\oone)$ by reconciling the normal forms point
of view used to define $\D(\oone)$ with the \lq action on $\partial
F_n$\rq\ point of view used to define $\A(\oone)$.  Theorem~\ref{disintegration groups are enough} is an immediate consequence of this result and the fact, mentioned above, that  $A_R \subset
\A(\oone)$ for generic $\oone \in A$.  

We make use of several important results from \cite{fh:recognition},
including the Recognition Theorem and the existence of relative train
track maps that are especially well suited to disintegrating an
element and forming $\D(\oone)$.  Section~\ref{s:background} reviews
this and other relevant material and sets notation for the paper.

\section{Background} \label{s:background}
Fix $n \ge 2$ and let $F_n$ be the free group of rank $n$.  Denote the
automorphism group of $F_n$ by $\Aut(F_n)$, the group of inner
automorphisms of $F_n$ by $\Inn(F_n)$ and the group of outer
automorphisms of $F_n$ by $\Out(F_n) = \Aut(F_n)/\Inn(F_n)$.  We
follow the convention that elements of $\Aut(F_n)$ are denoted by
upper case Greek letters and that the same Greek letter in lower case
denotes the corresponding element of $\Out(F_n)$.  Thus $\aone \in
\Aut(F_n)$ represents $\oone \in \Out(F_n)$.

\paragraph{Marked Graphs and Outer Automorphisms}
Identify $F_n$ with $\pi_1(R_n,*)$ where $R_n$ is the rose with one
vertex $*$ and $n$ edges.  A {\em marked graph} $G$ is a graph of rank
$n$ without valence one vertices, equipped with a homotopy equivalence
$m : R_n \to G$ called a {\em marking}.  Letting $b = m(*) \in G$, the
marking determines an identification of $F_n$ with $\pi_1(G,b)$.

A homotopy equivalence $\fG$ and a path $\sigma$ from $b$ to $f(b)$
determines an automorphism of $\pi_1(G,b)$ and hence an element of
$\Aut(F_n)$.  As the homotopy class of $\sigma$ varies, the
automorphism ranges over all representatives of the associated outer
automorphism $\oone$.  We say that {\em$\fG$ represents $\phi$}.  We
always assume that the restriction of $f$ to any edge is an immersion.

\paragraph{Paths, Circuits and Edge Paths}
Let $\Gamma$ be the universal cover of a marked graph $G$ and let $pr
: \Gamma\to G$ be the covering projection.  A proper map $\ti \sigma :
J \to \Gamma$ with domain a (possibly infinite) interval $J$ will be
called a {\it path in $\Gamma$} if it is an embedding or if $J$ is
finite and the image is a single point; in the latter case we say that
$\ti \sigma$ is {\it a trivial path}.  If $J$ is finite, then $\ti
\sigma : J \to \Gamma$ is homotopic rel endpoints to a unique
(possibly trivial) path $[\ti \sigma]$; we say that {\it $[\ti
\sigma]$ is obtained from $\ti \sigma$ by tightening}. If $\ti f
:\Gamma \to \Gamma$ is a lift of a homotopy equivalence $\fG$, we
denote $ [\ti f(\ti \sigma)]$ by $\ti f_\#(\ti \sigma)$.

We will not distinguish between paths in $\Gamma$ that differ only
by an orientation preserving change of parametrization. Thus we are
interested in the oriented image of $\ti \sigma$ and not $\ti \sigma$
itself.  If the domain of $\ti \sigma$ is finite, then the image of
$\ti \sigma$ has a natural decomposition as a concatenation $\wt E_1
\wt E_2 \ldots \wt E_{k-1} \wt E_k$ where $\wt E_i$, $1 < i < k$, is
an edge of $\Gamma$, $\wt E_1$ is the terminal segment of an edge and
$\ti E_k$ is the initial segment of an edge. If the endpoints of the
image of $\ti \sigma$ are vertices, then $\wt E_1$ and $\wt E_k $ are
full edges.  The sequence $\wt E_1 \wt E_2\ldots \wt E_k$ is called
{\it the edge path associated to $\ti \sigma$.}  This notation extends
naturally to the case that the interval of domain is half-infinite or
bi-infinite.  In the former case, an edge path has the form $\wt
E_1\wt E_2\ldots$ or $\ldots \wt E_{-2} \wt E_{-1}$ and in the latter
case has the form $\ldots\wt E_{-1}\wt E_0\wt E_1\wt E_2\ldots$.

A {\it path in $G$} is the composition of the projection map $pr$ with
a path in $\Gamma$.  Thus a map $\sigma : J \to G$ with domain a
(possibly infinite) interval will be called a path if it is an
immersion or if $J$ is finite and the image is a single point; paths
of the latter type are said to be {\em trivial}. If $J$ is finite,
then $ \sigma : J \to G$ is homotopic rel endpoints to a unique
(possibly trivial) path $[ \sigma]$; we say that {\it $[ \sigma]$ is
obtained from $ \sigma$ by tightening}. For any lift $\ti \sigma : J
\to \Gamma$ of $\sigma$, $[\sigma] = pr[\ti \sigma]$.  We denote
$[f(\sigma)]$ by $f_\#(\sigma)$.  We do not distinguish between paths
in $G$ that differ by an orientation preserving change of
parametrization.  The {\it edge path associated to $\sigma$} is the
projected image of the edge path associated to a lift $\ti \sigma$.
Thus the edge path associated to a path with finite domain has the
form $ E_1 E_2 \ldots E_{k-1} E_k$ where $ E_i$, $1 < i < k$, is an
edge of $G$, $ E_1$ is the terminal segment of an edge and $ E_k$ is
the initial segment of an edge.  We will identify paths with their
associated edge paths whenever it is convenient.

We reserve the word {\it circuit} for an immersion $\sigma : S^1 \to
G$.  Any homotopically non-trivial map $\sigma : S^1 \to G$ is
homotopic to a unique circuit $[\sigma]$. As was the case with paths,
we do not distinguish between circuits that differ only by an
orientation preserving change in parametrization and we identify a
circuit $\sigma$ with a {\it cyclically ordered edge path}
$E_1E_2\dots E_k$.

A path or circuit \emph{crosses} or \emph{contains} an edge if that
edge occurs in the associated edge path.  For any path $\sigma$ in $G$
define $\bar \sigma$ to be \lq $\sigma$ with its orientation
reversed\rq.  For notational simplicity, we sometimes refer to the
inverse of $\ti \sigma$ by $\ti \sigma^{-1}$.

A decomposition of a path or circuit into subpaths is a {\em
splitting} for $\fG$ and is denoted $\sigma = \ldots
\sigma_1\cdot\sigma_2\ldots $ if $f^k_\#(\sigma) = \ldots
f^k_\#(\sigma_1) f^k_\#(\sigma_2)\ldots $ for all $k \ge 0$.  In other
words, a decomposition of $\sigma$ into subpaths $\sigma_i$ is a
splitting if one can tighten the image of $\sigma$ under any iterate
of $f_\#$ by tightening the images of the $\sigma_i$'s.

A path $\sigma$ is a {\em periodic Nielsen path} if $f^k_\#(\sigma) =
\sigma$ for some $k \ge 1$.  The minimal such $k$ is the {\em period}
of $\sigma$ and if $k = 1$ then $\sigma$ is a {\em Nielsen path}.  A
(periodic) Nielsen path is {\em indivisible} if it does not decompose
as a concatenation of non-trivial (periodic) Nielsen subpaths. A path
is {\em primitive} if it is not multiple of a simpler path.

\paragraph{Automorphisms and Lifts}  \label{autos and lifts}
Section 1 of \cite{gjll} and section 2.1 of \cite{bfh:tits3} are good
sources for facts that we record below without specific references.
The universal cover $\Gamma$ of a marked graph $G$ with marking $m
:R_n \to G$, is a simplicial tree.  We always assume that a base point
$\ti b \in \Gamma$ projecting to $b = m(*) \in G$ has been chosen,
thereby defining an action of $F_n$ on $\Gamma$.  The set of ends
$\E(\Gamma)$ of $\Gamma$ is naturally identified with the boundary
$\partial F_n$ of $F_n$ and we make implicit use of this
identification throughout the paper.

Each $c \in F_n$ acts by a {\em covering translation} $T_c : \Gamma
\to \Gamma$ and each $T_c$ induces a homeomorphism $\hat T_c :
\partial F_n \to \partial F_n$ that fixes two points, a sink $T_c^+$
and a source $T_c^-$.  The line in $\Gamma$ whose ends converge to
$T_c^-$ and $T_c^+$ is called the {\em axis of $T_c$} and is denoted
$A_c$.  The image of $A_c$ in $G$ is the circuit corresponding to the
conjugacy class of $c$.

If $\fG$ represents $\phi \in \Out(F_n)$ then the choice of a path
$\sigma$ from $b$ to $f(b)$ determines both an automorphism
representing $\phi$ and a lift of $f$ to $\Gamma$.  This defines a
bijection between the set of lifts $\ti f : \Gamma \to \Gamma$ of
$\fG$ and the set of automorphisms $\Phi : F_n \to F_n$ representing
$\phi$.  Equivalently, this bijection is defined by $\ti f T_c =
T_{\Phi(c)} \ti f$ for all $c \in F_n$.  We say that $\ti f$ {\em
corresponds to} $\aone$ or {\em is determined by} $\aone$ and vice
versa.  Under the identification of $\E (\Gamma)$ with $\partial F_n$,
a lift $\ti f$ determines a homeomorphism $\hat f$ of $\partial F_n$.
An automorphism $\Phi$ also determines a homeomorphism $\hat \Phi$ of
$\partial F_n$ and $\hat f = \hat \aone$ if and only if $\ti f$
corresponds to $\aone$.  In particular, $\hat i_c = \hat T_c$ for all
$c \in F_n$ where $i_c(w) = cwc^{-1}$ is the inner automorphism of
$F_n$ determined by $c$.  We use the notation $\hat f$ and $\hat \Phi$
interchangeably depending on the context.

We are particularly interested in the dynamics of $\hat f = \hat
\aone$. The following two lemmas are contained in Lemma 2.3 and Lemma
2.4 of \cite{bfh:tits3} and in Proposition 1.1 of \cite{gjll}.

\begin{lemma} \label{l: first from bk3} Assume that   $\ti f : \Gamma \to 
\Gamma$ corresponds to $\Phi \in \Aut(F_n)$.  Then the following are
equivalent:
\begin{description}
\item[(i)] $c \in \Fix(\Phi)$.
\item[(ii)] $T_c$ commutes with $\ti f$.
\item[(iii)] $\hat T_c$ commutes with $\hat f$.
\item[(iv)]  $\Fix(\hat T_c) \subset \Fix(\hat f) =  \Fix(\hat \Phi)$.
\item[(v)] $\Fix(\hat f) =  \Fix(\hat \Phi)$ is $\hat T_c$-invariant.
\end{description}
\end{lemma}

\vskip .3in

\begin{lemma} \label{l: second from bk3}  Assume that   $\ti f : \Gamma \to 
\Gamma$ corresponds to $\Phi \in \Aut(F_n)$ and that $\Fix(\hat \Phi)
\subset \partial F_n$ contains at least three points.  Denote
$\Fix(\Phi)$ by $\mathbb F$. Then
\begin{description}
\item[(i)] $\partial \mathbb F$ is naturally identified with the
closure of $\{T_c^{\pm} : T_c \in \T(\Phi)\}$ in $\partial F_n$. None
of these points is isolated in $\Fix(\hat \Phi)$.
\item[(ii)] Each point in $\Fix(\hat \Phi) \setminus \partial \mathbb
F$ is isolated and is either an attractor or a repeller for the action
of $\hat \Phi$.
\item[(iii)] There are only finitely many $\mathbb F$-orbits in
$\Fix(\hat \Phi) \setminus \partial \mathbb F$.
\end{description}
\end{lemma} 

\paragraph{Lines and Laminations}   \label{laminations}   
Unoriented bi-infinite paths in $G$ or its universal cover $\Gamma$
are called {\em lines}.  There is a bijection between lines in
$\Gamma$ and unordered pairs of distinct elements of $\partial F_n$,
the latter being the endpoints of the former.  The advantage of the
$\partial F_n$ description is that it allows us to work with {\em
abstract lines} that are {\em realized} in any $\Gamma$ but are not
tied to any particular $\Gamma$.

A closed set of lines in $G$ or an equivariant closed set of lines in
$\Gamma$ is called a {\em lamination} and the lines that compose it
are called {\em leaves}.  If $\Lambda$ is a lamination in $G$ then we
denote its pre-image in $\Gamma$ by $\ti \Lambda$ and vice-versa.
 
Suppose that $\fG$ represents $\oone$ and that $\ti f$ is a lift of
$f$.  If $\ti \gamma$ is a line in $\Gamma$ with endpoints $P$ and
$Q$, then there is a bounded homotopy from $\ti f(\ti \gamma)$ to the
line $\ti f_\#(\gamma)$ with endpoints $\hat f(P)$ and $\hat f(Q)$.
This defines an action $\ti f_\#$ of $\ti f$ on lines in $\Gamma$.  If
$\aone \in \Aut(F_n)$ corresponds to $\ti f$ then $\Phi_\# = \ti f_\#$
is described on abstract lines by $(P,Q) \mapsto (\hat \aone(P), \hat
\aone(Q))$.  There is an induced action $\oone_\#$ of $\oone$ on lines
in $G$ and in particular on laminations in $G$.  The stabilizer
$\Stab(\Lambda)$ of a lamination $\Lambda$ is the subgroup of elements
of $\Out(F_n)$ that preserve $\Lambda$.

A point $P \in \partial F_n$ determines a lamination $\Lambda(P)$,
called {\em the accumulation set of $P$}, as follows.  Let $\Gamma$ be
the universal cover of a marked graph $G$ and let $\ti R$ be any ray
in $\Gamma$ converging to $P$.  A line $\ti \sigma \subset \Gamma$
belongs to $\widetilde{ \Lambda(P)}$ if every finite subpath of $\ti
\sigma$ is contained in some translate of $\ti R$.  Since any two rays
converging to $P$ have a common infinite end, this definition is
independent of the choice of $\ti R$.  The bounded cancellation lemma
implies (cf.~Lemma 3.1.4 of \cite{bfh:tits1}) that this definition is
independent of the choice of $G$ and $\Gamma$ and that $\hat
\aone_\#(\widetilde{\Lambda(P)}) = \widetilde{\Lambda(\hat
\aone(P))}$.  In particular, if $P \in \Fix(\hat \aone)$ then
$\Lambda(P)$ is $\oone_\#$-invariant.

For each $\oone \in \Out(F_n)$ there is an associated
$\oone$-invariant finite set $\L(\oone)$ of laminations called the set
of {\em attracting laminations for $\oone$}.  For each $\Lambda \in
\L(\oone)$ there is an {\em expansion factor homomorphism}
$\PF_{\Lambda}$ defined on $\Stab(\Lambda)$ and with image a discrete
subgroup of $\R$.  Each $\Lambda \in \L(\oone)$ has birecurrent leaves
called {\em generic leaves}.  See section~3.3 of \cite{bfh:tits1} for complete
details.

\paragraph{The Recognition Theorem (\cite{fh:recognition})} 

The set of non-repelling fixed [resp. periodic] points of $\hat \aone$
is denoted by $\Fix_N(\hat \aone)$ [resp. $\Per_N(\hat \aone)$].

\begin{definition}
If the cardinality of $\Fix_N(\hat \aone)$ is greater than two or if
$\Fix_N(\hat \aone)$ is a pair of points that does not cobound either
some axis $T_c$ or a generic leaf of an element of $\L(\oone)$ then
$\aone$ is a {\em principal automorphism} and we write $\aone \in
\PA(\oone)$.  The corresponding lift of $f$ is called a {\em principal
lift}.  
\end{definition}

There is a natural equivalence relation on automorphisms defined by
$\aone_1 \sim i_c\aone_2 i_c^{-1}$ for some $c \in F_n$.  There are
only finitely many such equivalence classes of principal automorphisms
- see Remark 3.9 of \cite{fh:recognition}.

\begin{definition} 
An outer automorphism $\oone$ is {\em forward rotationless} if $\Fix_N(\hat
\aone) = \Per_N(\hat \aone)$ for all $\aone \in \PA(\oone)$ and if for
each $k \ge 1$, $\Phi \mapsto \Phi^k$ defines a bijection between
$\PA(\oone)$ and $\PA(\oone^k)$.  Our standing assumption is that $n
\ge 2$.  For notational convenience we say that the identity element
of $\Out(F_1)$ is forward rotationless.
\end{definition} 

Every $\oone \in \Out(F_n)$ has a forward rotationless iterate $\oone^k$ by
Corollary~3.30 of \cite{fh:recognition}.   As an illustration of the utility of this property, and for convenient reference, we recall
Lemma~3.29 of \cite{fh:recognition}.

\begin{lemma} \label{no iterates necessary} The following properties hold for each forward rotationless $\oone \in \Out(F_n)$. 
\begin{enumerate}
\item Each periodic conjugacy class is fixed and each representative
of that conjugacy class is fixed by some principal automorphism
representing $\phi$.
\item Each $\Lambda \in \L(\oone)$ is $\oone$-invariant.
\item A free factor that is invariant under an iterate of $\oone$ is
$\oone$-invariant.
\end{enumerate}
\end{lemma}

Several of our constructions are motivated by the following theorem
from \cite{fh:recognition}.  We also use this theorem directly to
prove that $\A(\oone)$ is abelian.

\begin{thm}{\bf (Recognition Theorem)} \label{t:recognition} Suppose that 
$\oone,\otwo\in\Out(\f)$ are forward rotationless and that
\begin{enumerate}
\item  $\PF_{\Lambda}(\oone)=\PF_{\Lambda}(\otwo)$, for all $\Lambda\in 
\L(\oone)=\L(\otwo)$.
\item there is bijection $h : \PA(\oone) \to \PA(\otwo)$ such that:  
\begin{itemize} 
\item [(i)] {\bf (fixed sets preserved)} $\Fix_N(\hat \aone) =
\Fix_N(\widehat{h(\aone)})$
\item [(ii)] {\bf (twist coordinates preserved)} If $ w \in
\Fix(\aone)$ is primitive and $\aone, i_{w^d}\aone \in \PA(\oone)$,
then $h(i_{w^d}\aone)= i_{w^d}h(\aone)$.
\end{itemize}
\end{enumerate}
\end{thm}

\begin{remark}
In the special case that $\oone$ is realized as an element of
$\mcg(S)$, a $w$ that occurs 
in item 2-(ii) represents a reducing curve and $d$ is the degree of
Dehn twisting about that reducing curve.  See also the discussion of
\lq axes\rq\ at the end of this section.
\end{remark} 

We include the following  result for easy reference.

\begin{lemma} \label{conjugation basics}
The following properties hold for all $\aone$ representing $\oone$ and
$\atwo$ representing $\otwo$.
\begin{enumerate}
\item $\Fix(\widehat{\atwo \aone \atwo^{-1}}) = \hat \atwo(\Fix(\hat
\aone))$ and $\Fix_N(\widehat{\atwo \aone \atwo^{-1}}) = \hat
\atwo(\Fix_N(\hat \aone))$.
\item Conjugation by $\atwo$ defines a bijection $i_{\atwo} :
\PA(\oone) \mapsto \PA (\otwo \oone \otwo^{-1})$ that preserves
equivalence classes.  The induced bijection on the set of equivalence
classes depends only on $\otwo$ and not on the choice of $\atwo$.
\end{enumerate}
\end{lemma}

\proof (1) is standard and easily checked; it implies that $i_{\atwo}
: \PA(\oone) \mapsto \PA (\otwo \oone \otwo^{-1})$ is a bijection.
The rest of (2) follows from $\atwo (i_c \aone i_c^{-1})\atwo^{-1} =
i_{\atwo(c)} \atwo \aone \atwo^{-1} i_{\atwo(c)}^{-1}$ and $(i_d\atwo
) (\aone)(i_d \atwo)^{-1} = i_d( \atwo \aone\atwo^{-1})i_d^{-1}$.
\endproof

\paragraph{Free Factor Systems} 
The conjugacy class of a free factor $F^i$ of $F_n$ is denoted
$[[F^i]]$.  If $F^1,\ldots,F^k$ are non-trivial free factors and if
$F^1 \ast \ldots \ast F^k$ is a free factor then we say that the
collection $\{[[F^1]],\ldots,[[F^k]]\}$ is a {\em free factor system}.
For example, if $G$ is a marked graph and $G_r \subset G$ is a
subgraph with non-contractible components $C_1,\ldots, C_k$ then the
conjugacy class $[[\pi_1(C_i)]]$ of the fundamental group of $C_i$ is
well defined and the collection of these conjugacy classes is a free
factor system denoted $\F(G_r)$.

The image of a free factor $F$ under an element of $\Aut(F_n)$ is a
free factor.  This induces an action of $\Out(F_n)$ on the set of free
free systems.  We sometimes say that a free factor is
$\oone$-invariant when we really mean that its conjugacy class is
$\oone$-invariant.  If $[[F]]$ is $\oone$-invariant then $F$ is
$\aone$-invariant for some automorphism $\aone$ representing $\oone$
and $\aone|F$ determines a well defined element $\oone|F$ of
$\Out(F)$.

The conjugacy class $[a]$ of $a \in F_n$ is {\em carried by $[[F^i]]$}
if $F^i$ contains a representative of $[a]$. Sometimes we say that $a$
is carried by $F^i$ when we really mean that $[a]$ is carried by
$[[F^i]]$.  If $G$ is a marked graph and $G_r$ is a subgraph of $G$
such that $[[F^i]] = \F(G_r)$ then $[a]$ is carried by $[[F^i]]$ if
and only if the circuit in $G$ that represents $[a]$ is contained in
$G_r$.  We say that {\em an abstract line $ \ell$ is carried by
$[[F^i]]$} if its realization in $G$ is contained in $G_r$ for some,
and hence any, $G$ and $G_r$ as above.  Equivalently, $\ell$ is the
limit of periodic lines corresponding to $[c_i]$ where each $[c_i]$ is
carried by $[[F^i]]$.  A collection $W$ of abstract lines and
conjugacy classes in $F_n$ is carried by a free factor system $\F
=\{[[F^1]],\ldots,[[F^k]]\}$ if each element of $W$ is carried by some
$F^i$.

There is a partial order $\sqsubset$ on free factor systems generated
by inclusion. More precisely, $[[F^1]] \sqsubset [[F^2]]$ if $F^1$ is
conjugate to a free factor of $F^2$ and $\F_1 \sqsubset \F_2$ if for
each $[[F^i]] \in \F_1$ there exists $[[F^j]] \in \F_2$ such that
$[[F^i]] \sqsubset [[F^j]]$.

The {\em complexity} of a \ffs\ is defined on page 531 of
\cite{bfh:tits1}.  We include the following results for easy
reference.  The first is Corollary~2.6.5. of \cite{bfh:tits1}.  The
second is an immediate consequence of the uniqueness of $\F(W)$.

\begin{lemma}\label{ffs lemma}
For any collection $W$ of abstract lines there is a unique free factor
system $\F(W)$ of minimal complexity that carries every element of
$W$.  If $W$ is a single element then $\F(W)$ has a single element.
\end{lemma}

\begin{cor} \label{ffs invariance}
If a collection $W$ of abstract lines and conjugacy classes in $F_n$
is $\oone$-invariant then $\F(W)$ is $\oone$-invariant.
\end{cor} 

Further details on free factor systems can be found in section~2.6 of
\cite{bfh:tits1}.

\paragraph{Relative Train Track Maps} \label{ss:cs} 
We assume some familiarity with the basic definitions of relative
train track maps. Complete details can be found in
\cite{fh:recognition} and \cite{bfh:tits1}.

Suppose that $\fG$ is a \rtt\ defined with respect to a maximal
filtration $\filt$.  A path or circuit has {\em height} $r$ if it is
contained in $G_r$ but not $G_{r-1}$.  A lamination has height $r$ if
each leaf in its realization in $G$ has height at most $r$ and some
leaf has height $r$.  The $r^{th}$ stratum $H_r$ is defined to be the
closure of $G_r \setminus G_{r-1}$.  If $f(H_r) \subset G_{r-1}$ then
$H_r$ is called a {\em zero stratum}; all other strata have
irreducible transition matrices and are said to be {\em irreducible}.
If $H_r$ is irreducible and if the Perron-Frobenius eigenvalue of the
transition matrix for $H_r$ is greater than one, then $H_r$ is {\em
exponentially growing} or simply EG.  All other irreducible strata are
{\em non-EG} or simply {\em NEG}.

A {\em direction} $d$ at $x \in G$ is the germ of an initial segment
of an oriented edge (or partial edge if $x$ is not a vertex) based at
$x$.  There is an $f$-induced map $Df$ on directions and we say that
$d$ is {\em a periodic direction} if it is periodic under the action
of $Df$; if the period is one then $d$ is a {\em fixed direction}.
Thus the direction determined by an oriented edge $E$ is fixed if and
only if $E$ is the initial edge of $f(E)$.

A \emph{turn} is an unordered pair of directions with a common base
point.  The turn is \emph{nondegenerate} if is defined by distinct
directions and is \emph{degenerate} otherwise.  If $E_1E_2 \ldots
E_{k-1} E_k$ is the edge path associated to a path $\sigma$, then we
say that {\it $\sigma$ contains the turns $(\bar E_i,E_{i+1})$ for $1
\le i \le k-1$}.  A turn is {\it illegal} with respect to $\fG$ if its
image under some iterate of $Df$ is degenerate; a {turn is \it legal}
if it is not illegal.  A {\it path or circuit $\sigma \subset G$ is
legal} if it contains only legal turns.  If $\sigma \subset G_r$ does
not contain any illegal turns in $H_r$, meaning that both directions
correspond to edges of $H_r$, then $\sigma$ is $r$-legal.  It is
immediate from the definitions that $Df$ maps legal turns to legal
turns and that the restriction of $f$ to a legal path is an immersion.

A non-trivial path in a zero stratum $H_i$ whose endpoints belong to
EG strata is called a {\it connecting path}.

Suppose that $H_i$ and $H_j$ are distinct NEG strata consisting of
single edges $E_i$ and $E_j$, that $w$ is a primitive Nielsen path and
that $f(E_i) = E_iw^{d_i}$ and $f(E_j) = E_jw^{d_j}$ for some $d_i,d_j
> 0$.  Then a path of the form $E_i w^p \bar E_j$ is called an {\em
exceptional path of height $\max(i,j)$}  or just an {\em
exceptional path} if the height is not relevant.  The set of
exceptional paths of height $i$ is invariant under the action of
$f_\#$.

\begin{definition}
A non-trivial path or circuit $\sigma$ is {\em completely split}
if it has a splitting, called a {\it complete splitting}, into
subpaths, each of which is either a single edge in an irreducible
stratum, an \iNp, an exceptional path or a connecting path that is
maximal in the sense that it cannot be extended to a larger connecting
path in $\sigma$.
\end{definition}

\begin{definition}
A \rtt\ is {\em \cs}\ if
\begin{enumerate}
\item $f(E)$ is completely split for each edge $E$ in an irreducible
stratum.
\item $f|\sigma$ is an immersion and $f(\sigma)$ is completely split
for each connecting path $\sigma$.
\end{enumerate}
\end{definition}

\begin{remark}  \label{cs is unique}
For $\fG$ satisfying the conclusions of Theorem~\ref{comp sp exists}
below, a completely split path or circuit has a unique complete
splitting by Lemma~4.14 of \cite{fh:recognition}.
\end{remark}

\begin{definition}\label{rotRttDef}
A periodic vertex $w$ that does not satisfy one of the following two
conditions is {\em principal}.
\begin{itemize}
\item $w$ is the only element of $\Per(f)$ in its Nielsen class and
there are exactly two periodic directions at $w$, both of which are
contained in the same EG stratum.
\item $w$ is contained in a component $C$ of $\Per(f)$ that is
topologically a circle and each point in $C$ has exactly two periodic
directions.
\end{itemize}   
We also say that a lift of a principal vertex to the universal cover
is a principal vertex.
\end{definition}

\begin{remark} \label{some essential}
It is immediate from the definition that the initial endpoint of an
NEG edge is a principal vertex.  By Lemma~3.18 of
\cite{fh:recognition} every EG stratum $H_r$ contains a principal
vertex that is the basepoint for a periodic direction in $H_r$.
\end{remark}

\begin{definition}
If the endpoints of all indivisible periodic Nielsen paths are
vertices and if each principal vertex and each periodic direction at a
principal vertex has period one then we say that $\fG$ is {\em forward
rotationless}.
\end{definition}  

  Proposition~3.28 of \cite{fh:recognition} states that for \rtt s
satisfying certain mild assumptions the two notions of forward
rotationless coincide.  Namely $\fG$ is forward rotationless if and
only if $\oone$ is forward rotationless.  Assuming that $\fG$ is
forward rotationless, Corollary~3.21 and Lemma~3.26 of
\cite{fh:recognition} imply that $\ti f$ is a principal lift if and
only if some (and hence every) element of $\Fix(\ti f)$ is a principal
vertex.

A vertex in $G$ is an {\em attaching vertex} if it belongs to a
non-contractible component of $G_r$ and is the endpoint of an edge in
$H_s$ for $s > r$.  We recall Theorem~4.6 of \cite{fh:recognition}.

\begin{thm} \label{comp sp exists}
Every forward rotationless $\oone \in \Out(F_n)$ is represented by a
forward rotationless \cs\ \rtt\ $\fG$.  If $\F$ is a $\oone$-invariant \ffs,
then one may choose the associated $f$-invariant filtration $\filt$ so
that $\F = \F(G_l)$ for some filtration element $G_l$. Moreover,
\begin{description}
\item [(V)] Each attaching vertex is principal (and hence fixed).
\item [(NEG)] Each non-fixed \noneg\ stratum $H_r$ is a single edge $E_r$ oriented so that 
  $f(E_r) = E_r\cdot u_r$ for some nontrivial closed path $u_r \subset
G_{r-1}$.  The initial vertex of $E_r$ is principal. 
\end{description}
If $u_r$ is a non-trivial Nielsen path then $H_r$ and $E_r$ are said
to be {\em linear} and we sometimes write $u_r = w_r^{d_r}$ where
$w_r$ is primitive.
\begin{description}
\item [(L)] If $E_r$ and $E_s$ are distinct linear edges and if $w_r$
and $w_s$ agree as unoriented loops, then $w_r = w_s$ and $d_r \ne
d_s$.
\item[(N)] Every periodic Nielsen path has period one.  The endpoints
of each \iNp\ $\sigma$ are vertices.  For each EG stratum $H_r$ there is at most one $\sigma$ of height $r$.  If $\sigma$ has height $r$ and
if $H_r$ is not \EG\ then $H_r$ is linear and $\sigma = E_r w_r^k \bar
E_r$ for some $k \ne 0$.
\item [(Per)] The vertices in any non-trivial component $C$ of
$\Per(f)$ are principal.  In particular $C \subset \Fix(f)$.  If $C$
is contractible and contains an edge in $H_r$, $ 1 \le r \le N$, then
some vertex of $C$ has valence at least two in $G_{r-1}$.
\item[(Z)] $H_i$ is a zero stratum if and only if it is a contractible
component of $G_i$.  If $H_j$ is the first irreducible stratum
following a zero stratum then $H_j$ is \eg\ and all components of
$G_j$ are non-contractible.   If
$H_i$ is a zero stratum then $f|H_i$ is an immersion. If the link of a vertex $v$ is contained in a zero stratum $H_i$ then $v$ has valence at least three in $H_i$.
\end{description}
\end{thm}

{\em We assume throughout the remainder of this paper that $\fG$
satisfies the conclusions of Theorem~\ref{comp sp exists}.}

\paragraph{Iterating an Edge}
We make frequent use of isolated points in $\Fix_N(\hat f)$ for
principal lifts $\ti f$.  We quote two results that we refer to
several times for the reader's convenience.  The first is a
combination of Lemma~3.25 and Lemma~4.19 of \cite{fh:recognition}. The
second is Lemma~4.21 of that same paper.

\begin{lemma} \label{iterates to}
Assume that $\fG$ satisfies the conclusions of Theorem~\ref{comp sp
exists}. The following properties hold for every principal lift $\ti f
: \Gamma \to \Gamma$.
\begin{enumerate}
\item If $\ti v \in \Fix(\ti f)$ and a non-fixed edge $\ti E$
determines a fixed direction at $\ti v$, then $\ti E\ \subset \ti
f_\#(\ti E) \subset \ti f_\#^2(\ti E) \subset \ldots$ is an increasing
sequence of paths whose union is a ray $\ti R$ that converges to some
$P \in \Fix_N(\hat f)$ and whose interior is fixed point free.  If
$\ti E$ is a lift of an edge in an EG stratum then accumulation set of
$P$ is the element in $\L(\oone)$ corresponding to that stratum.
\item For every isolated $P \in \Fix_N(\hat f)$ there exists $\ti E$
and $\ti R$ as in (1) that converges to $P$.
\end{enumerate}
\end{lemma}

If $\ti E$ and $P$ are as in Lemma~\ref{iterates to} then we say that
$\ti E$ {\em iterates to} $P$ and that $P$ is {\em associated to }
$\ti E$.

\begin{lemma} \label{isolated for egs}
Suppose that $\otwo \in \ofn$ is forward rotationless and that $P \in
\Fix_N(\hat \atwo)$ for some $\atwo \in \PA(\otwo)$.  Suppose further
that $\Lambda$ is an attracting lamination for some element of $\ofn$,
that $\Lambda$ is $\otwo$-invariant and that $\Lambda$ is contained in
the accumulation set of $P$.  Then $\PF_{\Lambda}(\otwo) \ge 0$ and
$\PF_{\Lambda}(\otwo) > 0$ if and only if $P$ is isolated in
$\Fix_N(\hat \atwo)$.
\end{lemma}

\paragraph{Axes}
Assume that $\oone$ is forward rotationless and that $\fG$ is as in
Theorem~\ref{comp sp exists}.  Following the notation of
\cite{bfh:tits3} we say that an unoriented conjugacy class $\mu$ of a
primitive element of $F_n$ is an {\em axis for $\oone$} if for some
(and hence any) representative $c \in F_n$ there exist distinct
$\aone_1, \aone_2 \in \PA(\oone)$ that fix $c$.  Equivalently
$\Fix_N(\hat \aone_1) \cap \Fix_N(\hat \aone_2)$ is the endpoint set
of the axis $ A_c$ for $T_c$.  The number of distinct elements of
$\PA(\oone)$ that fix $c$ is called the {\em multiplicity} of $\mu$.
It is a consequence of Lemma~\ref{axes and lifts} below that both the
number of axes and the multiplicity of each axis is finite.

Lemma~\ref{no iterates necessary} implies that the oriented conjugacy
class of $c$ is $\oone$-invariant.  By Lemmas~4.1.4 and 4.2.6 of
\cite{bfh:tits1}, the circuit $\gamma$ representing $c$ splits into a
concatenation of periodic, and hence fixed, Nielsen paths.  There is
an induced decomposition of $A_c$ into subpaths $\ti \alpha_i$ that
project to either fixed edges or \iNp s.  The lift $\ti f_0 : \Gamma
\to \Gamma$ that fixes the endpoints of each $\ti \alpha_i$ is a
principal lift and commutes with $T_c$.  We say that $\ti f_0$ and the
corresponding $\aone_0 \in \PA(\oone)$ are the {\em base lift} and
{\em base principal automorphism} associated to $\mu$ and the choices
of $T_c$ and $\fG$.  (If $\mu$ is not represented by a basis element
then $\aone_0$ is independent of the choice of $\fG$ but otherwise it
is not; see Example~\ref{D depends on f} for ramifications of this
fact.)  Remark~\ref{cs is unique} implies that $\ti f_0$ is the only
lift that commutes with $T_c$ and has fixed points in $A_c$.

We recall Lemma~4.23 of \cite{fh:recognition}.

\begin{lemma}  \label{axes and lifts}
Assume notation as above and that $\fG$ satisfies the conclusions of
Theorem~\ref{comp sp exists}.  There is a bijection between principal
lifts [principal automorphisms] $\ti f_j \ne \ti f_0$ [respectively
$\aone_j \ne \aone_0 \in \PA(\oone)$] that commute with $T_c$ [fix
$c$] and the linear edges $\{E_j\}$ with $w_j$ representing $\mu$.
Moreover, if $f(E_j) = E_j w_j^{d_j}$ then $\ti f_j = T_c^{d_j} \ti
f_0$ [$\aone_j = i_c^{d_j}\aone_0$].
\end{lemma}

\section{Rotationless Abelian Subgroups} \label{s:lifts}
The Recognition Theorem is stated purely in terms of $\oone$ and its
forward iterates.  No condition on $\oone^{-1}$ is required.  In the
context of abelian subgroups, it is more natural to give $\oone$ and
$\oone^{-1}$ equal footing.
 
\begin{definition} \label{d:birot} $\PA^\pm(\oone)=\PA(\oone)\cup\PA(\oone^{-1})^{-1}$.
An outer automorphism $\oone$ is
{\em rotationless} if $\Fix(\hat \aone) = \Per(\hat \aone)$ for all
$\aone \in \PA^{\pm}(\oone)$ and if for each $k \ge 1$, $\Phi \mapsto
\Phi^k$ defines a bijection (see Remark~\ref{surjection}) between
$\PA^{\pm}(\oone)$ and $\PA^{\pm}(\oone^k)$.  A subgroup of
$\Out(F_n)$ is {\em rotationless} if each of its elements
is.  
\end{definition}

\begin{remark} \label{surjection}
There is no loss in replacing the assumption that $\Phi \mapsto
\Phi^k$ defines a bijection with the a priori weaker assumption that
$\Phi \mapsto \Phi^k$ defines a surjection.  Indeed if $\Phi \mapsto
\Phi^k$ is not injective then there exist distinct $\aone_1,\aone_2
\in \PA^{\pm}(\oone)$ and $k \ge 1$ such that $\aone_1^k = \aone_2^k$.
This contradicts the fact that $\aone_2 \aone_1^{-1}$ is a non-trivial
covering translation and the fact that $\Fix(\hat \Phi_2 \hat
\Phi_1^{-1})$ contains $\Fix(\hat \aone_1) = \Fix(\hat \aone_1^k) =
\Fix(\hat \aone_2)$ and so contains at least three points.
\end{remark}

The natural guess is that $\oone$ is rotationless if and only if
$\oone$ and $\oone^{-1}$ are forward rotationless.  The following lemma and
corollary fall short of proving this but is sufficient for our needs.

\begin{lemma}  \label{l:rotationless} 
\begin{enumerate}
\item If $\oone$ is rotationless then $\oone$ and $\oone^{-1}$ are
forward rotationless.
\item If $\oone$ and $\oone^{-1}$ are forward rotationless and $(\ast)$ is
satisfied for $\theta = \oone$ and $\theta = \oone^{-1}$ then $\oone$
is rotationless.
\begin{itemize}
\item [$(\ast)$] For all $\Theta \in \PA(\theta)$, the set of
repelling periodic points for $\hat \Theta$ is not a period two orbit
that is the endpoint set of a lift of a generic leaf $\gamma$ of an
element of $\L(\theta^{-1})$ .
\end{itemize}
\end{enumerate}
\end{lemma}

\proof Assume that $\oone$ is rotationless.  For $k > 0$, each
element of $\PA(\oone^k)$ has the form $\aone^k$ where $\Fix(\hat
\aone) = \Per(\hat \aone)$ and hence $\Fix_N(\hat \aone) = \Per_N(\hat
\aone^k)$.  Thus $\aone \in \PA(\oone)$ proving that $\oone$ is
forward rotationless.  The symmetric argument showing that $\oone^{-1}$ is
forward rotationless completes the proof of (1).

Assume now that the hypotheses of (2) are satisfied, that $k \ge 1$
and that $\aone_k \in \PA^{\pm}(\oone^k)$.  The plus and minus cases
are symmetric so we may assume that $\aone_k \in \PA(\oone^k)$.  Since
$\oone$ is forward rotationless, $\aone_k = \aone^k$ for some $\aone \in
\PA(\oone)$ satisfying $\Fix_N(\hat \aone) = \Per_N(\hat \aone)$.  To
prove that $\Fix(\hat \aone) = \Per(\hat \aone)$ it suffices to show
that all periodic repelling points for $\hat \aone$ have period one.
Since $\oone^{-1}$ is forward rotationless, the only way this could fail would
be if the repelling set is a period two orbit and if $\aone^2 \not \in
\PA(\oone^{-1})$.  This possibility is ruled out by $(\ast)$.
\endproof

\begin{corollary}  
If $\oone$ and $\oone^{-1}$ are forward rotationless then $\oone^2$ is
 rotationless.  There exists $k > 0$ so that $\oone^{2k}$ is
 rotationless for every $\oone \in \Out(F_n)$.
\end{corollary}

\proof The first statement follows from Lemma~\ref{l:rotationless}.
The second follows from the first and from the fact (Corollary~4.26 of \cite{fh:recognition}) that there is
a uniform $k > 0$ such that $\oone^k$ are forward rotationless for all
$\oone$.  \endproof

\begin{ex}  Let $G$ be the graph with one vertex $v$ and edges labelled $A$, $B$ and $C$.  Let $f: G \to G$ be the homotopy equivalence defined by 
$$
A \mapsto B^3 A \qquad \qquad B \mapsto C^3B \qquad \qquad C \mapsto (B^3A)^3 C.
$$ 

The directions at $v$ determined by $\bar A, \bar B$ and $\bar C$ are
fixed by $Df$ and those determined by $B$ and $C$ are interchanged by
$Df$.  Thus $f$ is not rotationless and the outer automorphism $\phi$
that it determined is neither forward rotationless nor rotationless.
The map $f$ factors as $f_3f_2f_1$ where $f_1$ fixes $A$ and $B$ and
$f_1(C) = A^3C$, $f_2$ fixes $B$ and $C$ and $f_2(B)=C^3B$ and $f_3$
fixes $A$ and $B$ and $f_3(A) = B^3A$.  It is easy to check that each
of these homotopy equivalence determines a rotationless element of
$\Out(F_n)$. This shows that the composition of rotationless elements
need not be rotationless.  Obviously, $\phi$ induces the identity on
$H_1(G,\Z_3)$ and so illustrates that not every such element is
rotationless.
\end{ex} 

\begin{lemma} \label{not empty}
If $\oone$ is   rotationless and $\aone
\in \PA(\oone)$ then $\Fix_N(\hat \aone^{-1}) \ne \emptyset$.  
\end{lemma} 

\proof Choose $\fG$ representing $\oone^{-1}$ and let $\ti f : \Gamma
\to \Gamma$ be the lift corresponding to $\aone^{-1}$.  It suffices to
show that $\Fix_N(\hat f^k) \ne \emptyset$ for some $k \ge 1$.  This
follows from Lemma~3.23 of \cite{fh:recognition} if $\Fix(\ti f) =
\emptyset$ and from Lemma~3.26 of \cite{fh:recognition} otherwise.
\endproof

Abelian subgroups of $\Out(F_n)$ are finitely generated
\cite{bl:niltech}.  Thus given any generating set for an abelian
subgroup, there is a finite subset which also generates.  At the end
of this section (Corollary~\ref{is birot}) we prove that an abelian
subgroup $A$ of $\Out(F_n)$ that is generated by rotationless
elements, is rotationless.

Many of our arguments proceed by induction on the cardinality  
of a given set of rotationless generators.

\begin{lemma}
If $\oone$ is rotationless and $F$ is a $\oone$-invariant free
factor of rank at least two then $\oone|F$ is rotationless.
\end{lemma}

\proof This follows from the definitions and the fact that an element
of $\PA^{\pm}(\oone|F)$ extends to an element of $\PA^{\pm}(\oone)$.
\endproof

We produce lifts of an abelian subgroup of $\Out(F_n)$ to $\Aut(F_n)$
that is generated by rotationless elements via the following
definition and lemma.

\begin{definition} 
A set $X \subset \partial F_n$ with at least three points is {\em a
principal set} for a subgroup $A$ of $\Out(F_n)$ if each $\otwo \in A$
is represented by $\atwo \in \Aut(F_n)$ satisfying $X \subset
\Fix(\hat \atwo)$ and if this necessarily unique $\atwo$ is an element
of $\PA^{\pm}(\otwo)$.  The assignment $\otwo \mapsto \atwo$ is a lift
of $A$ from $\Out(F_n)$ to $\Aut(F_n)$.
\end{definition}

\begin{lemma} \label{principal sets}
Suppose that $A$ is an  abelian subgroup of $\Out(F_n)$ that is generated by rotationless elements, 
that $\oone \in A$ is
rotationless and that $\aone \in \PA^{\pm}(\oone)$.
\begin{enumerate}
\item If $\Fix(\aone)$ has rank zero then $\Fix(\hat \aone)$ is a
principal set for $A$.
\item If $\Fix(\aone)$ has rank one with generator $c$ and if $P$ is
an isolated point in $\Fix(\hat \aone)$ then $\{P, T_c^{\pm}\}$ is a
principal set for $A$.
\item If $\Fix(\aone)$ has rank at least two then $\partial
\Fix(\aone)$ contains at least one principal set $X$ for $A$ and one
can choose $X$ to contain $T_c^{\pm}$ for any given $A$-invariant
$[c]$ for $c \in \Fix(\aone)$.  Moreover, for every isolated
point $P$ in $\Fix(\hat \aone)$ there is a principal set $Y$ for $A$
that contains $P$ and at least two elements of $\partial \Fix(\aone)$.
\end{enumerate}
In particular, $\Fix(\hat \aone)$ contains at least one principal set
for $A$ and every isolated point in $\Fix(\hat \aone)$ is contained
in a principal set.  If $s : A \to \Out(F_n)$ is the lift determined
by a principal set contained in $\Fix(\hat \aone)$ then $s(\oone) =
\aone$.
\end{lemma}

\proof Let $\otwo$ be an element of a a finite rotationless generating
set $S$ for $A$ and let $\atwo$ represent $\otwo$.
Lemma~\ref{conjugation basics} implies that conjugation by $\atwo$
defines a permutation of the finite set of equivalence classes in
$\PA^{\pm}(\oone)$.  Choose $k>0$ so that the permutation induced by
$\atwo^k$ is trivial.  Then $\atwo^k \aone {\atwo}^{-k} = i_c \aone
i_c^{-1}$ for some $c \in F_n$ and $\atwo_k := i_c^{-1} \atwo^k$
commutes with $\aone$.  In particular, $\mathbb F := \Fix(\aone)$ is
$\atwo_k$-invariant.

Assume at first that $\mathbb F$ has rank zero.  By Lemma~\ref{l:
second from bk3}, $\Fix(\hat \aone)$ is a finite union of attractors
and repellers and by Lemma~\ref{not empty} there is at least one of
each.  Since $\aone \in \PA^{\pm}(\oone)$, there are at least three
points in $\Fix(\hat \aone)$.

 We claim that if $\Theta$ represents $\theta \in A$ and if $\Fix(\hat
\aone)\subset \Fix(\hat \Theta)$ then $\Theta \in \PA^{\pm}(\theta)$.
If $\Fix(\hat \Theta)$ contains at least five points then this is
obvious.  After replacing $\theta$ with its inverse if necessary,
there are two potentially bad cases.  The first is that $\Fix(\hat
\Theta)$ has exactly one repelling point and exactly two attracting
points and that the attractors bound a lift $\ti \gamma$ of a generic
leaf of some $\Lambda \in \L(\theta)$.  Since the endpoints of $\ti
\gamma$ are isolated fixed points of $\hat \aone$, $\Lambda \in
\L(\oone) \cup \L(\oone^{-1})$ by Lemma~\ref{isolated for egs}.  After
replacing $\oone$ with its inverse if necessary, we may assume that
$\Lambda \in \L(\oone)$ and that the endpoints of $\ti \gamma$ are
attractors for $\aone$.  Since $\Fix(\hat \aone)$ contains only three
points and by Lemma~\ref{not empty} has at least one $\hat
\aone$-repeller, this contradicts the assumption that $\aone \in
\PA^{\pm}(\oone)$.

The other bad possibility is that $\Fix(\hat \Theta)$  is
a four point set with two repelling points that bound a lift of a leaf
of an element of $\L(\theta^{-1})$ and two attracting points that
bound a lift of a leaf of an element of $\L(\theta)$.  As in the
previous case, this description also applies to $\aone$ in
contradiction to the assumption that $\aone \in \PA^{\pm}(\oone)$.
This completes the proof that $\Theta \in \PA^{\pm}(\theta)$.

After replacing $\atwo_k$ with an iterate,  we may assume that $\Fix(\hat
\aone) \subset \Fix(\hat \atwo_k)$ and hence that $\atwo_k \in
\PA^{\pm}(\otwo^k)$.  Since $\otwo$ is rotationless, there exists
$\atwo \in \PA^{\pm}(\otwo)$ with $\Fix(\hat \aone) \subset
\Fix(\hat \atwo)$.  As this
holds for every element of $S$, we have proved (1).

Suppose next that $\mathbb F $ has rank one with generator $c$ and
that $P$ is an isolated point in $\Fix(\hat \aone)$.  Lemma~\ref{l:
second from bk3} implies that there are only finitely many
$i_c$-orbits of isolated points in $\Fix(\hat \aone)$.  After
increasing $k$ if necessary, we may assume that $c \in \Fix(\atwo_k)$
and that $\atwo_k$ preserves each such $i_c$-orbit.  In particular,
$\hat \atwo_k(P) = \hat T_c^q(P)$  for some $q$.  Let $\atwo_k':=
i_c^{-q}\atwo_k$.  Then $\{T_c^{\pm},P\} \subset \Fix(\hat \atwo_k')$
and $\atwo_k' \in \PA^{\pm}(\otwo)$.  Since $\otwo$ is rotationless,
there exists $\atwo \in \PA^{\pm}(\otwo)$ such that
$\{T_c^{\pm},P\}\subset \Fix(\hat\atwo)$.  As this holds for every
element of $S$, it follows that for each $\theta \in A$ there exists
$\Theta$ such that $\{T_c^{\pm},P\} \subset \Fix(\hat \Theta)$.  In
this case it is obvious that $\Theta \in \PA^{\pm}(\theta)$.  This
completes the proof of (2).

We turn next to the moreover part of (3).  Assume that $P$ is an
isolated point in $\Fix(\hat \aone)$.  As in the rank one case, the
fact that there are only finitely many $\mathbb F$-orbits of isolated
points in $\Fix(\aone)$ allows us to choose $\atwo_k^{\ast}$
representing an iterate $\otwo^k$ of $\otwo$ such that $P \in
\Fix(\widehat {\atwo_k^{\ast}})$ and such that $\mathbb F$ is
$\atwo_k^{\ast}$-invariant.  We claim that $\atwo^{\ast}_k \in
\PA^{\pm}(\otwo^k)$.  Assuming without loss that $\Fix(\widehat
{\atwo^{\ast}_k | \mathbb F})$ is finite, Lemma~\ref{not empty}
implies, after replacing $\atwo^{\ast}_k $ by an iterate if necessary,
that $\Fix(\widehat {\atwo^{\ast}_k | \mathbb F})$ has at least one
non-attractor $Q_-$ and one non-repeller $Q_+$.  Lemma~\ref{isolated
for egs} implies that $Q_+$ and $Q_-$ do not cobound a lift of a
generic leaf of an attracting lamination. (This method for proving
that a pair of points do not cobound a lift of a generic leaf of an
attracting lamination is used implicitly throughout the rest of the
proof.)  Generic leaves of an attracting lamination are birecurrent
and so either have both endpoints in $\partial \mathbb F$ or neither
endpoint in $\partial \mathbb F$.  Thus $P$ and $Q_{\pm}$ do not
cobound a lift of a generic leaf of an attracting lamination.  This
verifies our claim.  Since $\otwo$ is rotationless, there exists
$\atwo \in \PA^{\pm}(\otwo)$ with $\{P,Q_+,Q_-\} \subset \Fix(\hat
\atwo)$.  These three points are also in $\Fix(\hat \aone)$.  It
follows that $\atwo$ commutes with $\aone$ and hence that $\mathbb F$
is $\atwo$-invariant.

We have shown that if $S = \{\otwo_1,\ldots,\otwo_K\}$ then for all $1
\le j\le K$ there exists $\atwo_j \in \PA^{\pm}(\otwo_j$) such that $P
\in \Fix(\hat \atwo_j)$ and such that $\mathbb F$ is
$\atwo_j$-invariant.  Since $P$ is not fixed by any covering
translation, the $\atwo_j$'s commute.

We produce the desired principal set $Y$ by induction on $j$.  To this
end, let $Y_j = (\bigcap_{i=1}^j \Fix(\widehat{ \atwo_i)}) \cap
\partial \mathbb F $ and let $I_j$ be the statement that $Y_j$ either
contains three points or contains two points that do not cobound a
lift of a generic leaf of any attracting lamination.  As noted above,
$P$ and an element of $\partial \mathbb F$ can not cobound a generic
leaf of an attracting lamination. Thus $I_K$ completes the proof of
the moreover part of (3).

$I_1$ follows from Lemma~\ref{not empty} applied to $\atwo_1|\mathbb
F$.  Assume that $I_{j-1}$ holds.  $Y_{j-1}$ is $\hat
\atwo_j$-invariant.  If $Y_{j-1}$ is finite then it is fixed by an
iterate of $\hat \atwo_j$ and hence by $\hat \atwo_j$.  If $Y_{j-1}$
contains $T_b^{\pm}$ for some unique primitive unoriented $b$ then
$T_b^{\pm}$ is fixed by an iterate of $\hat \atwo_j$ and hence by
$\hat \atwo_j$.  In either case $I_j$ holds.  In the remaining case
$\cap_{i=1}^{j-1}\Fix(\atwo_j)$ intersects $\mathbb F$ in a subgroup
$\mathbb F_{j-1}$ of rank at least two and $I_j$ follows from
Lemma~\ref{not empty} applied $\hat \aone_j|\mathbb F_{j-1}$, keeping
in mind that $\Fix(\hat \aone_j|\mathbb F_{j-1}) \subset Y_j$. This
completes the induction step and so proves $I_K$.

It remains to prove the main statement of (3).  We argue by induction
on the cardinality $K$ of a given rotationless generating set for
$A$.  If $K = 1$ and $S=\{\otwo\}$ then there exists $\atwo \in
\PA^{\pm}(\otwo)$ such that $\Fix(\hat \atwo) = \Fix(\hat \aone)$ and
$\Fix(\hat \atwo)$ is obviously a principal set for $A$.  We now
assume that $K \ge 2$ and that (3) holds for subgroups that are
generated by fewer than $K$   rotationless
elements.

The defining property of $\atwo_k$ is that it commutes with $\aone$.
We may therefore replace our current $\atwo_k$ with any lift of any
iterate of $\otwo$ that preserves $\mathbb F$.  By Lemma 5.2 of
\cite{bfh:tits3} or Proposition I.5 of \cite{ll:ends}, there is such a
lift, still called $\atwo_k$, such that $\atwo_k |\mathbb F \in
\PA^{\pm}(\otwo^k|\mathbb F)$; moreover if $c \in \Fix(\aone)$ is
$A$-invariant then we may choose $\atwo_k$ so that $c \in
\Fix(\atwo_k)$.  Since $\otwo$ is rotationless, there exists $\atwo
\in \PA^{\pm}(\otwo)$ such that $\Fix(\hat \atwo) = \Fix(\hat
\atwo_k)$.  Thus $\Fix(\hat \atwo) \cap \Fix(\hat \aone)$ contains at
least three points which implies that $\atwo$ commutes with $\aone$.
To summarize, we have $\atwo \in \PA^{\pm}(\otwo)$ that preserves
$\mathbb F$ and such that $\atwo|\mathbb F \in \PA^{\pm}(\otwo|\mathbb
F)$; if $c \in \Fix(\hat \aone)$ is $A$-invariant then we may assume
that $c \in \Fix(\atwo)$.  As each $\atwo$ preserves $\mathbb F$, it
follows that $[\mathbb F]$ is $A$-invariant.

Let $A' = A|\mathbb F$, let $\otwo' = \otwo|\mathbb F$ and let $\atwo'
= \atwo|\mathbb F$. A principal set for $A'$ is also a principal set
for $A$ because an automorphism of $\mathbb F$ representing
$\theta|\mathbb F\in A'$ extends uniquely to an automorphism of $F_n$
representing $\theta$.  To prove the existence of a principal set $X$
(containing $T_c^{\pm})$ for $A$ it suffices to prove the existence of
a principal set $X'$ (containing $T_c^{\pm})$ for $A'$.  If
$\Fix(\atwo')$ has rank less than two then the existence of $X'$
follows from (1) and (2) applied to $\atwo' \in A'$.  Suppose then
that $\Fix(\atwo')$ has rank at least two.  By the same logic, it is
sufficient to find a principal set $X''$ (containing $T_c^{\pm})$ for
$A'|\Fix(\atwo')$ and this exists by the inductive hypothesis and the
fact that $A'|\Fix(\atwo') $ has a generating set with fewer than $K$
elements.  \endproof

\begin{lemma} \label{torsion free}
An abelian subgroup $A$ that is generated by rotationless elements is
torsion free.
\end{lemma}

\proof If $\theta \in A$ is a torsion element then it is represented
by a finite order homeomorphism $f':G'\to G'$ of a marked graph $G'$.
Suppose that $X$ is a principal set for $A$ and that $P_1,P_2,P_3 \in
X$.  There is a lift $\ti f' : \Gamma' \to \Gamma'$ such that each
$P_i \in \Fix(\widehat {f'})$.  The line $L_{12}$ with endpoints $P_1$
and $P_2$ and the line $L_{13}$ with endpoints $P_1$ and $P_3$ are
$\ti f'_\#$-invariant and since $\ti f'$ is a homeomorphism they are
$\ti f'$-invariant.  The intersection $L_{12} \cap L_{13}$ is an $\ti
f'$- invariant ray and so is contained in $\Fix(\ti f')$.  It follows
that $L_{12} \subset \Fix(\ti f')$ and that the image of $L_{12}$ in
$G'$ is contained in $\Fix(f')$.  It therefore suffices to show that
every edge of $G'$ is crossed by at least one such line.

For any set $Y \subset \partial F_n$, let $C_Y$ be the set of
bi-infinite lines cobounded by pairs of elements of $Y$.  Let
$W_A = \cup C_X$ where the union is over all principal sets $X$ for
$A$ and let $\F$ be the smallest free factor system that carries
$W_A$.  It suffices to show that $\F =\{[[F_n]]\}$.  The proof of this
assertion is by induction on the cardinality $K$ of a given
rotationless generating set $S$ for $A$.

Assume to the contrary that $\F$ is proper and choose $\otwo \in S$.
There is a homotopy equivalence $\fG$ representing $\otwo$ as in
Theorem~\ref{comp sp exists} in which $\F$ is realized as a filtration
element $G_r$.  Lemma~3.25(2) implies that each $\Lambda\in \L(\otwo)$
is the accumulation set of an isolated point in $\Fix_N(\hat \atwo)$
for some $\atwo \in \PA(\otwo)$.  By Lemma~\ref{principal sets},
$\Lambda$ is carried by $\F$. Thus each stratum above $G_r$ is NEG.
Items (NEG) and (PER) of Theorem~\ref{comp sp exists} imply that every
edge $E$ of $G\setminus G_r$ has an orientation so that its initial
vertex is principal and so that its initial direction is fixed.
Choose a lift $\ti E$ of $E$ and a principal lift $\ti f : \Gamma \to
\Gamma$ that fixes the initial direction determined by $\ti E$.  There
is a ray that begins with $\ti E$ and converges to a point in
$\Fix_N(\hat f)$.  This follows from Lemma~\ref{iterates to} if $E$ is
not a fixed edge and from Lemma~3.26 of \cite{fh:recognition}
otherwise.  Let $\atwo$ be the principal automorphism corresponding to
$\ti f$. It suffices to show that each element of $C_{\Fix(\hat
\atwo)}$ is carried by $\F$.  This is obvious if $K=1$.  We have now
proved the basis step of our induction argument and may assume that $K
> 1$ and that $\F =\{[[F_n]]\}$ when $A$ has a rotationless generating
set with fewer than $K$ elements.

If $\Fix(\atwo)$ has rank zero then $\Fix(\hat \atwo)$ is contained in
a principal set for $A$ by Lemma~\ref{principal sets}(1) and
$C_{\Fix(\hat \atwo)}$ is carried by $\F$.  If $\Fix(\atwo)$ has rank
one with generator $c$ then Lemma~\ref{principal sets}(2) implies that
the line connecting $P$ to $T_c^+$ is carried by $\F$ for each $P\in
\Fix(\hat \atwo)$.  It follows that the line connecting any two points
of $\Fix(\hat \atwo)$ is carried by $\F$.

We may therefore assume that $\Fix(\atwo)$ has rank at least two.  Let
us show that $\Fix(\atwo)$ is carried by $\F$.  The inductive
hypothesis and the fact that $A|\Fix(\atwo)$ has a generating set with
fewer than $K$ elements implies that no proper free factor system of
$\Fix(\atwo)$ carries $W_{A|\Fix(\atwo)}$.  The Kurosh subgroup
theorem therefore implies that any free factor system of $F_n$ that
carries $W_{A|\Fix(\atwo)}$ also carries all of $\Fix(\atwo)$.  Since
$W_{A|\Fix(\atwo)} \subset W_A$ we conclude that $\Fix(\atwo)$ is
carried by $\F$.

Lemma~\ref{principal sets}(3) implies that for each $P \in \Fix(\hat
\atwo)$ there exists $Q \in \partial \Fix(\atwo)$ so that the line
connecting $P$ to $Q$ is carried by $\F$.  Since the line connecting
any two points in $\partial \Fix(\atwo)$ is carried by $\F$ it follows
that the line connecting any two points in $ \Fix(\hat \atwo)$ is
carried by $\F$.  \endproof

\begin{corollary} \label{is birot}  An abelian subgroup $A$ that is generated by rotationless elements is rotationless.
\end{corollary}

\proof Suppose that $\oone \in A$, that $ k > 1$ and that $\aone_k \in
\PA^{\pm}(\oone^k)$.  Choose $m \ge 1$ so that $\oone^{km}$ is
rotationless.  By Lemma~\ref{principal sets} there is a principal set
$X$ for $A$ with $X \subset \Fix(\hat \aone_k^m)$.  Let $s : A \to
\Aut(F_n)$ be the lift determined by $X$ and let $\aone= s(\oone)$.
Then $\aone^{km} = s(\oone^k)^m = \aone_k^m$ and so $\aone^k =
s(\oone_k) = \aone_k$ by Lemma~\ref{torsion free}.  To complete the
proof it suffices by Remark~\ref{surjection} to show that $\Fix(\hat
\aone) = \Fix(\hat \aone^{km})$.  Let $\mathbb F = \Fix(\aone^{km})$
and note that $\mathbb F$ is $s(A)$-invariant. Lemma~\ref{torsion
free} implies that $\aone$ is uniquely characterized by $\aone^{km} =
\aone_k^m$ and hence that $\aone$ is independent of the choice of $X$.
Parts (1) and (2) of Lemma~\ref{principal sets} therefore imply that
$\Fix(\hat \aone)$ contains each isolated point in $\Fix(\hat
\aone_k^m)$ and contains $\partial \mathbb F$ if $\mathbb F$ has rank
less than two.  If $\mathbb F$ has rank at least two then $\mathbb F
\subset \Fix(\aone)$ by Lemma~\ref{torsion free} applied to
$\oone|\mathbb F$.  \endproof

\begin{corollary} \label{rotationless subgroup}
For each abelian subgroup $A$ of $\Out(F_n)$, the set of
rotationless elements is a rotationless subgroup $A_R$ that has
finite index in $A$.
\end{corollary}

\proof This is an immediate corollary of Corollary~\ref{is birot} and
the fact that every element of $A$ has a rotationless iterate.
\endproof

\section{Generic Elements of rotationless abelian subgroups}  \label{s:generic}

In this section we define an embedding of a given rotationless
abelian subgroup $A$ into an integer lattice $\mathbb Z^N$ and say
what it means for an element of $A$ to be generic with respect to this
embedding.

\begin{definition}
Suppose that $X_1$ and $X_2$ are principal sets for $A$ that define
distinct lifts $s_1$ and $s_2$ of $A$ to $\Aut(F_n)$ and that
$T_c^{\pm} \in X_1 \cap X_2$.  Then $s_2(\otwo) =
i_c^{d(\otwo)}s_1(\otwo)$  for all $\otwo \in A$
and some $d(\otwo) \in \Z$; the assignment $\otwo \mapsto d(\otwo)$
defines a homomorphism that we call the {\em comparison homomorphism}
$\omega : A \to \Z$ determined by $X_1$ and $X_2$.
\end{definition} 

\begin{lemma}
For any rotationless abelian subgroup $A$ there are only finitely
many comparison homomorphisms $\omega : A \to \Z$.
\end{lemma}

\proof Distinct comparison homomorphisms must disagree on some basis
element of $A$ so we can restrict attention to those comparison
homomorphisms that disagree on a single element $\otwo \in A$. If
$\omega$ is defined with respect to $X_1, X_2$ and $c$ then $[c]_u$,
the unoriented conjugacy class of $c$, is an axis of $\otwo$.  As
$\otwo$ has only finitely many axes, we may restrict attention to
those comparison homomorphisms that are defined with respect to the
same $[c]_u$.  If $a \in F_n$ and $\omega'$ is defined with respect to
$\hat i_aX_1, \hat i_aX_2$ and $i_a(c)$ then $\omega'=\omega$.  We may
therefore restrict attention to comparison homomorphisms that are
defined with respect to the same $c$.  The number of such comparison
homomorphisms is bounded by the multiplicity of $[c]_u$ as an axis for
$\oone$ by Lemma~\ref{axes and lifts}.  \endproof

\begin{lemma} \label{ef homo}
If $A$ is a rotationless abelian subgroup then $\L(A) = \cup_{\oone
\in A} \L(\oone)$ is a finite collection of $A$-invariant laminations.
\end{lemma}

\proof Let $\{\otwo_1,\ldots,\otwo_K\}$ be a rotationless basis for
$A$.  If $\L(\oone) =\{\Lambda_1,\ldots,\Lambda_q\}$ and
$F(\Lambda_i)$ is the smallest free factor that carries $\Lambda_i$
then the $F(\Lambda_i)$'s are distinct by Lemma~3.2.4 of
\cite{bfh:tits1}.  Each $\otwo_j$ permutes the $\Lambda_i$'s by
Lemma~3.1.6 of \cite{bfh:tits1} and so permutes the $F(\Lambda_i)$'s
by Lemma~\ref{ffs invariance}.  Since $\otwo_j$ is rotationless,
each $F(\Lambda_i)$, and hence each $\Lambda_i$, is
$\otwo_j$-invariant by Lemma~\ref{no iterates necessary}.  This proves
that $\Lambda_i$ is $A$-invariant and hence that $PF_{\Lambda_i}$ is
defined on $A$.  Each $PF_{\Lambda_i}$ must be non-zero when applied
to some $\otwo_j$ and by Lemma~3.3.1 of \cite{bfh:tits1} this is
equivalent to $\Lambda \in \L(\otwo_j) \cup \L(\otwo_j^{-1})$, which
is a finite set.  
\endproof

\begin{definition}
For each $\Lambda \in \L(A)$, we say that $PF_{\Lambda}|A$ is the {\em
expansion factor homomorphism for $A$} determined by $\Lambda$.  Let
$N$ be the number of distinct comparison and expansion factor
homomorphisms for $A$.  Define $\Omega : A \to \Z^N$ to be the product
of these homomorphisms.  We say that $\Omega$ is the {\em coordinate
homomorphism for $A$} and that each comparison homomorphism
 and expansion factor homomorphism is a
{\em coordinate} of $\Omega$.
\end{definition}

\begin{lemma} \label{embeds}
If $A$ is a rotationless abelian subgroup then $\Omega : A \to \Z^N$
is injective.
\end{lemma}

\proof Given non-trivial $\theta \in A$, choose $\fG$ and $\filt$
representing $\theta$ as in Theorem~\ref{comp sp exists} and let $H_l$
be the lowest non-fixed irreducible stratum.  If $H_l$ is \EG\ then
$PF_{\Lambda}(\theta) \ne 0$ for the attracting lamination $\Lambda
\in \L(\theta)$ associated to $H_l$.  Otherwise $H_l$ is a single edge
$E$ and $f(E) = E\cdot u$ where $u \subset G_{l-1}$ is a loop that is
fixed by $f$.  Lemma~\ref{axes and lifts} implies that there are
distinct principal lifts $\Theta_1$ and $\Theta_2$ of $\theta$ that
fix the primitive element $c\in F_n$ determined by $u$.  Thus
$\Theta_2 = T_c^d \Theta_1$ for some $d \ne 0$.  By
Lemma~\ref{principal sets} there exists principal sets $X_1 \subset
\Fix(\hat \Theta_1)$ and $X_2 \subset \Fix(\hat \Theta_2)$ that
contain $c$.  These determine a comparison homomorphism $\omega$ such
that $\omega(\theta) = d$. We have shown that some coordinate of
$\Omega(\theta) \ne 0$ and since $\theta$ was arbitrary, $\Omega$ is
injective.  \endproof

\begin{definition}
Assume that $A$ is a rotationless abelian subgroup and
that $\Omega : A \to \Z^N$ is its coordinate homomorphism. Then $\oone
\in A$ is {\em generic} if all coordinates of $\Omega(\oone)$ are
non-zero.
\end{definition}

\begin{com}
$\oone$ is generic in $A$ if and only if $\L(\oone) = \L(A)$ and \lq
$\oone$ has the same axes and multiplicity as $A$\rq.
\end{com} 

\begin{lemma} \label{generic basis}
Every rotationless abelian subgroup $A$ has a basis of generic
elements.
\end{lemma}

\proof Given a basis $\otwo_1,\ldots,\otwo_K$  for $A$ and $\theta \in A$ let $NZ(\theta) \subset
\{1,\ldots,N\}$ be the non-zero coordinates of
$\Omega(\theta)$. For all but finitely many
positive integers $a_2$, $NZ(\otwo_1 \otwo_2^{a_2}) = NZ(\otwo_1) \cup
NZ(\otwo_2)$.  Inductively choose positive integers $a_i$ for $i > 1$
so that $\atwo_1' := \atwo_1 \atwo_2^{a_2}\cdot \ldots \cdot
\atwo_K^{a_K}$ satisfies $NZ(\otwo_1') = \cup_{i=1}^K NZ(\otwo_i) =
\{1,\ldots,N\}$.  Replacing
$\otwo_1$ with $\otwo_1'$ produces a new basis in which the first
element is generic.  Repeating this step $K$ times produces the
desired basis or just replace $\otwo_2$ with
$\otwo_2\otwo_1^{big}$ and so on.
 \endproof

A principal set $X$ for $A$ determines a lift of $A$ to $\Aut(F_n)$.
If $X_1 \subset X_2$ are principal sets for $A$ then $X_1$ and $X_2$
determine the same lift.  It therefore makes sense to consider
principal sets that are maximal with respect to inclusion.

\begin{lemma} \label{generic fixed sets}
If $\oone \in A $ is generic then $\{\Fix(\hat \aone): \aone \in
\PA^{\pm}(\oone)\}$ is the set of maximal (with respect to inclusion)
principal sets for $A$.
\end{lemma}

\proof Each principal set $X'$ for $A$ determines a lift $s : A \to
\Aut(F_n)$.  If $\aone \in \PA^{\pm}(\oone)$ and $\Fix(\hat \aone)
\subset X'$ then $s(\oone)=\aone$ and $X' \subset \Fix(\hat \aone)$.
This proves that $\Fix(\aone)$ is a maximal principal set if it is a
principal set.  It therefore suffices to show that each $\Fix(\hat
\aone)$ is a principal set.

If $ \mathbb F :=\Fix(\aone)$ has rank zero then $\Fix(\hat \aone)$ is
a principal set by Lemma~\ref{principal sets}(1).  If $ \mathbb F$ has
rank one with generator $c$ and with isolated points $P, Q \in
\Fix(\hat \aone)$ then by Lemma~\ref{principal sets}(2) there is a
maximal principal set $X_P$ that contains $P$ and $T_c^{\pm}$ and a
maximal principal set $X_Q$ that contains contains $Q$ and
$T_c^{\pm}$.  If $X_P \ne X_Q$ then the comparison homomorphism that
they determine evaluates to zero on $\oone$ since $\aone_P =\aone_Q
=\aone$ in contradiction to the assumption that $\oone$ is
generic. Thus $X_P = X_Q $.  Since $P$ and $Q$ are arbitrary, $X_P
=\Fix(\hat \aone)$.

Suppose finally that $\mathbb F$ has rank at least two.  We claim that
$A|\mathbb F$ is trivial.  If not, let $\Omega'$ be the homomorphism
defined on $A|\mathbb F$ as the product of expansion factor and
comparison homomorphisms that occur for $A|\mathbb F$.  Each
coordinate $\omega'$ of $\Omega'$ extends to a coordinate $\omega$ of
$\Omega$. Since $\oone|\mathbb F$ is the identity, $\omega(\oone) = 0$
in contradiction to the assumption that $\oone$ is generic.  Thus
$A|\mathbb F$ is trivial and $\partial \mathbb F$ is contained in a
maximal principal set $X$ for $A$.

By Lemma~\ref{principal sets}(3), each isolated point $P$ in
$\Fix(\hat \aone)$ is contained in a maximal principal
set $X_P$ whose intersection $Y$ with $\partial \mathbb F$ contains at
least two points.  If $X_P \ne X$ then $Y$ has exactly two points and
in fact equals $\{T_b^{\pm}\}$ for some $b \in F_n$ since every lift
of the identity outer automorphism is an inner automorphism.  The
comparison homomorphism $\omega$ determined by $X_P$ and $X$ evaluates
to $0$ on $\oone$ in contradiction to the assumption that $\oone$ is
generic.  Thus $X_P = X$ for all isolated points $P$ and $\Fix(\hat
\aone) = X$ as desired.  \endproof

It is an immediate corollary, that from the point of view of fixed
points of principal lifts, generic elements are
indistinguishable.  

\begin{corollary} \label{same fixed points}
For any generic $\oone,\otwo \in A$ there is a bijection $h :
\PA^{\pm}(\oone) \to \PA^{\pm}(\otwo)$ such that $\Fix(\hat \aone) =
\Fix(\widehat{h(\aone)})$ for all $\aone \in \PA^{\pm}(\oone)$.
\end{corollary}  

\section{$\A(\oone)$}  \label{s:Aone}
The data required in the Recognition theorem has both qualitative and
quantitative components.  If we fix the qualitative part and allow the
quantitative part to vary then we generate an abelian group that is
naturally associated to the outer automorphism being considered.  This
section contains a formal treatment of this observation.  A more
computational friendly approach in terms of relative train track maps
is given in the next section.
 
\begin{definition} \label{associated abelian}
Assume that $\oone$ is rotationless.  $\A(\oone)$ is the subgroup of
$\Out(F_n)$ generated by rotationless elements $\othree$ for which
there is a bijection $h: \PA^{\pm}(\oone) \to \PA^{\pm}(\othree)$
satisfying $\Fix(\widehat{h(\aone)}) = \Fix(\hat \aone)$ for all $\hat
\aone \in \PA(\oone)$.
\end{definition}

\begin{remark}
It is an immediate consequence of the definitions that
$\A(\oone) = \A(\oone^k)$ for all $k \ne 0$ for all rotationless
$\oone$.
\end{remark} 

\begin{lemma}  \label{synthetic is everything}
If $A$ is a rotationless abelian subgroup and $\oone$ is generic in
$A$, then $A \subset \A(\oone)$.
\end{lemma}

\proof Lemma~\ref{generic basis} and Corollary~\ref{same fixed points}
imply that there is a generating set of $A$ that is contained in
$\A(\oone)$.  \endproof

To prove that $\A(\oone)$ is abelian we appeal to the following
characterization of the rotationless elements in the centralizer
$C(\oone)$ of $\oone$.

\begin{lemma} \label{criteria for commuting}  
If $\oone,\otwo \in \Out(F_n)$ are rotationless, then $\otwo \in
C(\oone)$ if and only if all for $\aone \in \PA^{\pm}(\oone)$:
\begin{description}  
\item [($\Phi-1$)] there exists $\atwo \in \PA^{\pm}(\otwo)$ such that
$\Fix(\hat \aone)$ is $\hat \atwo$-invariant.
\item [($\Phi-2$)] If $P \in \Fix(\hat \aone)$ is isolated then one
may choose $\atwo$ in ($\Phi-1$) such that $P \in \Fix(\hat \atwo)$.
\item [($\Phi-3$)] If $a \in \Fix(\aone)$ and $[a]_u$ is an axis of
$\oone$ then one may choose $\atwo$ in ($\Phi - 1$) such that $a \in
\Fix(\atwo)$.
\end{description}
Moreover, if $\otwo \in C(\oone)$ and $\atwo$ is as in ($\Phi-1$) then
$\atwo$ commutes with $\aone$.
\end{lemma}

\proof If $\otwo \in C(\oone)$, let $A = \langle \oone,\otwo \rangle$.
Lemma~\ref{principal sets} implies that for each $\aone \in
\PA^{\pm}(\oone)$, there is a principal set $X$ for $A$ whose
associated lift $s:A \to \Aut(F_n)$ satisfies $s(\oone) = \aone$.
Then $ s(\otwo) \in \PA^{\pm}(\otwo)$ commutes with $\aone$ and
($\Phi-1$) is satisfied.  ($\Phi-2$) follows from Lemma~\ref{principal
sets}.  If $[a]_u$ is an axis of $\oone$ then $[a]_u$ is
$\otwo^k$-invariant for some $k >0$ and so is $\otwo$-invariant by
Lemma~\ref{no iterates necessary}.  Items (2) and (3) of
Lemma~\ref{principal sets} allow us to choose $X$ to contain
$T_a^{\pm}$ which implies ($\Phi-3$).  This completes the only if
direction of the lemma.

For the if direction, we assume that $\otwo$ satisfies the three
items, define $\oone' := \otwo \oone \otwo^{-1}$ and prove that
$\oone' = \oone$ by applying the Recognition Theorem.
 
For each $\aone \in \PA(\oone)$ choose $\atwo_1$ satisfying $(\Phi -
1)$ and define $\aone' = \atwo_1 \aone \atwo_1^{-1} \in \PA(\oone')$.
If $\atwo_2$ also satisfies $(\Phi - 1)$ then $\atwo_2 = \atwo_1 i_x $
where $\Fix(\hat \aone)$ is $\widehat{i_x}$-invariant.  By
Lemma~\ref{l: first from bk3}, $x \in \Fix(\aone)$.  Thus $\atwo_2
\aone \atwo_2^{-1} = \atwo_1 i_x\aone i_x^{-1}\atwo_1^{-1} = \atwo_1
\aone \atwo_1^{-1}$ and $\aone'$ is independent of the choice of
$\atwo_1$.  We denote $\aone \mapsto \aone'$ by $h : \PA(\oone) \to
\PA(\oone')$ and note that $\Fix(\widehat{h(\aone)}) = \hat
\atwo_1(\Fix(\hat \aone)) = \Fix(\hat \aone)$ and that
$\Fix_N(\widehat{h(\aone)}) = \Fix_N(\hat \aone)$.  In particular, $h$
is injective.  If $\aone$ is replaced by $i_c \aone i_c^{-1}$ then
$\atwo_1$ can be replaced by $i_c \atwo_1 i_c^{-1}$ and $\aone'$ is
replaced by $i_c \aone' i_c^{-1}$.  Thus the restriction of $h$ to an
equivalence class in $\PA(\oone)$ is a bijection onto an equivalence
class in $\PA(\oone')$.  Lemma~\ref{conjugation basics}(2) implies
that $\PA(\oone)$ and $\PA(\oone')$ have the same number of
equivalence classes and hence that $h$ is a bijection.

Suppose that $\aone_1 \in \PA(\oone)$, that $a \in \Fix(\aone_1)$ is
primitive and that $\aone_2:= i_{a}^d\aone_1 \in \PA(\oone)$ for some
$d \ne 0$.  Then $[a]_u$ is an axis for $\oone$ and by $(\aone_1- 3)$
and $(\aone_2- 3)$ we may choose $\atwo_1$ for $\aone_1$ and $\atwo_2$
for $\aone_2$ to fix $a$.  Thus $\atwo_2 = i_{a}^m\atwo_1$ for some
$m$ and $\aone_2' = i_a^m \atwo_1 i_a^d\aone_1 \atwo_1^{-1} i_{a}^{-m}
= i_a^d\atwo_1 \aone_1 \atwo_1^{-1} = i_a^d\aone_1'$ which proves that
$h$ satisfies Theorem~\ref{t:recognition}-2(ii).

By Lemma~\ref{iterates to}, for each $\Lambda \in \L(\oone)$ there
exists $\aone \in \PA(\oone)$ and an isolated point $P \in \Fix_N(\hat
\aone)$ whose accumulation set equals $\Lambda$.  By $(\aone -2)$, we
may assume that $P$ is $\hat \atwo_1$-invariant and hence that
$\Lambda$ is $\otwo$-invariant.  It follows that $\Lambda$ is
$\oone'$-invariant and that $\PF_{\Lambda}(\oone') =
\PF_{\Lambda}(\oone)$.  Theorem~\ref{t:recognition} implies that
$\oone = \oone'$ and since $\Fix_N(\aone') = \Fix_N(\aone)$, $\aone =
\aone'$, which proves that $\atwo$ commutes with $\aone$.  \endproof

We denote the {\em center} of a group $H$ by $Z(H)$ and define the
{\em weak center} $WZ(H)$ to be the subgroup of $H$ consisting of
elements that commute with some iterate of each element of $H$.

\begin{corollary} \label{in the center}
If $\oone \in \Out(F_n)$ is rotationless then $\A(\oone)$ is an
abelian subgroup of $C(\oone)$.  Moreover, each element of $\A(\oone)$
commutes with each rotationless element of $C(\oone)$ and so
$\A(\oone) \subset WZ(C(\oone))$.  
\end{corollary}

\proof Lemma~\ref{criteria for commuting} implies that $\theta \in
C(\oone)$ for each $\theta$ in the defining generating set of
$\A(\oone)$ and that $C(\oone)$ and $C(\theta)$ contains the same
rotationless elements.  The corollary follows.  \endproof

\begin{remark}
In general, $\A(\oone)$ is not contained in the center of $C(\oone)$.
For example, if $n = 2k$ and $\aone \in \PA^{\pm}(\oone)$ commutes
with an order two automorphism $\Theta$ that interchanges the free
factor generated by the first $k$ elements in a basis with the free
factor generated by the last $k$ elements of that basis, then
$\A(\oone)$ will contain elements that do not commute with
$\theta$.  
\end{remark}
  
It is natural to ask if $\oone$ is generic in $\A(\oone)$.

\begin{lemma}\label{is generic}
If $\oone$ is rotationless then $\oone$ is generic in $\A(\oone)$.
\end{lemma}

\proof We must show that if $\omega$ is a coordinate of $\Omega :
\A(\oone) \to \Z^N$ then $\omega(\oone) \ne 0$.  Choose an element
$\othree$ of the defining generating set for $\A(\oone)$ such that
$\omega(\othree) \ne 0$.  If $\omega= PF_{\Lambda}$ then, after
replacing $\othree$ with $\othree^{-1}$ if necessary, $\Lambda \in
\L(\othree)$.  By Lemmas~\ref{some essential} and \ref{iterates to},
there exist $\athree \in \PA(\othree)$ and an isolated point $P \in
\Fix_N(\athree)$ whose accumulation set is $\Lambda$.  After replacing
$\oone$ with $\oone^{-1}$ if necessary, there exists $\aone \in
\PA(\oone)$ such that $\Fix(\hat \aone) = \Fix(\hat \athree)$ and such
that $P$ is an isolated point in $\Fix_N(\atwo)$.  Lemma~\ref{isolated
for egs} implies that $\omega(\oone) \ne 0$.

If $\omega$ is a comparison homomorphism determined by lifts $s,t
:\A(\oone) \to \Aut(F_n)$ then $s(\othree) \ne t(\othree)$.  Thus
$\Fix(\widehat{s(\oone)})= \Fix(\widehat{s(\othree)}) \ne
\Fix(\widehat{t(\othree)}) =\Fix(\widehat{t(\oone)}$ which implies
that $\omega(\oone) \ne 0$.  \endproof

\section{Disintegrating $\phi$} \label{s:disintegration}
We have reduced the study of rotationless abelian subgroups of
$\Out(F_n)$, and so of abelian subgroups of $\Out(F_n)$ up to finite
index, to the study of $\A(\oone)$ for rotationless $\oone \in
\Out(F_n)$.  In this section we construct the subgroup $\D(\oone)$ of
$\A(\oone)$ described in the introduction.  In section~\ref{s:finite
index} we show that $\D(\phi)$ has finite index in $\A(\phi)$.

Choose $\fG$ representing $\oone$ as in Theorem~\ref{comp sp exists}.
We will need a coarsening of the complete splitting of
a path. For each axis $\mu$ of $\phi$ there exists a primitive closed
Nielsen path $w$ and edges $\{E_i\}$ as in Lemma~\ref{axes and lifts}
such that $f(E_i) = E_i \cdot w^{d_i}$; we say that these edges are
{\em associated to $\mu$} and that $d_i$ is the {\em exponent of
$E_i$}.  For distinct $E_i$ and $E_j$ associated to $\mu$, paths of
the form $E_i w^* \bar E_j$ are said to {\em belong to the same \qe\
family} or to {\em be \qe}.  By assumption $d_i \ne d_j$.  If $d_i$
and $d_j$ have the same sign then $E_i w^* \bar E_j$ is an exceptional
path but otherwise it is not.

Assume that $\sigma = \sigma_1 \cdot \ldots \cdot \sigma_s$ is the
unique complete splitting of $\sigma$.  If $\sigma_{ab}:= \sigma_a
\cdot \ldots \cdot \sigma_b$ is quasi-exceptional then we say that
$\sigma_{ab}$ is a {\em $QE$-subpath} of $\sigma$.

\begin{lemma}   For any completely split path $\sigma$, distinct $QE$-subpaths 
of $\sigma$ have disjoint interiors.
\end{lemma}

\proof Suppose that $\sigma = \sigma_1 \cdot \ldots \cdot \sigma_s$ is
the complete splitting of $\sigma$ and that there exist $1 \le a < b
\le s$ and $1 \le a \le c < d \le s$ such that $\sigma_{ab}:= \sigma_a
\cdot \ldots \cdot \sigma_b$ and $\sigma_{cd}:=\sigma_c \cdot \ldots
\cdot \sigma_d$ are distinct quasi-exceptional paths.  We must show
that $c >b$.

Since $\sigma_{ab}$ is \qe, $\sigma_a$ and $\sigma_b$ are linear edges
and $\sigma_{a+1}\cdot\ldots\cdot \sigma_{b-1}$ is a Nielsen path.
Each $\sigma_l$, $a< l < b$ must be a Nielsen path, which implies that
$c \ge b$.  Since $\bar E_j$ is the not an initial segment of any \qe\
path, $c > b$.  \endproof

\begin{definition}
The {\em \ds}\ of a completely split path $\sigma$ is the coarsening
of the complete splitting of $\sigma$ obtained by declaring each
$QE$-subpath to be a single element.  Thus the \ds\ is a splitting
into single edges, connecting subpaths Nielsen paths
and quasi-exceptional paths. These subpaths are the {\em terms of
the \ds.}
\end{definition}
  
\begin{definition}  \label{almost invariant}
For a stratum $H_{\alpha}$ whose edges are not fixed by $f$, we let
$A_{\alpha} \subset H_{\alpha}$ denote an edge if the stratum
$H_{\alpha}$ is irreducible and a connecting path if $H_{\alpha}$ is a
zero stratum.  The rule
\begin{itemize}
\item $H_i \sim H_j$ if there exist $A_i \subset H_i$ and $A_j \subset
H_j$ such that $A_j$ occurs as a  term in the \ds\ of
$f(A_i)$.
\end{itemize}
 generates an equivalence relation on those strata on which $f$ is not the 
identity.  The equivalence classes   $X_1, 
\dots, X_M$ are called  {\em almost invariant subgraphs}.

\vspace{.1in}

For each $M$-tuple $\va$ of non-negative integers, define $f_{\va} : G
\to G$ by
$$f_{\va}(E) = \left\{ \begin{array}{ll} f_{\#}^{{a}_i}(E) & \mbox
{for each edge $E \subset X_i$} \\ E & \mbox{for each edge $E$ that is
fixed by $f$}\end{array} \right.$$  
\end{definition}

\begin{remark} \label{zero strata with EG strata}
If $H_i$ is a zero stratum and $H_j$ is the first irreducible stratum
above $H_i$ then $H_i$ and $H_j$ belong to the same almost invariant
subgraph.  This follows from the fact that every edge in $H_i$ is
contained in a connecting path in $H_i$ that is in the image of
either an edge in $H_j$ or a connecting path in some zero stratum
between $H_i$ and $H_j$.
\end{remark}

\begin{lemma} \label{homeq}
$f_{\va} : G \to G$ is a homotopy equivalence for all ${\va}$.
\end{lemma}

\proof
Let $NI$ be the number of irreducible strata in the filtration and for
each $0 \le m \le NI$, let $G_{i(m)}$ be the smallest filtration
element containing the first $m$ irreducible strata.  We will prove by
induction that each $f_{\va}|G_{i(m)}$ is a homotopy equivalence.

Since $H_1$ is never a zero stratum, $i(1) = 1$. If $G_1$ is not a
single edge fixed by $f$, then every edge in $G_1$ is contained in a
single almost invariant subgraph $X_i$.  Thus $f_{\va}|G_1$ is either
the identity or is homotopic to $f^{{a}_i}|G_1$; in either case it is
a homotopy equivalence.

We assume now that $f_{\va}|G_{i(m)}$ is a homotopy equivalence.
Define $g_1 : G_{i(m+1)} \to G_{i(m+1)}$ on edges by
$$
g_1(E)  = \left\{ \begin{array}{ll} f_{\va}(E) & \mbox {if $E 
\subset G_{i(m)}$}  \\  E & \mbox{if $E \subset G_{i(m+1)} \setminus 
G_{i(m)}$}\end{array} \right.
$$ Every vertex in $G_{i(m)}$ whose link is not entirely contained in
$G_{i(m)}$ is  an
attaching vertex (see Theorem 2.18(V)) and so is fixed by $f$.  This
guarantees that $g_1$ is well defined.  It is easy to check that $g_1$
is a homotopy equivalence.  If the edges of $H_{i(m+1)}$ are fixed by
$f$, then $g_1 = f_{\va}|G_{i(m+1)}$ and we are done.

If $f|H_{i(m+1)}$ is not the identity, then  
Remark 6.4 implies that the edges in $G_{i(m+1)} \setminus G_{i(m)}$
are contained in a single almost invariant subgraph, say $X_k$. Define
$g_2 : G_{i(m+1)} \to G_{i(m+1)}$ on edges by
$$
g_2(E)  = \left\{ \begin{array}{ll} f_{\#}^{a_k}(E) & \mbox {if $E 
\subset G_{i(m)}$}  \\  E & \mbox{if $E \subset G_{i(m+1)} \setminus    
G_{i(m)}$}\end{array} \right.
$$ 
and $g_3 : G_{i(m+1)} \to G_{i(m+1)}$ on edges by
$$
g_3(E)  = \left\{ \begin{array}{ll} E & \mbox {if $E 
\subset G_{i(m)}$}  \\  f_{\#}^{a_k}(E) & \mbox{if $E \subset G_{i(m+1)} 
\setminus G_{i(m)}$}\end{array} \right.
$$ Then $g_2$ is a homotopy equivalence and $f^{a_k}|G_{i(m+1)} = g_3
g_2$. Each component of $G_{i(m+1)}$ is non-contractible by item (Z)
of Theorem~\ref{comp sp exists}, so $f^{a_k}|G_{i(m+1)}$ is a homotopy
equivalence.  It follows that $g_3$, and hence also
$f_{\va}|G_{i(m+1)} = g_3 g_1$ is a homotopy equivalence. \qed

\vspace{.1in}  

Almost invariant subgraphs are defined
without reference to the \qe\ paths in the \ds\ of edge images.  The
next definition brings these into the discussion.

\begin{definition} \label{consistency}
If $\{X_1,\ldots,X_M\}$ are the almost invariant subgraphs of $\fG$
then an $M$-tuple $a = (a_1,\dots, a_M)$ of non-negative integers is
{\em admissible} if for all axes $\mu$, whenever:
\begin{itemize}
\item $X_s$ contains an edge $E_i$ associated to $\mu$ with exponent
$d_i$
\item $X_t$ contains an edge $E_j$ associated to $\mu$ with exponent
$d_j$
\item $X_r$ contains an edge $E_k$ such that some term of
the \ds\ of $f(E_k)$ is in the same \qe\ family as $E_i \bar E_j$
\end{itemize}
then  $a_r(d_i-d_j) = a_sd_i-a_td_j$.
\end{definition}

\begin{ex} 
  Suppose that $G$ is the rose with edges
$E_1,E_2,E_3$ and $E_4$ and that $f : G \to G$ is defined by $E_1
\mapsto E_1$, $E_2 \mapsto E_2E_1^2$, $E_3 \mapsto E_3E_1$ and $E_4
\mapsto E_4E_3E_3\bar E_2$.  Then $M=2$ with $X_1$ having the
single edge $E_2$ and $X_2$ consisting of $E_3$ and $E_4$.  The pair
$(a_1,a_2)$ is admissible if $a_2 = 2a_2-a_1$ or equivalently $a_2 =
a_1$.  Thus $f_{\va} = f^{a_1}$ for each admissible 
  $\va$. 
\end{ex}

\begin{definition}
   Each $f_{\va}$ determines an element an element $\phi_{\va} \in \ofn$ and also an element 
$[f_{\va}]$ in  the semigroup of homotopy
equivalences of $G$ that respect the filtration  modulo
homotopy relative to the set of vertices of $G$.  
Define $\D(\oone) = \langle \phi_{\va} : a \mbox{ is
admissible}\rangle$.  Both $\oone_{\va}$ and $\D(\oone)$ depend on the
choice of $\fG$; see Example~\ref{D depends on f} below.  Since we
work with a single $\fG$ throughout the paper and since $\D(\oone)$ is
well defined up to finite index by Theorem~\ref{D has finite index},
we suppress this dependence in the notation.
\end{definition}

\begin{ex} \label{D depends on f}
 Let $G$ be the rose with edges $E_1,E_2$ and
$E_3$.  Define $f_1 :G \to G$ by
$$ 
E_1 \mapsto E_1 \qquad  E_2 \mapsto E_1E_2 \qquad  E_3 \mapsto E_1^2E_3E_1
$$
and
$f_2 :G \to G$ by
$$ 
E_1 \mapsto E_1 \qquad  E_2 \mapsto E_2E_1 \qquad  E_3 \mapsto E_1E_3E_1^2.
$$ These maps differ by $i_{E_1}$ and so determine the same element
$\oone \in \Out(F_n)$.  The homotopy equivalence of $G$ that fixes
$E_1$ and $E_3$ and maps $E_2$ to $E_2E_1$ represents an element of
$\D(\oone)$ if $f_2$ is used but not if $f_1$ is used.
\end{ex}

\begin{lemma} \label{still Nielsen}
Suppose that $a$ is admissible and that $\sigma$ is a path in $G$.
\begin{enumerate}
\item If $\sigma$ is Nielsen path for $f$ then $\sigma$ is a Nielsen
path for $f_{\va}$.
\item If $\sigma$ is \qe\ and if some path in the same \qe\ family as
$\sigma$ occurs as a term in the \ds\ of $f(E)$ for some
edge $E$ in $X_k$ then $(f_{\va})_\#(\sigma) = f^{a_k}_\#(\sigma)$.
\end{enumerate}
\end{lemma}

\proof The proof is by induction on the height $r$ of $\sigma$. In
the context of (1), we may assume that $\sigma$ is either indivisible
or a single fixed edge.

$G_1$ is either a single fixed edge or is contained in a single almost
invariant subgraph.  Thus $f_{\va}|G_1$ is either the identity or an
iterate of $f|G_1$.  In either case (1) is obvious for $\sigma \subset
G_1$. Since $G_1$ does not contain any \qe\ paths, the lemma holds for
$\sigma \subset G_1$.  We assume now that $r \ge 2$, that the lemma
holds for paths in $G_{r-1}$ and that $\sigma$ has height $r$ and is
either an indivisible Nielsen path or a quasi-exceptional
path. Property (N) of Theorem~\ref{comp sp exists} implies that $H_r$
is either \eg\ or linear.
  
Let $X_s$ be the almost invariant subgraph containing $H_r$.  If $H_r$
is linear then it is a single edge $E_r$ and $f(E_r) = E_r w^{d_r}$
for some non-trivial Nielsen path $w$.  If $\sigma$ is an \iNp, then
$\sigma = E_r w^p \bar E_r$ for some integer $p$.  By the inductive
hypothesis, $(f_{\va})_\#(w) = w$ so
    
$$(f_{\va})_\#(\sigma) = [(E_r w^{a_sd_r}) w^p (\bar w^{a_sd_r}\bar
E_r)] = E_r w^p \bar E_r = \sigma.
$$

If $\sigma $ is as in (2), then up to a reversal of orientation,
$\sigma = E_r w^p \bar E_j$ for some edge $E_j$ associated to the same
axis as $E_r$.  Let $X_t$ be the almost invariant subgraph containing
$E_j$.  Since $a$ is admissible, $a_k(d_r-d_j) = a_sd_r -a_td_j$.

Thus \newline\noindent
\begin{eqnarray*}
(f_{\va})_\#(\sigma) & = & [f_\#^{a_s}(E_r) (f_{\va}(w))^p
f_\#^{a_t}(\bar E_j)] \\& = & [E_rw^{a_sd_r} w^p \bar w^{a_td_j}\bar
E_j] \\& = & [E_rw^{a_sd_r-a_td_j+p}\bar E_j] \\& = & [E_r
w^{a_k(d_r-d_j)+p}E_j] \\& = & [ E_rw^{a_kd_r} w^p \bar w^{a_kd_j}
\bar E_j] \\& = & [f_\#^{a_k}(E_r) (f_\#^{a_k}(w))^p f_\#^{a_k}(\bar
E_j)]\\& = & f_\#^{a_k}(\sigma).
\end{eqnarray*}

Suppose now that $H_r$ is \eg. There are no \qe\ paths of height $r$
so $\sigma$ is an \iNp\ of height $r$.  By item (2) of Lemma~5.11 of
\cite{bh:tracks}, $\sigma = \alpha \beta$ where $\alpha$ and $\beta$
are $r$-legal paths that begin and end with edges in $H_r$.  It
suffices to show that $(f_{\va})_\#(\alpha) = f^{a_s}_\#(\alpha)$ and
$(f_{\va})_\#(\beta) = f^{a_s}_\#(\beta)$. The argument is the same
for both $\alpha$ and $\beta$. If $E_{\alpha}$ is the initial edge of
$\alpha$, then there exists $m > 0$ such that $\alpha \subset
f^m_\#(E_{\alpha})$.  The terms in the \qe\ splitting of
$f^m_\#(E_{\alpha})$ are edges and connecting paths in $X_s$, Nielsen
paths in $G_{r-1}$ or \qe\ paths in $G_{r-1}$.  Since $\alpha$ begins
and ends with an edge in $H_r$, the \qe\ splitting of
$f^m_\#(E_{\alpha})$ restricts to a \qe\ splitting of $\alpha$.  By
definition and by the inductive hypothesis, $(f_{\va})_\#$ equals
$f^{a_s}_\#$ on all four types of subpath.  Thus $(f_{\va})_\#$ equals
$f^{a_s}_\#$ on $\alpha$ as desired.  \qed\

\vspace{.1in}

\begin{corollary} \label{local control}
For $1 \le s \le M$, let $P_s$ be the set of completely split paths
whose \qe\ splittings are composed of: (i) edges and connecting paths
in $X_s$;\; (ii) \iNp s; and (iii) \qe\ paths in the same \qe\ family
as a term in the \ds\ of $f(E)$ for some edge $E$ in $X_s$.  Then
$P_s$ is preserved by both $f_\#$ and $(f_{\va})_\#$ and moreover
$(f_{\va})_\#(\sigma)= f^{a_s}_\#(\sigma)$ for all $\sigma \in P_s$.
\end{corollary} 

\proof This is an immediate consequence of Lemma~\ref{still Nielsen},
the definition of $X_s$ and the definition of $f_{\va}$.  \endproof

\begin{cor} \label{commute} 
For each admissible $a$ and $b$, $[f_{\va}][f_{\vb}] =
[f_{\vb}][f_{\va}] = [f_{\va+\vb}].$ In particular, $\D(\oone)$ is
abelian.
\end{cor}    

\proof It suffices to check that $(f_{\va+\vb})_\#(E) =
(f_{\va})_\#(f_{\vb})_\#(E)$ for each edge $E$.  If $E$ is fixed by
$f$, then $E$ is also fixed by $f_{\va}$, $f_{\vb}$ and $f_{\va+\vb}$.
Suppose that $E \subset X_k$.  Then
$$ (f_{\va+\vb})_\#(E) = f_\#^{(\va+\vb)_k}(E) = f_\#^{a_k+b_k}(E)
= f_\#^{a_k}f_\#^{b_k}(E) = (f_{\va})_\#(f_\#^{b_k}(E) )=
(f_{\va})_\#(f_{\vb})_\#(E)
$$   where the next to the last equality comes
from Corollary~\ref{local control}. \qed

\begin{definition}
An admissible $\va$ is {\em generic} if each $a_i > 0$ and if whenever
$E_i \in X_r$ and $E_j \in X_s$ are distinct linear edges associated
to the same axis, then $a_rd_i \ne a_sd_j$ where $d_i$ and $d_j$ are
the exponents of $E_i$ and $E_j$ respectively.
\end{definition}

\begin{lemma} \label{is cs}
If $\va$ is generic then $f_{\va} :G \to G$ satisfies the conclusions
of Theorem~\ref{comp sp exists} and $f_{\va}$ has the same principal
vertices and Nielsen paths as $f$.
\end{lemma}

\proof Corollary~\ref{local control} implies that $f_{\va}$ is a \cs\
\rtt\ for $\phi_{\va}$ with respect to the filtration $\filt$, that
$f_{\va}$ has the same principal vertices as $f$ and that $f_{\va}$
satisfies all of the items listed in the statement of
Theorem~\ref{comp sp exists} except perhaps for (L) and (N).  Property
(L) follows from the genericity assumption and the observation that if
$E_j \in X_s$ is a linear edge for $f$ with exponent $d_j$ then $E_j$
is a linear edge for $f_{\va}$ with exponent $a_sd_j$.

We show below that if $\rho$ is an \iNp\ for $f_{\va}$ then it is a
Nielsen path for $f$.  Combined with Lemma~\ref{still Nielsen}(1),
this proves that $f$ and $f_{\va}$ have the same Nielsen paths and
hence that (N) is satisfied.  Corollary~\ref{local control} then
implies that $f_{\va}$ is forward rotationless and completes the proof.
(Nielsen paths are relevant to this because it is part of the
definition of forward rotationless that the endpoints of all indivisible
periodic Nielsen paths be vertices.)

Suppose then that $\rho$ is an \iNp\ for $f_{\va}$. Let $i$ be the
height of $\rho$ and let $X_r$ be the almost invariant subgraph that
contains $H_i$.  If $H_i$ is EG then by Lemma~5.11 of \cite{bh:tracks}
$\rho = \alpha \beta$ where $\alpha$ and $\beta$ are $i$-legal paths
for $f_{\va}$ that begin and end in $H_i$. Let $E_{\alpha}$ be the
edge whose interior contains an initial segment of $\alpha$.  If the
initial endpoint $x$ of $\alpha$ is a vertex let $\alpha' = \alpha$;
otherwise $\alpha'$ is the extension of $\alpha$ that contains all of
$E_{\alpha}$.  Choose $k \ge 1$ so that $(f^k_{\va})_\#(E_{\alpha})$
contains $\alpha$.  Since both $E_{\alpha}$ and the terminal edge of
$\alpha$ are edges of height $i$ (see Theorem 2.18(N)), the \qe\
splitting of $(f^k_{\va})_\#(E_{\alpha})$ restricts to a \qe\
splitting of $\alpha'$.  Corollary~\ref{local control} implies that
$(f_{\va})_\#(\alpha') = f_\#^{a_r}(\alpha')$ and since $\alpha$ and
$\alpha'$ are $i$-legal it follows that $(f_{\va})_\#(\alpha) =
f_\#^{a_r}(\alpha)$.  The analogous argument applies to $\beta$ and we
conclude that $\rho$ is an indivisible periodic Nielsen path for $f$,
and so by (N), an indivisible Nielsen path for $f$.

Suppose next that $H_i$ is a single NEG edge $E_i$.  Choose lifts $\ti
\rho$ and $\ti f_{\va}$ such that $\ti f_{\va}$ fixes the endpoints of
$\ti \rho$.  Let $\ti f$ be the lift of $f$ that fixes the initial
endpoint and direction determined by $\ti \rho$.  By
Lemma~\ref{iterates to} there is a ray $\ti R_1$ with the same initial
vertex and direction as $\ti \rho_1$ and satisfying the following
properties.
\begin{itemize}
\item $\Fix(\ti f) \cap \ti R_1$ is the initial endpoint of $\ti R_1$.
\item If $\ti R_1 = \ti \tau_1 \cdot \ti \tau_2 \cdot \ldots$ is the
\qe\ splitting of $\ti R_1$ and if $\ti x_l$ is the terminal endpoint
of $\ti \tau_l$ then $\ti f(\ti x_l)= \ti x_k$ for some $k > l$ and
$Df$ maps the turn taken by $\ti R_1$ at $\ti x_l$ to the turn taken
by $\ti R_1$ at $\ti x_k$.
\item $\ti R_1$ converges to some $P_1 \in \Fix_N(\hat f)$.
\end{itemize}

We claim that these three items also hold with $\ti f$ replaced by
$\ti f_{\va}$.  The second and third items follow from
Corollary~\ref{local control}. If the first item fails then a
fixed point in the interior of $\ti R_1$ must be in the interior of
some $\ti \tau_l$ that is not a connecting subpath and so there exists
an initial subpath $\ti \mu$ of $\ti \tau_l$ such that $(\ti
f_a)_\#(\ti \mu)$ is trivial.  But no such $\ti \mu$ can exist.  This
follows from Corollary~\ref{local control} if $\tau_l$ is a single
edge and is easy to check by inspection if $\tau_l$ is an exceptional
path or a Nielsen path.  This completes the proof of the claim.  We
now know that $\ti R_1$ does not contain the terminal endpoint of $\ti
\rho$.

Define $\ti R_2$ and $P_2$ similarly using the initial vertex and
direction of $\ti \rho^{-1}$.  If $P_1 \ne P_2$ let $\ti L_{12}$ be
the line connecting $P_1$ to $P_2$.  Then $\ti L_{12}$ is contained in
$\ti R_1 \cup \ti \rho \cup \ti R_2$ and does not contain the
endpoints of $\ti \rho$, which are also the endpoints of $\ti R_1$ and
$\ti R_2$.  It follows that $\ti L_{12} \cap \Fix(\ti f_{\va})=
\emptyset$ which contradicts (Lemma~3.15 of \cite{fh:recognition}) the
fact that the two endpoints of $\ti L_{12}$ are attracting.  We
conclude that $P_1=P_2$.

If $P_1 \ne T_b^{\pm}$ for some $b \in F_n$, then there is a unique
$\hat f$ that fixes $P_1$.  In this case, the lifts of $f$ that fix
the initial and terminal endpoints of $\ti \rho$ are equal and $\rho$
is a Nielsen path for $f$.  We may therefore assume that $P_1 =
T_b^{\pm}$ in which case $E_i$ and $E_j$ are linear edges associated
to the same axis for $f$ and $\rho$ is exceptional for $f$. Property
(L) for $f_{\va}$ implies that $E_i = E_j$ and hence that $\rho$ is a
Nielsen path for $f$.  \endproof

We now relate $\D(\oone)$ to $\A(\oone)$, using the correspondence
between principal lifts of relative train track maps and principal
automorphisms.

\begin{corollary} \label{distinct suffixes}
For each generic $\va$ there is a bijection $h : \PA(\phi) \to
\PA(\phi_{\va})$ such that $\Fix_N(\widehat{h(\aone)}) =\Fix_N(\hat
\aone)$ for all $\hat \aone \in \PA(\oone)$.  If $\ti f$ corresponds
to $\aone$ and $\ti f_{\va}$  corresponds to $h(\aone)$
then $\Fix(\ti f) = \Fix(\ti f_{\va})$.
\end{corollary}

\proof By Lemma~\ref{is cs}, $f$ and $f_{\va}$ have the same Nielsen
classes of principal vertices. There is an induced bijection $h$
between principal lifts of $ f_{\va}$ and principal lifts of $f$; if
$\ti f_{\va} = h(\ti f)$ then $\Fix(\ti f) = \Fix(\ti
f_{\va})$. Lemma~\ref{l: second from bk3} implies that $\Fix_N(\hat
f)$ and $\Fix_N(\hat f_{\va})$ have the same non-isolated points.
Lemma~\ref{iterates to} and Corollary~\ref{local control} imply that
$\Fix_N(\hat f)$ and $\Fix_N(\hat f_{\va})$ have the same isolated
points. \endproof

Let $\D_R(\oone)$ be the finite index rotationless subgroup of
$\D(\oone)$ given by Corollary~\ref{rotationless subgroup}.

\begin{corollary}  \label{D is contained in A} 
$\D_R(\oone)$ is contained in $\A(\oone)$ and is generated by elements
of the form $\oone_{\va}$ with $\va$ generic.
\end{corollary}

\proof Corollary~\ref{distinct suffixes} implies that if $\oone_{\va}
\in \D_R(\oone)$ is generic then $\phi_{\va} \in \A(\oone)$.  It
therefore suffices to find a generating set $S$ for $\D_R (\oone)$ in
which each element has this form.  Let $S' = \{\oone_{\va}\}$ be any
generating set for $\D_R(\oone)$.  If $\mathbf I$ is the $M$-tuple
with $1$'s in each coordinate then $\oone_{\mathbf I}=\oone$ is
generic and represented by $f_{\mathbf I}=f$.  Corollary~\ref{commute}
implies that if $k$ is sufficiently large then $\oone^k \oone_{\va}$
is represented by $f_{\vb}$ where $\vb= \va + k{\mathbf I}$ is
projectively close to ${\mathbf I}$ and so is generic.  Thus $S =
\{\oone, \oone^k\oone_{\va}: \oone_{\va} \in S'\}$ is the desired
generating set for $\D_R(\oone)$.  \endproof

The definition of $\D(\oone)$ is not symmetric in $\oone$ and
$\oone^{-1}$ leaving open the following questions.

\begin{question}
Is each element of $\D(\oone)$ rotationless? Is $\D(\oone) =
\D(\oone^{-1})$ ?
\end{question}

\begin{thm}  \label{WZC} 
$\D_R(\phi) \subset WZ(C(\phi))$ for all rotationless $\phi$.
\end{thm}

\proof  $\D_R(\oone) \subset \A(\phi) \subset WZ(C(\phi))$ by Corollary~\ref{D is contained in A} and Corollary~\ref{in the center} .    
\endproof

\section{Finite Index}  \label{s:finite index}
Our goal in this section is to prove  

\begin{thm} \label{D has finite index}
$\D_R(\oone)$ has finite index in $\A(\oone)$ for all rotationless
$\oone$.
\end{thm} 

Before turning to the proof of Theorem~\ref{D has finite index}  we use it to prove one of our main results.

\begin{thm}  \label{disintegration groups are enough} 
For every abelian subgroup $A$ of $\Out(F_n)$ there exists $\phi \in
A$ such that $A \cap D(\oone)$ has finite index in $A$.
\end{thm}

\proof Corollary~\ref{rotationless subgroup} and Lemma~\ref{synthetic is everything} imply that  $A \cap \A(\phi)$ has finite index in $A$ for each generic $\phi \in A_R$.  Theorem~\ref{D has finite index} therefore completes the proof.  
\endproof

Choose once and for all $\fG$ representing $\oone$ as in
Theorem~\ref{comp sp exists}.

We set notation for the linear edges associated to an axis $[c]_u$ of
$\oone$ as follows.  If $[c]_u$ has multiplicity $m+1$ then there is a
primitive closed path $w$ whose circuit represents $c$ and for $1 \le
j \le m$, there are linear edges $E_j$ and distinct non-zero integers
$d_j$ such that $f(E_j) = E_j \cdot w^{d_j}$.  Choose a lift $\ti E_j$
whose terminal endpoint is in the axis $A_c \subset \Gamma$.
Following Lemma~\ref{axes and lifts}, the principal lift of $f$ that
fixes the initial endpoint of $\ti E_j$ is denoted $\ti f_j$ and the
associated principal automorphism is denoted $\aone_j$; both $\ti f_j$
and $\aone_j$ are independent of the choice of $\ti E_j$.  By
Lemma~\ref{generic fixed sets} and Lemma~\ref{is generic}, $\Fix(\hat
\aone_j)$ is a maximal principal set $X_j$.  The lift $s_j$ of
$\A(\oone)$ to $\Aut(F_n)$ determined by $X_j$ satisfies $s_j(\oone) =
\aone_j$.  The principal lift of $f$ that fixes the terminal endpoint
of $\ti E_j$ is denoted $\ti f_0$, its associated principal
automorphism is denoted $\aone_0$, the maximal principal set
$\Fix(\hat \aone_0)$ is denoted $X_0$ and the lift to $\Aut(F_n)$
determined by $X_0$ is denoted $s_0$.  The automorphisms
$\aone_0,\ldots, \aone_m$ are the only elements of $\PA(\oone)$ that
commute with $T_c$ (Lemma~\ref{axes and lifts}).

For $1 \le j \ne k \le m$, let $\omega_{c,j}$ be the comparison
homomorphism determined by $X_0$ and $X_j$ and let $\omega_{c,j,k}$ be
the comparison homomorphism determined by $X_j$ and $X_k$.  Thus
$\omega_{c,j,k} = \omega_{c,j} - \omega_{c,k}$.  There is an obvious
bijection between the $\omega_{c,j}$'s and the linear edges $E_j$
associated to $c$.  There is also a bijection between the
$\omega_{c,j,k}$'s and the families of quasi-exceptional paths
$E_jw^*\bar E_k$ associated to to $c$. We make use these bijections
without further notice.

For each $\Lambda \in \L(\oone)$ let $\omega_{\Lambda}
=\PF_{\Lambda}|\A(\oone)$.  We also identify $\Lambda$ with
$\omega_{\Lambda}$ when convenient.

We define a new homomorphism $\Omega^{\oone}: \A(\oone)\to \Z^K$ whose
coordinates are in one to one correspondence with the the linear and
EG strata of $\fG$ by removing extraneous coordinates from $\Omega:
\A(\oone)\to \Z^N$.

\begin{definition} \label{better Omega}
$\Omega^{\oone} : \A(\oone) \to \Z^K$ is the product of the
$\omega_{c,j}$'s and the $\omega_{\Lambda}$'s as $[c]_u$ varies over
the axes of $\oone$ and as $\Lambda$ varies over $\L(\oone)$.
\end{definition}  

\begin{lemma} \label{still embeds}
$\Omega^{\oone} : \A(\oone) \to \Z^K$ is injective.
\end{lemma}

\proof The coordinates of $\Omega^{\oone}$ are coordinates of the
injective homomorphism $\Omega$.  It therefore suffices to assume that
$\omega(\otwo) \ne 0$ for a coordinate $\omega$ of $\Omega$ and prove
that the image of $\otwo$ under some coordinate of $\Omega^{\oone}$ is
non-zero.  There is no loss in assuming that $\omega$ is not a
coordinate of $\Omega^{\oone}$ and is therefore either some
$\omega_{c,j,k}$ or $\omega_{\Lambda}$ for some $\Lambda \in
\L(\oone^{-1})$.  In the former case, $\omega_{c,j}(\otwo) \ne 0$ or
$\omega_{c,k}(\otwo) \ne 0$ and we are done.  In the latter case,
Lemma~\ref{is generic} implies that $\Lambda \in \L(\otwo) \cup
\L(\otwo^{-1})$.  By Lemma~3.2.4 of \cite{bfh:tits1} there is a unique
$\Lambda' \ne \Lambda \in \L(\oone) \cup \L(\oone^{-1})$ such that
$\Lambda$ and $\Lambda'$ are carried by the same minimal rank free
factor; moreover, $\Lambda' \in \L(\oone)$.  Similarly, there is a
unique $\Lambda'' \ne \Lambda\in \L(\otwo) \cup \L(\otwo^{-1})$ such
that $\Lambda$ and $\Lambda''$ are carried by the same minimal rank
free factor.  Since $\oone$ is generic, $\omega_{\Lambda''}(\oone) \ne
0$ which implies that $\Lambda'' \in \L(\oone) \cup \L(\oone^{-1})$
and hence that $\Lambda' = \Lambda''$.  Thus $\omega_{\Lambda'}$ is a
coordinate of $\Omega^{\oone}$ and $\omega_{\Lambda'}(\otwo) \ne 0$.
\endproof

\begin{lemma} \label{coordinate map evaluation}
If a coordinate $\omega$ of $\Omega^{\oone}$ corresponds to a stratum
in the almost invariant subgraph $X_s$ then $\omega(\oone_{\va}) =
a_s\omega(\oone)$ for all $\phi_{\va} \in \A(\oone)$.
\end{lemma}

\proof We may assume by Corollary~\ref{D is contained in A} that $\va$
is generic.  If $\omega = \omega_{\Lambda}$ then the lemma follows
from Corollary~\ref{local control} and the definition of the expansion
factor homomorphism.  Suppose then that $\omega = \omega_{c,j}$.
Lemma~\ref{distinct suffixes} implies that $s_j(\oone_{\va})$
corresponds to the principal lift of $f_{\va}$ that fixes the initial
endpoint of $\ti E_j$ and $s_0(\oone_{\va})$ corresponds to the
principal lift of $f_{\va}$ that fixes the terminal endpoint of $\ti
E_j$.  Since $f_{\va}(E_j) = E_j\cdot w^{a_sd_j}$ we have
$\omega_{c,j}(\oone_{\va}) = a_sd_j$.  \endproof
 
\begin{cor}  The rank of $D_R(\oone)$ is equal to the rank  of the sublattice $L$ of $\Z^M$ generated by the admissible $M$-tuples for $\fG$.
\end{cor}

\proof Let $\rho:L \to \A(\phi)$ be the homomorphism determined by
$\va \to \phi_{\va}$.  It suffices to show that $\rho$ is injective
and for this it suffices to show that $\Omega^{\oone} \rho$ is
injective.  The lattice $L$ contains the $M$-tuple $I$, all of whose
coordinates are $1$. Given $x,y \in L$ there exists $k \ge 0$ so that
$x+kI$ and $y+kI$ are admissible.  Lemma~\ref{still embeds} and
Lemma~\ref{coordinate map evaluation} imply that $\Omega^{\oone}
\rho(x+kI) \ne \Omega^{\oone} \rho(y+kI)$ and hence that
$\Omega^{\oone} \rho(x) \ne \Omega^{\oone} \rho(y)$.  \endproof
 
We now come to our main technical proposition, a generalization of
Lemma~\ref{isolated for egs}. Before proving it we show that it
implies Theorem~\ref{D has finite index}.  (The process of iterating
an edge is discussed in section~\ref{s:background}.)

\begin{proposition}\label{main technical}
Suppose that $E$ is a nonlinear edge of $G$, that $k >0$ and that
$\mu$ is a term in the \qe\ splitting of $f^k(E)$ that is either a
linear edge, an EG edge or a \qe\ subpath.  Suppose further that $\ti
f : \Gamma \to \Gamma$ is a principal lift of $f$ that fixes the
initial vertex of a lift $\ti E$ of $E$ and that the ray $\ti R$
determined by iterating $\ti E$ converges to $P \in \Fix_N(\hat f)$.
Let $\omega$ be the homomorphism associated to $\mu$ and let $s$ be
the lift of $\A(\oone)$ to $\Aut(F_n)$ determined by the maximal
principal set $\Fix(\hat f)$.  Then the following are equivalent for
all $\otwo \in \A(\oone)$.
\begin{enumerate}
\item  $P$ is isolated in $\Fix(\widehat{s(\otwo)})$
\item $\omega(\psi) \ne 0$.
\end{enumerate} 
\end{proposition}

\noindent{\bf Proof of Theorem~\ref{D has finite index}} Each linear
or EG stratum $H_i$ determines a coordinate of $\Omega^{\oone}$ that
we denote $\omega_i$.  For all $\otwo \in \A(\oone)$, define
$a_i(\otwo) = \omega_i(\otwo)/\omega_i(\oone)$.

We first observe that there is a virtual basis $\{\otwo_l\}$ for
$\A(\oone)$, meaning that it is a basis of a finite index subgroup of
$\A(\oone)$, such that $a_i(\otwo_l)$ is a positive integer for all
$i$ and $l$.  To construct $\{\otwo_l\}$, start with any basis
$\{\eta_l\}$ of $\A(\oone)$.  Choose $m \ge 1$ so that each
$a_i(\eta_l^m)$ is an integer.  Then $\{\eta_l^m\}$ is a virtual basis
and for all but at most one value of $s$, the set obtained from
$\{\eta_l^m\}$ by replacing $\eta_1^m$ with $\otwo_1 =
\oone^{s}\eta_1^m$ is also a virtual basis for $\A(\oone)$.  If $s$ is
sufficiently large then $a_i(\otwo_1)$ is a positive integer.  Repeat
this, focusing on the second basis element and so on to arrive at the
desired virtual basis.

In what follows we restrict to a single $\otwo_l$ so we refer to
$\otwo_l$ simply as $\otwo$ and to $a_i(\otwo_l)$ simply as $a_i$.

We show next that if $H_i$ and $H_j$ are linear or EG strata that
belong to the same almost invariant subgraph then $a_i= a_j$.  Define
$\theta = \otwo \oone^{-{a_i}}$.  Then $\omega_j(\theta) =
\omega_j(\otwo)-a_i \omega_j(\oone)$ and $\omega_i(\theta) =
\omega_i(\otwo)-a_i \omega_i(\oone) = 0$ so it suffices to show that
$\omega_j(\theta) = 0$.

As a first case, suppose that $H_i$ is EG, that some, and hence every,
edge in $H_j$ occurs as a term in the \qe\ splitting of an iterate of
some, and hence every, edge in $H_i$.  By Remark~\ref{some essential}
there is an edge $E_i$ in $H_i$ whose initial vertex is principal and
whose initial direction is fixed.  Choose a lift $\ti E_i$, let $\ti
f$ be the principal lift that fixes the initial endpoint of $\ti E_i$,
let $P \in \Fix(\hat f)$ be the terminal endpoint of the ray obtained
by iterating $\ti E_i$ by $\ti f$ and let $s : \A(\oone) \to
\Aut(F_n)$ be the lift determined by the maximal principal set
$\Fix(\hat f)$.  Proposition~\ref{main technical} with $E = \mu= E_i$
and the assumption that $\omega_i(\theta) = 0$ imply that $P$ is not
isolated in $\Fix(\widehat{s(\theta)})$.  A second application of
Proposition~\ref{main technical}, this time with $E =E_i$ and $\mu$ an
edge in $H_j$ implies that $\omega_j(\theta)= 0$.

As a second case, suppose that there is a non-linear NEG edge $E_l$
and that edges $E_i$ of $H_i$ and $E_j$ of $H_j$ occur as terms in the
\qe\ splitting of an iterate of $E_l$.  Define $P \in \Fix(\hat f)$
and $s : \A(\oone) \to \Aut(F_n)$ as in the previous case using a lift
of $E_l$ instead of a lift of $E_i$.  As in the previous case
Proposition~\ref{main technical} can be applied with $E= E_l$ and with
either $\mu= E_i$ or $\mu= E_j$.  Thus $\omega_i(\theta)= 0$ if and
only if $\omega_j(\theta)= 0$ as desired.

The equivalence relation on strata that defines almost invariant
subgraphs is generated by these two cases so we have shown that the
$a_i$'s determine a well defined $M$-tuple $\hat{\va} = (\hat
a_1,\ldots,\hat a_M)$ with one $\hat a_s$ for each almost invariant
subgraph $X_s$.  To show that $\hat {\va}$ is admissible, assume that
$E_i \in X_s, E_j\in X_t$ and $E_k \in X_r $ are as in
Definition~\ref{consistency}.  As in the previous cases, we may assume
that the initial vertex of $E_k$ is principal and the initial
direction of $E_k$ is fixed.

Define $\eta = \otwo \oone^{-a_k}$. As in the previous cases,
Proposition~\ref{main technical} can be applied with $E= E_k$ and with
either $\mu= E_k$ or $\mu$ equal to an element in the \qe\ family
determined by $E_i \bar E_j$.  Since $\omega_k(\eta) =
\omega_k(\otwo)-a_k \omega_k(\oone) = 0$,  it follows that $0 = \omega_{c,i,j}(\eta) =
\omega_{c,i,j}(\otwo)-a_k \omega_{c,i,j}(\oone)$.  Keeping in mind
that $\hat{a}_r =a_k$, we have
$$
\omega_{c,i,j}(\otwo) = \hat a_r \omega_{c,i,j}(\oone). 
$$
Combining this with

$$ \hat a_r(d_i-d_j) = \hat a_r(\omega_i(\oone)-\omega_j(\oone)) =
\hat a_r(\omega_{c,i,j}(\oone))
$$
and 
$$ \omega_{c,i,j}(\otwo) = \omega_i(\otwo)-\omega_j(\otwo) = \hat a_s
\omega_i(\oone) - \hat a_t \omega_j(\oone) = \hat a_s d_i -\hat a_t
d_j
$$ proves that $\hat a$ is admissible.  Choose $K\ge 1$  so that
$\oone^K_{\hat {\va}} = \oone_{K\va}$ is rotationless.
Corollary~\ref{D is contained in A} implies that $\oone^K_{\hat {\va}}
\in \A(\oone)$ and Lemma~\ref{still embeds} then implies that $\otwo^K
= \oone^K_{\hat {\va}} \in \D_R(\oone)$. Thus $\{\otwo_l^K\} \subset
\D_R(\oone)$ is a virtual basis for $\A(\oone)$.  \endproof

The remainder of the section is devoted to the proof of
Proposition~\ref{main technical}.  For motivation we consider the
proof as it applies to a simple example.

\begin{ex}
Suppose that $G = R_3$ with edges $A,B$ and $C$ and that $\fG$
representing $\oone$ is defined by $A \mapsto A$, $B \mapsto BA$\ and
$C \mapsto CB$.

Let $T_A$   be the covering
translation corresponding to $A$ and let $\ti B$ be a lift of $B$ with
terminal endpoint in the axis of $T_A$.  Denote the principal lifts of $f$ that
fix the initial and terminal endpoints of $\ti B$ by $\ti f_-$ and
$\ti f_+$ respectively.  The fixed point sets $X_{\pm}$ of $\hat
f_{\pm}$ are maximal principal sets for $A(\oone)$ and so determine
lifts $s_{\pm} : \A(\oone) \to \Aut(F_n)$ such that $X_{\pm} \subset
\Fix(\widehat{s_{\pm}(\otwo)})$ for all $\otwo \in \A(\oone)$.  The
coordinate homomorphism $\omega$ corresponding to $B$ satisfies
$\omega(\otwo) = 0$ if and only if $s_+(\otwo) = s_-(\otwo)$.  Note
that $T_A^{\pm}$ is contained in both $X_+$ and $X_-$.

Choose a lift $\ti C$ of $C$ and let $\ti f$ be the principal lift
that fixes its initial endpoint.  Iterating $\ti E$ by $\ti f$
produces a ray $\ti R$ that converges to some $P \in \Fix(\hat f)$ and
that projects to an $f$-invariant ray $R = CBBABA^2\ldots BA^l
BA^{l+1} BA^{l+2}\ldots$.  The maximal principal set $\Fix(\hat f)$
determines a lift $s:\A(\oone)$ to $\Aut(F_n)$.  Denote the subpath
$BBABA^2$ of $R$ that follows the initial $C$ by $\sigma_0$ and the
subpath $f_\#^l(\sigma_0) = BA^l BA^{l+1} BA^{l+2}$ of $R$ by
$\sigma_l$.  There are lifts $\ti \sigma_l \subset \ti R$ of
$\sigma_l$, $l \to \infty$, that are cofinal in $\ti R$ and so limit
on $P$.

There are also lifts $\ti \delta_l$ of $\sigma_l$ for which $\ti B$ is
the edge that projects to the middle $B$ in $\sigma_l$.  The endpoints
of $\ti \delta_l$ are denoted $\ti x_l$ and $\ti y_l$.  The path
connecting $\ti x_l$ to the initial endpoint of $\ti B$ is a lift of
$BA^l$ and the path connecting the terminal endpoint of $\ti B$ to
$T_A^{-l}\ti y_l$ is a lift of $ABA^{l+2}$.  Thus $\ti x_l \to Q_- \in
X_{-}\setminus T_A^{\pm}$ and \ $T_A^{-l}\ti y_l \to Q_+ \in
X_{+}\setminus T_A^{\pm}$.  The line connecting $Q_-$ to $Q_+$
projects to $A^{\infty}BABA^{\infty}$.

Choose $g : G' \to G'$ representing $\otwo$ as in Theorem~\ref{comp sp
exists}.  The lift $\ti g: \Gamma' \to \Gamma'$ corresponding to
$s(\otwo)$ satisfies $P \in \Fix(\hat g)$.  For simplicity, we
suppress the equivariant map that identifies $\Gamma$ with $\Gamma'$.

If $P$ is not isolated in $\Fix(\hat g)$ then Lemma~\ref{l: second
from bk3} implies that $\ti g$ moves the endpoints of $\ti \sigma_l$
by an amount that is bounded independently of $l$.  Since $\ti
\delta_l$ is a translate of $\ti \sigma_l$ there is a lift $\ti g_l$
of $g$ that moves $\ti x_l$ and $\ti y_l$ by a uniformly bounded
amount, say $\kappa$.  In Lemma~\ref{for the non-isolated case} below
we show that under these circumstances, $\ti g_l$ commutes with $T_A$.
Since there is a lift of $g$ that commutes with $T_A$ and fixes $Q_-$,
it follows that $\ti g_l(Q_-) = T_A^{d_l}(Q_-)$.  If $d_l \ne 0$ and
$x_l$ is sufficiently close to $Q_-$ then the distance between $x_l$
and $\ti g_l(\ti x_l)$ would be greater than $\kappa$ which is a
contradiction.  Thus $d_l = 0$ and $Q_- \in \Fix(\widehat {g_l})$ for
all sufficiently large $l$.  A second consequence of the fact that
$\ti g_l$ commutes with $T_A$ is that $\ti g_l$ moves $T_A^{-l}\ti
y_l$ by a uniformly bounded amount.  Arguing as in the previous case
we conclude that $Q_+ \in \Fix(\widehat {g_l})$ for all sufficiently
large $l$.  For these $l$, $\Fix(\widehat{g_l})$ intersects both $X_+$
and $X_-$ in at least three points which implies that $\ti g_l$ is the
lift associated to both $s_-(\otwo)$ and $s_+(\otwo)$ and hence that
$\omega(\otwo) = 0$.

If $P$ is isolated in $\Fix(\hat g)$ then by Lemma~\ref{iterates to}
there is an edge $\ti E'$ of $\Gamma'$ that iterates toward $P$ under
the action of $\ti g$.  The ray $\ti R'$ connecting $\ti E'$ to $P$
eventually agrees with $\ti R$ and so contains $\ti \sigma_l$ for
large $l$.  Lemma~\ref{increasing length} below states, roughly
speaking, that since iterating $E'$ by $g$ produces segments of the
form $BA^lB$ for arbitrarily large $l$, it must be that $g_\#(BAB) =
BA^kB$ for some $k > 0$.  This implies that $A^{\infty}BABA^{\infty}$
is not $g_\#$-invariant and hence that the lifts of $g'$ corresponding
to $s_-(\otwo)$ and to $s_+(\otwo)$  are distinct.
Equivalently, $\omega(\otwo) \ne 0$.
 \end{ex}

We now turn to the formal proof. 

\begin{com} \label{long overlaps commute}
For the following lemmas it is useful to recall that if the circuits
representing $[b]$ and $[c]$ have edge length $L_b$ and $L_c$ and if
$A_b \cap A_c$ has edge length at least $L_b + L_c$ then $T_c$
commutes with $T_b$ because the initial endpoint $\ti x$ of $A_b \cap
A_c$ satisfies $T_bT_c(\ti x) = T_cT_b(\ti x)$.  It follows that $A_b
= A_c$ and that $T_b = T_c^{\pm}$.
\end{com}   

\begin{lemma} \label{for the non-isolated case}
Suppose that $\otwo \in \ofn$ is rotationless and that $g : G' \to
G'$ represents $\otwo$ and satisfies the conclusions of
Theorem~\ref{comp sp exists}.  Then for any primitive covering
translation $T_c$ of the universal cover $\Gamma'$ of $G'$, there
exists $K>0$ with the following property.  If $\tau \subset G'$ is a
Nielsen path for $g$ and $\ti \tau \subset \Gamma'$ is a lift whose
intersection with the axis $A_c$ of $T_c$ contains at least $K$ edges,
then the lift $\ti g$ that fixes the endpoints of $\ti \tau$ commutes
with $T_c$.
\end{lemma}

\proof Choose $L$ greater than the number of edges in each of the
following:
\begin{itemize}
\item [(1)] the loop in $G'$ that represents $c$
\item [(2)] each of the loops in $G'$ representing an  axis of $\otwo$
\item [(3)] any indivisible Nielsen path associated to an EG 
stratum for $g :G' \to G'$. 
\end{itemize}

There is a decomposition $\ti \tau = \ti \tau_1\cdot \ldots \cdot \ti
\tau_N$ into subpaths $\ti \tau_i$ that are either fixed edges or
indivisible Nielsen paths.  The endpoints of the $\ti \tau_i$'s are
fixed by $\ti g$.  There is no loss in assuming that each $\ti \tau_i$
intersects $A_c$ in at least an edge.

If $N \ge L$ then by (1), there exist $\ti \tau_i$ with initial
endpoint $\ti x$ and $\ti \tau_j$ with initial endpoint $T_c^l(\ti x)$
for some some $l \ne 0$.  Thus $\ti g T_c^l(\ti x) = T_c^l(\ti x)=
T_c^l \ti g(\ti x)$.  Since lifts of a map that agree on a point are
identical, $\ti g T^l_c = T_c^l \ti g$.  It follows that $\hat g$
fixes $T_c^{\pm}$ which then implies that $\ti g$ commutes with $T_c$.

We may therefore assume $N < L$.  In fact we may assume that $N =1$ :
if $K$ works in this case then $(L+2)K$ works in the general case. If
$\tau$ is a fixed edge then $K=2$ vacuously works.  We may
therefore assume that $\tau$ is indivisible.

Let $K =2L+2$. We may assume by (3) that $\tau$ is not associated to
an \eg\ stratum and so by Theorem~\ref{comp sp exists}(N), $\ti \tau =
\ti E_i \ti w^p \ti E_i^{-1}$ for some linear edge $E_i$ satisfying
$f(E_i) = E_iw^{d_i}$ where $w$ represents an axis $\mu$ of $\otwo$
and therefore has fewer than $L$ edges.  There is an axis $A_b$ for a
primitive $b \in F_n$ that contains $\ti w^p$ and whose projection
into $G$ is the loop determined by $w$.  Remark~\ref{long overlaps
commute} and our choice of $K$ imply that $T_b= T_c^{\pm}$.  It is
obvious that $\ti g$ commutes with $T_b$ so $\ti g$ also commutes with
$T_c$.  \endproof

Suppose that $E_i$ is a linear edge and that $f(E_i) = E_i w^{d_i}$.
If either $E_i$ or a \qe\ path $E_i w^*\bar E_j$ occurs as a term
in the \qe \ splitting of some $f^l_\#(\sigma)$ then $f^m_\#(\sigma)$
contains subpaths of the form $w^k$ where $k \to \pm\infty$ as $m
\to \infty$.  This is essentially the only way that such paths develop
under iteration.  Lemma~\ref{increasing length} below is an
application of this observation stated in the way that it is applied
in the proof of Proposition~\ref{main technical}.

We use $\EL(\cdot)$ to denote edge length of a path or circuit.  By
extension, for $c \in F_n$, we use $\EL(c)$ to denote the edge length
of the circuit representing $[c]$.

We isolate the following observation for easy reference.

\begin{lemma}  \label{stabilize}
Suppose that $g : G' \to G'$ satisfies the conclusions of
Theorem~\ref{comp sp exists} and that $\tau \subset G'$ is a
completely split path such that $EL(g_\#^m(\tau))$ is not uniformly
bounded. Then for all $L > 0$ there exists $M> 0$ so that for all $m
\ge M$, \ $EL(g_\#^m(\tau)) > 2L$ and the initial and terminal subpaths
of $g_\#^m(\tau)$ with edge length $L$ is are independent of $m$.
\end{lemma}

\proof The proof is by induction on the height $r$ of $\tau$.  The
$r=0$ case is vacuous so we may assume that the lemma holds for paths
of height less than $r$.  By symmetry it is sufficient to show that
$EL(g_\#^m(\tau))\to \infty$ and that initial segment of
$g_\#^m(\tau)$ with edge length $L$ stabilizes under iteration.

Let $\tau = \tau_1 \cdot \ldots \cdot \tau_s$ be the complete
splitting of $\tau$ and let $\tau_i$ be the first term such that
$EL(g_\#^m(\tau_i))$ is not uniformly bounded.  The terms preceding
$\tau_i$, if any, are Nielsen paths or pre-Nielsen connecting
paths. Their iterates stabilize so there is no loss in truncating
$\tau$ by removing them.  We may therefore assume that $i=1$.  It now
suffices to show that $EL(g_\#^m(\tau_1))\to \infty$ and that initial
segment of $g_\#^m(\tau_1)$ with edge length $L$ stabilizes under
iteration.  If $\tau_1$ is a connecting path this follows by induction
on $r$.  The remaining cases are that $\tau_1$ is a non-fixed edge in
an irreducible stratum or a \qe\ path and the result is clear in both
these cases.  \endproof

The following lemma is a case-by-case analysis of
the occurrence of long periodic segments in iterates of a single path.
The basic observation is that once a periodic segment reaches a
certain length it continues to get longer under further iteration.

\begin{lemma} \label{increasing length}
Suppose that $g : G' \to G'$ satisfies the conclusions of
Theorem~\ref{comp sp exists}, that $c \in F_n$ is primitive and that
$\ti g: \Gamma' \to \Gamma'$ is a lift of $g$ that commutes with $T_c$.
Then for all completely split paths $\sigma \subset G'$, there exists
$L_{\sigma} > 0$ so that if $m \ge 1$ and $\ti \rho_m$ is a lift of
$\rho_m = g^m_\#(\sigma)$ such that $\EL(\ti \rho_m \cap A_c) >
L_{\sigma}$ then $\EL(\ti g_\#(\ti \rho_m) \cap A_c) > \EL(\ti \rho_m
\cap A_c)$.
\end{lemma}

\proof Lemma~\ref{l: first from bk3} implies that the circuit
corresponding to $c$ is $g_\#$-invariant and hence that $A_c$
decomposes as a concatenation of subpaths that project to $g$-fixed
edges and indivisible Nielsen paths for $g$.  After composing $\ti g$
with an iterate of $T_c$ if necessary, we may assume that these
subpaths are $\ti g$-Nielsen paths. The endpoints of these paths are
called {\em splitting vertices} and their union is the set of $\ti
g$-fixed vertices in $A_c$.

The proof is by induction on the height $r$ of $\sigma$.  The
induction statement is enhanced to include the following property: if
$\EL(\ti \rho_m \cap A_c) > L_{\sigma}$ and if $\ti \rho_m \cap A_c$
contains an endpoint $\ti v$ of $\ti \rho_m$ then $\ti v$ is a
splitting vertex.

In certain cases we will show that $A_c \cap \ti \rho_m$ is uniformly
bounded, meaning that it is bounded independently of $m$.  One then
chooses $L_{\sigma}$ greater than that bound.  The $r=0$ case is
vacuously true so we may assume that the inductive statement holds for
all paths of height less than $r$.

Assume for now that there is only one term in the \ds\ of $\sigma$.
There are five cases, two of which are immediate.  If $\sigma$ is a
Nielsen path then $A_c \cap \ti \rho_m$ is uniformly bounded and we
are done.  If $\sigma$ is a connecting path then we let $L_{\sigma} =
L_{g_\#(\sigma)}$ where the latter exists by the inductive hypothesis
and the fact that $g_\#(\sigma)$ has height less than $r$.

If $\sigma$ is a linear edge $E$ then $\rho_m = Ew^{dm}$ for some
Nielsen path $w$ that forms a primitive circuit and some $d > 0$.  Let
$L_{\sigma} =EL(c) + EL(w)$.  If $EL(\ti \rho_m \cap A_c )>
L_{\sigma}$ then by Remark~\ref{long overlaps commute} there is a lift
$\ti w$ of $w$ such that $\ti \rho_m \cap A_c = \ti w^{dm}$ contains
all of $\ti \rho_m$ but the initial edge and $\ti g_\#(\ti \rho_m)
\cap A_c = \ti w^{d(m+1)}$.  Since $w$ is a Nielsen path and $\ti w$
is a fundamental domain of $A_c$ the endpoints of $\ti w$ are
splitting vertices.

If $\sigma$ is an exceptional path $E_i w^p \bar E_j$ where $g(E_i) =
E_iw^{d_i}$ and $g(E_j) = E_jw^{d_j}$, the proof is similar to the
linear case and we can use the same value of $L_{\sigma}$. If $EL(\ti
\rho_m \cap A_c )> L_{\sigma}$ then there is a lift $\ti w$ of $w$
such that $\ti \rho_m \cap A_c = \ti w^{m(d_i-d_j)+p}$ contains all of
$\ti \rho_m$ but the initial and terminal edges and $\ti g_\#(\ti
\rho_m) \cap A_c = \ti w^{(m+1)(d_i-d_j)+p}$.  In this case the
endpoints of $\ti \rho_m$ are not contained in $A_c$.

The fifth and hardest case is that $\sigma$ is a single edge $E$ in a
non-linear irreducible stratum $H_r$.  If the height of $A_c$ is
greater than $r$ then $A_c \cap \ti \rho_m$ has uniformly bounded
length.  We may therefore assume that $A_c$ has height at most $r$.
We consider the EG and NEG subcases separately.

Suppose that $H_r$ is EG. If $A_c$ has height $r$ then it has an
illegal turn in the $r$-stratum which implies that $EL(A_c \cap
\ti \rho_m) < EL(c)$.  We may therefore assume $A_c$ has height
less than $r$.  In particular, endpoints of $\ti \rho_m$ are not
contained in $A_c$.  For each edge $E'$ of $H_r$ there is a coarsening
of the \ds\ of $g(E')$ into an alternating concatenation of subpaths
in $ H_r$ subpaths in $G_{r-1}$.  Let $\{\mu_j\}$ be the set of paths
of $G_{r-1}$ that occur as subpaths in this decomposition as $E'$
varies over all edges of $H_r$.  The path $g^m_\#(E)$ also splits as
an alternating concatenation of subpaths in $ H_r$ and subpaths in
$G_{r-1}$; each of the subpaths in $G_{r-1}$ equals $g^l_\#(\nu_j)$
for some $\nu_j$ and some $0 \le l \le m$.  We may therefore choose
$L_{\sigma} = \max\{L_{\mu_j}\}$.

Finally, suppose that $H_r$ is non-linear and \noneg.  There is a path
$u \subset G_{r-1}$ such that $g^m(E) = E \cdot u \cdot g_\#(u)
\cdot\ldots\cdot g_\#^{m}(u)$ for all $m$ and such $EL(g_\#^{j}(u))
\to \infty$.  We may assume without loss that $A_c$ has height less
than $r$ and hence that $\ti \rho_m \cap A_c $ projects into $u \cdot
g_\#(u) \cdot\ldots\cdot g_\#^{m}(u)$ .  We claim that if $r$ is
sufficiently large, say $r > R$, then the projection of $\ti \rho_m
\cap A_c$ does not contain $g_\#^{r}(u)$ for any $m$.  Assume the
claim for now.  If $\EL(\ti \rho_m \cap A_c) > \EL(u \cdot g_\#(u)
\cdot\ldots\cdot g_\#^{R+1}(u))$ then the projection of $\ti \rho_m
\cap A_c$ is contained in $g_\#^{q-1}( u)\cdot g_\#^{q}( u) =
g^{q-1}_\#( u \cdot g_\#(u))$ for some $q$.  We may therefore choose
$L_{\sigma}$ to be the maximum of $\EL(u \cdot g_\#(u)
\cdot\ldots\cdot g_\#^{R+1}(u))$ and $L_{u \cdot g_\#(u)}$.
 
The claim is obvious is unless $u$ and $A_c$ have the same height, say
$t$, so assume that this is the case.  The claim is also obvious if
the maximal length of a subpath of $g_\#^{q}(u)$ with height less than
$t$ goes to infinity with $q$.  We may therefore assume that the
number of height $t$ edges in $g_\#^{q}(u)$ is unbounded.  Thus $H_t$
is EG and $g^r_\#(u)$ contains $t$-legal subpaths of length greater
than $EL(c)$ for all sufficiently large $r$. Since no such subpath is
contained in $A_c$ this completes the proof of the claim and so also
the induction step when there is only one term in the \ds\ of
$\sigma$.

Assume now that $\sigma= \sigma_1 \cdot \ldots \cdot \sigma_s$ is the
\ds\ of $\sigma$ and that $s> 1$.  Let $L_1 = \max\{L_{\sigma_i}\}$.
By Lemma~\ref{stabilize} there exists $M>0$ so that for all $m> M$ and
all $\sigma_i$, either $g_\#^m(\sigma_i)$ is independent of $m$ or
$EL(g^m_\#(\sigma_i)) > 2L_1$ and the initial and terminal segments of
$g^m_\#(\sigma_i)$ with edge length $L_1$ are independent of $m$.  The
former corresponds to $\sigma_i$ being a Nielsen path or a pre-Nielsen
connecting path and the latter to all remaining cases.  Choose
$L_{\sigma} > sL_1$ so that $EL(g^m_\#(\sigma)) < L_{\sigma}$ for all
$m \le M$.

Denote $g^m_\#(\sigma_i)$ by $\rho_{i,m}$ and write $\ti \rho_m = \ti
\rho_{1,m}\cdot \ldots \cdot \ti \rho_{s,m}$.  If $\EL(A_c \cap \ti
\rho) \ge L_{\sigma}$ then $\EL(A_c \cap \ti \rho_{i,m}) \ge
L_{\sigma_i}$ for some $1 \le i \le s$.  Thus $EL(g_\#(A_c \cap \ti
\rho_{i,m})) > EL(A_c \cap \ti \rho_{i,m})$.  If $A_c \cap \ti \rho
\subset \ti \rho_{i,m}$ we are done.  Otherwise we may assume that
$A_c \cap \ti \rho_{i+1,m}$ is a non-trivial initial segment of $\ti
\rho_{i+1,m}$ that begins at a splitting vertex of $A_c$.  If
$\rho_{i+1,m}$ is a Nielsen path or a pre-Nielsen connecting path then
$\ti g_\#(\ti \rho_{i+1,m}) = \ti \rho_{i+1,m}$ so $g_\#(A_c \cap \ti
\rho_{i+1,m}) = A_c \cap \ti \rho_{i+1,m}$.  This same equality holds
if $EL(A_c \cap \ti \rho_{i+1,m}) \le L_1$ by our choice of $M$.
Finally, if $EL(A_c \cap \ti \rho_{i+1,m}) > L_1$ then $EL(g_\#(A_c
\cap \ti \rho_{i+1,m})) > EL(A_c \cap \ti \rho_{i+1,m})$.  This
completes the proof if $A_c \cap \ti \rho_{m} \subset \ti
\rho_{i,m}\ti \rho_{i+1,m}$.  Iterating this argument completes the
proof in general.  \endproof

We need one more lemma before proving the main proposition.

\begin{lemma}  \label{limit exists}
Suppose that $g : G' \to G'$ satisfies the conclusions of
Theorem~\ref{comp sp exists}, that $\sigma$ is a completely split
non-Nielsen path for $g$ and that $\ti \sigma \subset \Gamma'$ is a
lift of $\sigma$ with endpoints at vertices $\ti x$ and $\ti y$.  If
$\ti g': \Gamma' \to \Gamma'$ fixes $\ti x$ then $\lim_{k\to
\infty}\ti g'^k(\ti y) \to Q$ for some $Q \in \Fix_N(\hat g)$.
\end{lemma}

\proof There is no loss in assuming that $\sigma$ is either a single
non-fixed edge or an exceptional path $E\tau^l \bar E'$.  In the
former case the lemma follows from Lemma~\ref{iterates to}.  In the
latter case, $Q$ is an endpoint of the axis of a covering translation
corresponding to $\tau$.  \endproof

\noindent{\bf Proof of Proposition~\ref{main technical}} The case that
$\mu$ is an EG edge follows from Lemma~\ref{isolated for egs}.  In the
remaining cases there is an axis $[c]_u$ associated to $\mu$ and we
let $T_c$, $\aone_0$, $\{\aone_i\}$, $\{E_i\}$ and $\{d_i\}$ be as in
Lemma~\ref{axes and lifts}.  Thus $\mu$ is either $E_j$ for some $j$
or an element of the \qe\ family determined by $E_j \bar E_{j'}$ for
some $j$ and $j'$.

Letting $\ti u$ be the path such that $\ti f(\ti E) = \ti E\cdot \ti
u$, we have $\ti R = \ti E \cdot \ti R_0$ where $\ti R_0 = \ti u \cdot
\ti f_\#(\ti u) \cdot \ti f^2_\#(\ti u)\ldots$.  Since $\ti E$ is not
linear, $\ti \mu$ occurs infinitely often as a term in the \ds\ of
$\ti R_0$, where we do not distinguish between elements of the same
\qe\ family of subpaths.  There is a completely split subpath $\ti
\sigma_0$ of $\ti R_0$ and a coarsening $\ti \sigma_0 = \ti \tau_1
\cdot \ti \mu \cdot \ti \tau_2$ of the \ds\ of $\sigma_0$ where $\ti
\mu$ is a lift of $\mu$ and where $ \tau_1$ and $ \tau_2$ are not
Nielsen paths.  Denote the initial and terminal endpoints of $\ti
\sigma_0$ by $ \ti a_0$ and $\ti b_0$ and for $l \ge 1$, let $\ti
\sigma_l = \ti f^l_\#(\ti \sigma_0)$, $\ti a_l = \ti f^l(\ti a_0)$,
and $\ti b_l = \ti f^l(\ti b_0)$.  Then
\begin{itemize}
\item[(1)] $\ti \sigma_l \subset \ti R_0$ and $\ti \sigma_l \to P$.
\end{itemize}
Let $\ti f_j$ be the lift of $f$ corresponding to $\aone_j$ and let
$\ti E_j$ be a lift of $E_j$ whose initial endpoint is fixed by $\ti
f_j$ and whose terminal endpoint is contained in $A_c$.  There is a
covering translation $S_0: \Gamma \to \Gamma$ such that $\ti E_j$ is
the initial edge of $S_0(\ti \mu)$.  Let $\ti \delta_0 = S_0(\ti
\sigma_0)$. For $l \ge 1$, let $S_l:\Gamma \to \Gamma$ be the covering
translation such that $\ti f_j^l S_0 = S_l \ti f^l$, let $\ti \delta_l
= S_l(\ti \sigma_l)$ and let $\ti x_l$ and $\ti y_l$ be the endpoints
of $\ti \delta_l$.  It is immediate that
\begin{itemize}
\item[(2)] $\ti E_j \subset \ti \delta_l$.
\item[(3)]$\ti \delta_l = {\ti {f_j}^l}_\#(\ti \delta_0)$.
\item[(4)] the length of $\ti \delta_l \cap A_c$ goes to infinity with
$l$.
\end{itemize}

Lemma~\ref{limit exists} applied to $\ti f_j$ and $S_0(\ti \tau_1)$
implies that
\begin{itemize}  
\item[(5)] $\ti x_l \to Q_- \in \Fix_N(\hat \aone_j) \setminus
\{T_c^{\pm}\}$.
\end{itemize}

If $\mu$ corresponds to $E_j$, let $m=d_j$ and $t=0$.  If $\mu$
corresponds to $E_j\bar E_{j'}$, let $m=d_j-d_{j'}$ and $t=j'$.  Then
the terminal endpoint of $S_0(\ti \mu)$ is fixed by $T_c^{-m}\ti f_j$.
Lemma~\ref{limit exists} applied to $\ti T_c^{-m}\ti f_j$ and $S_0(\ti
\tau_2)$ implies that
\begin{itemize}  
\item[(6)] $T_c^{-ml}\ti y_l \to Q_+ \in \Fix_N(\hat \aone_t)
\setminus \{T_c^{\pm}\}$.
\end{itemize}
The maximal principal sets $X_j = \Fix(\aone_j)$ and $X_t =
\Fix(\aone_t)$ contain $T_c^{\pm}$ and determine lifts $s_j,
s_t:\A(\oone) \to \Aut(F_n)$.

We have so far only focused on $\oone$.  We now bring in $\otwo$.  Let
$g : G' \to G'$ be a representative of $\otwo$ satisfying the
conclusions of Theorem~\ref{comp sp exists} and let $\ti g, \ti g_j$
and $\ti g_t$ be lifts of $g$ to the universal cover $\Gamma'$
corresponding to $\atwo=s(\otwo)$,  $\atwo_j = s_j(\otwo)$ and
$\atwo_t= s_t(\otwo)$ respectively.  The following are equivalent.
\begin{itemize}
\item  $\omega(\otwo) = 0$.
\item  $\atwo_j = \atwo_t$.
\item   $Q_+  \in \Fix(\hat \atwo_j)$.
\end{itemize}
 It suffices to show that $P$ is isolated in $\Fix(\hat \atwo)$ if and
only if $Q_+ \not \in \Fix(\hat \atwo_j)$.

To compare points in $\Gamma$ and $\Gamma'$, choose an equivariant map
$h : \Gamma \to \Gamma'$ that preserves the markings; equivalently,
when $\partial \Gamma$ and $\partial \Gamma'$ are identified with
$\partial F_n$ then $\hat h : \partial \Gamma \to \partial \Gamma'$ is
the identity.  Let $C$ be the bounded cancellation constant for $h :
\Gamma \to \Gamma'$ and let $\ti R' = h_\#(\ti R)$.  We use prime
notation for covering translations and axes of $\Gamma'$.  Thus $S_l':
\Gamma' \to \Gamma'$ is the covering translation such that $S_l' h = h
S_l$.  Denote $h(\ti a_l)$, $h(\ti b_l)$ and the path that they bound
by $\ti a_l',\ti b_l'$ and $\ti \sigma_l'$.  Let $\ti x_l' = S'_l(\ti
a'_l) =S_l'h(\ti a_l) = hS_l(\ti a_l) = h(\ti x_l)$, let $\ti y_l' =
S'_l(\ti b'_l) =h(\ti y_l)$ and let $\ti \delta_l' = S'_l(\ti
\sigma'_l) =h_\#(\delta_l)$ be the path connecting $\ti x_l'$ to $\ti
y_l'$.  We have
 
\begin{itemize}
\item[($1'$)] $\ti \sigma_l'$ is $C$-close to $R'$ and $\ti \sigma_l'
\to P$.
   
\item[($4'$)] the length of $\ti \delta'_l \cap A'_c$ goes to infinity
with $l$.
\item[($5'$)] $\ti x'_l \to Q_- \in \Fix(\hat \atwo_j) \setminus 
\{{T'}_c^{\pm}\}$.
\item[($6'$)] ${T'}_c^{-ml}\ti y'_l \to Q_+ \in \Fix(\hat
\atwo_t)\setminus \{ {T'}_c^{\pm}\}$.
\end{itemize}

If $P$ is not isolated in $\Fix(\hat \atwo)$ then Lemma~\ref{l: second
from bk3} implies, after increasing $C$ if necessary, that $\ti a_l'$
and $\ti b'_l$ are $C$-close to $\Fix(\ti g)$ for all sufficiently
large $l$.  After replacing $\ti a'_l$ and $\ti b'_l$ with $C$-close
elements of $\Fix(\ti g)$, replacing $\sigma'_l$ with the path
connecting the new values of $\ti a'_l$ and $\ti b'_l$, and replacing
$C$ by $2C$, properties $(1'),(4'), (5')$ and $(6')$ still hold and
each $\sigma_l'$ is a Nielsen path for $g$.  Since $\ti \delta_l'$ is
a lift of $\sigma_l'$, Lemma~\ref{for the non-isolated case} implies
that for all sufficiently large $l$, the lift of $g$ that fixes $\ti
x'_l$ and $\ti y'_l$ commutes with $T'_c$ and so equals ${T'}_c^{d_l}
\ti g_j$ for some $d_l$.  Since $Q_- \in \Fix(\hat g_j)$ there is a
neighborhood of $Q_-$ in $\Gamma'$ that is disjoint from
$\Fix({T'}_c^m \ti g_j)$ for all $m \ne 0$.  Since $\ti x'_l \to Q_-$,
it follows that $d_l = 0$ and hence that $\ti y_l' \in \Fix(\ti g_j)$
for all sufficiently large $l$.  Since $ \Fix(\ti g_j)$ is
$T_c'$-invariant, ${T'}_c^{- ml} \ti y_l' \in \Fix(\ti g_j)$ and so
$Q_+ \in \Fix(\hat g_j)$ as desired.

Suppose then that $P$ is isolated in $\Fix(\hat{\atwo})$.  After
replacing $\ti a_l'$ and $\ti b_l'$ by their nearest points in $\ti
R'$, we may assume that $\ti \sigma_l' \subset R'$ and that properties
$(1'),(4'), (5')$ and $(6')$ still hold.  Lemma~\ref{iterates to}
implies that there is a non-linear edge $\ti E'$ that iterates toward
$P$ under the action of $\ti g$.  Denoting $g^m_\#(E')$ by $\rho_m$ we
have that for all sufficiently large $l$ there exists $m > 0$ such
that $\sigma_l'$ is a subpath of $\rho_m$.  There is a lift $\ti
\rho_m$ of $\rho_m$ that contains $\ti \delta_l'$ and so has endpoints
$\partial_{\pm}\ti \rho_m$ such that $\partial_{-}\ti \rho_m \to Q_-$
and ${T'}_c^{-ml}\partial_{+}\ti \rho_m \to Q_+$.  The former implies
that for sufficiently large $m$, the initial endpoints of $\ti \rho_m
\cap A'_c$ and $\ti {g_j}_\#(\ti \rho_m) \cap A'_c$ are equal and
the latter implies that if $Q_+ \in \Fix(\hat g_j) = \Fix(\hat
\atwo_j)$ then the terminal endpoints of $\ti \rho_m \cap A'_c$ and
$\ti {g_j}_\#(\ti \rho_m) \cap A'_c$ are equal.  On the other hand,
$\ti \rho_m \cap A'_c$ and $\ti g{_j}_\#(\ti \rho_m) \cap A'_c$ have
different lengths by Lemma~\ref{increasing length} so we conclude that
$Q_+ \not \in \Fix(\atwo_j)$.
\qed

\section{Abelian Subgroups of Maximal Rank} \label{s:maximal rank}
By Theorem~\ref{disintegration groups are enough}, all abelian
subgroups are realized, up to finite index, as subgroups of some
$\D_R(\phi)$.  In this section we describe those $\phi$ for which
$\D_R(\phi)$ has maximal rank.  As usual, $\phi$ is represented by a
\rtt\ $\fG$ and filtration $\filt$ satisfying the conclusions of
Theorem~\ref{comp sp exists}.

For the simplest example, start with $G_2$ having one vertex $v_1$,
two edges $E_1$ and $E_2$ and with $f$ defined by $f(E_1) = E_1$ and
$f(E_2) = E_2E_1^{m_1}$ for some $m_1 \in \Z$.  For $k=1,\ldots,n- 2$,
add pairs of linear edges, $E_{2k+1}$ and $E_{2k+2}$, initiating at a
new common vertex $v_{k+1}$, terminating at $v_1$ and satisfying
$f(E_j)=E_j E_1^{m_j}$ for distinct $m_j$.  Thus $G$ has $2n-3$ linear
edges and the resulting $\D_R(\phi)$ has rank $2n-3$, which is known
\cite{cv:moduli} to be maximal.  In this example all edges terminate
at the same vertex and there is only one axis, but this is just for
simplicity. One could, for example, take the terminal vertex of $E_5$
equal to $v_2$ and define $f(E_5) = E_5w_5$ where $w_5$ is a closed
Nielsen path based at $v_2$.  Similar modifications can be done to the
other edges as well.

Another simple modification is to redefine $f|G_2$ so that $G_2$ is a
single EG stratum with Nielsen path $\rho$ and redefine $f$ on the
other edges to be linear with axis represented by $\rho$.  We may view
the original example as being built over a Dehn twist of the punctured
torus and this modification as being built over a pseudo-Anosov
homeomorphism of the punctured torus.

A perhaps more surprising example of a maximal rank abelian subgroup
is constructed as follows.  Let $S$ be the genus zero surface with
four boundary components $\beta_1,\dots,\beta_4$ and let $h :S \to S$
be a homeomorphism that represents a pseudo-Anosov mapping class and
that pointwise fixes each $\beta_m$.  Let $A$ be an annulus with
boundary components $\alpha_1$ and $\alpha_2$ and with its central
circle labeled $\alpha_3$. Define $D_{jk}:A \to A$ to be the
homeomorphism that restricts to a Dehn twist of order $j$ on the
subannulus bounded by $\alpha_1$ and $\alpha_3$ and to a Dehn twist of
order $k$ on the subannulus bounded by $\alpha_2$ and $\alpha_3$.
Finally, define $Y = S \cup A/\sim$ where $\sim$ identifies $\alpha_m$
to $\beta_m$ for $1 \le m \le 3$.  The homeomorphisms $g_{ijk} : Y \to
Y$ induced by $h^i$ and $D_{jk}$ for $i,j,k \in \Z$ define a rank
three abelian subgroup $\A'$.  The fundamental group of $Y$ is a free
group of rank three and the image of $\A'$ in $\Out(F_3)$ is an
abelian subgroup $\A$ of maximal rank.

We present a slight generalization of this example in terms of
relative train tracks as follows.

\begin{ex} \label{rank three maximal}
Suppose that $G$ is a rank three marked graph  with vertices
$v_1,\dots, v_4$, that $\emptyset = G_0 \subset G_1 \subset
\cdot\ldots\cdot \subset G_4 = G$ is a filtration and that $f : G \to
G$ is a relative train track map such that
\begin{itemize}
\item $G_1 $ is a single fixed edge $E_1$ with both ends attached to
$v_1$.
\item for $m=2,3$, $H_m$ is a single edge $E_m$ with terminal endpoint
$v_1$ and initial endpoint $v_m$; $f(E_m) = E_m E_1^{d_m}$ where $d_2$
and $d_3$ are distinct non-zero integers.
\item $H_4$ is an EG stratum with three edges, one connecting $v_4$ to
$v_l$ for each $l=1,2,3$; for each edge $E$ of $H_4$, $f(E)$ is a
concatenation of edges in $H_{4}$ and Nielsen paths of the form
$E_1^*$, $E_{2} E_1^* \bar E_{2}$ and $E_{3} E_1^* \bar E_{3}$.
\end{itemize}
Then $f$ determines an element $\oone \in \Out(F_3)$ such that
$\D_R(\phi)$ has rank three.  The example described above using a four
times punctured sphere is a special case of this construction with
$\ast = \pm1$.  In general, $H_4$ is not a geometric stratum in the
sense of \cite{bfh:tits1}.
\end{ex} 

We think of the strata $H_2 \cup H_3 \cup H_4$ in Example~\ref{rank
three maximal} as being a single unit added on to the lower filtration
element, which in this case is a single circle.  If the lower
filtration element has higher rank then we have the option of adding
an additional linear edge. In the geometric case this amounts to Dehn
twisting on three of components of the four times punctured sphere
instead of just two.  We formalize this as follows, where the acronym
FPS is chosen to remind the reader of the four times punctured sphere.

\begin{notn}
We say that $H_{l+1} \cup H_{l+2} \cup H_{l+3}$ is a {\em partial
\fpss} if
\begin{enumerate}
\item [(1)]There are (not necessarily distinct) closed Nielsen paths
$\alpha_1,\alpha_2 \subset G_{l}$.
\item [(2)] For $j=1,2$ the stratum $H_{l+j}$ is a single linear edge
$E_{l+j}$ such that $f(E_{l+j}) = E_{l+j}\alpha_j^{d_j}$ for some
non-zero $d_j$.  The initial endpoints $v_{l+j}$ of $E_{l+j}$ are
distinct and are not contained in $G_l$. (Equivalently, $G_{l+2}$
deformation retracts to $G_l$.)
\item [(3a)] $H_{l+3}$ is \eg.
\item [(3b)] $H_{l+3} \cap G_{l+2} = \{v_{l+1},v_{l+2},v\}$  for some vertex $v \in G_l$.
\item [(3c)] $H_{l+3}$ is either a pair of arcs joined at a common
endpoint or a {\em triad}, which is three arcs joined at a common endpoint.  whose valence three vertex is not contained in $G_{l+2}$.
For each edge $E$ of $H_{l+3}$ the edge path $f(E)$ is a concatenation
of edges in $H_{l+3}$, Nielsen paths of the form $E_{l+j} \alpha_j^*
\bar E_{l+j}$ for $j =1,2$ and Nielsen paths in $G_l$.
\end{enumerate}
\end{notn}

\begin{remark}\label{first fpss}
In Example~\ref{rank three maximal}, $H_2 \cup H_3 \cup H_4$ is
partial \fpss.
\end{remark}

\begin{notn} \label{n:fpss}
We say that $H_{l+1} \cup \dots \cup H_{l+4}$ is an {\em \fpss} if
\begin{enumerate}
\item [(1)]There are (not necessarily distinct) closed Nielsen paths
$\alpha_1,\alpha_2, \alpha_3 \subset G_{l}$.
\item [(2)] For $j=1,2,3$ the stratum $H_{l+j}$ is a single linear
edge $E_{l+j}$ such that $f(E_{l+j}) = E_{l+j}\alpha_j^{d_j}$ for some
non-zero $d_j$.  The initial endpoints $v_{l+j}$ of $E_{l+j}$ are
distinct and are not contained in $G_l$. (Equivalently, $G_{l+3}$
deformation retracts to $G_l$.)
\item [(3a)] $H_{l+4}$ is \eg.
\item [(3b)] $H_{l+4} \cap G_{l+2} = \{v_{l+1},v_{l+2},v_{l+3}\}$.
 \item [(3c)] $H_{l+4}$ is either a pair of arcs joined at a common
endpoint or a triad whose valence three vertex is not contained in
$G_{l+3}$.  For each edge $E$ of $H_{l+4}$ the edge path $f(E)$ is a
concatenation of edges in $H_{l+4}$ and Nielsen paths of the form
$E_{l+j} \alpha_j^* \bar E_{l+j}$ for $j =1,2,3$.
\end{enumerate}
\end{notn}

\begin{com}
If $H_{l+1} \cup H_{l+2} \cup H_{l+3}$ is a is partial \fpss\ then
$\chi(G_{l}) - \chi(G_{l+3}) =2$.  If $H_{l+1} \cup \dots \cup
H_{l+4}$ is a \fpss\ then $\chi(G_{l}) - \chi(G_{l+4}) =2$.
\end{com}  

We can now state the main results of this section.  We assume the
existence of $\fG$ satisfying the conclusions of Theorem~\ref{comp sp
exists} applied without reference to a particular $\F$ and satisfying
the additional condition that there are no non-trivial invariant
forests.  This is always possible by Remark~4.7 of
\cite{fh:recognition}.  Partial \fpss s arise in the proof of the
propositions but not in their statements.

\begin{prop}\label{max abel}
Suppose that $\oone \in\Out(F_n)$ is rotationless, that
$\D_R(\phi)$ has rank $2n-3$ and that $\oone$ is represented by $\fG$
and $\filt$ as in Theorem~\ref{comp sp exists}.  Suppose further that
$\fG$ has no non-trivial invariant forests.  Then after
reordering the filtration if necessary, there are $1 \le l_1 < \dots <
l_K \le N$ such that
\begin{description}
\item [(A)]$G_{l_1}$ either:
\begin{enumerate}
\item has rank two and is a single \eg\ stratum.
\item has rank two and consists of two edges $E_1,E_2$ where $f(E_1) = E_1$ and 
$f(E_2) = E_2E_1^m$ for some $m\ge 1$.
\item has rank three and $f|G_{l_1}$ is as in Example~\ref{rank three maximal}.
\end{enumerate}
\item [(B)]for $m> 1$, $\cup_{j=l_m+1}^{l_{m+1}}H_{i_j}$ is either 
\begin{enumerate}
\item a pair of linear edges with a common initial vertex that is not contained in $G_{l_m}$ or
\item an \fpss.
\end{enumerate}
\end{description}
\end{prop}

There is an analogous result for abelian subgroups of the subgroup
$\IA$ of $\Out(F_n)$ consisting of elements that act trivially in
homology.

\begin{prop}\label{IA max abel}
Suppose that $\oone\in\Out(F_n)$ is rotationless, that
$\D_R(\phi) \subset \IA$ has rank $2n-4$ and that $\oone$ is
represented by $\fG$ and $\filt$ as in Theorem~\ref{comp sp exists}.
Suppose further that $\fG$ has no non-trivial invariant forests.
Then after reordering the filtration if necessary, there are $1 \le
l_1 < \dots < l_K \le N$ such that
\begin{description}
\item [(A)] $l_1=2$ and $G_{2}$ is connected, has rank two and is
contained in $\Fix(f)$.
\item [(B)]for $m> 1$, $\cup_{j=l_m+1}^{l_{m+1}}H_{i_j}$ is either 
\begin{enumerate}
\item a pair of linear edges with homologically trivial axes and a
common initial vertex that is not contained in $G_{l_m}$.
\item an \fpss\ with homologically trivial axes.
\end{enumerate}
\end{description}
\end{prop} 

Recall that one uses the \ds\ of the $f$-image of edges of $G$ to
define almost invariant subgraphs $X_1,\ldots, X_M$ of $G$ and that if
$a_i$ is a non-negative integer assigned to $X_i$ then
$(a_1,\ldots,a_M)$ is admissible (Definition~\ref{consistency}) if it
satisfies certain linear relations involving three of the $a_i$'s.
The rank of $D_R(\oone)$ is equal to the rank of the subspace of
$\R^M$ generated by the admissible $M$-tuples for $\fG$.

For induction purposes, it is useful to consider admissible sequences
of $f|G_j$ for each filtration element $G_j$.  Let $M_j$ be the number
of \ias s for $f|G_j$ and let $R_j$ be the rank of the subspace of
$\R^{M_j}$ generated by the admissible $M_j$-tuples defined with
respect to $f|G_j$.  Each \ias\ for $f|G_j$ is contained in a \ias\
for $f|G_{j+1}$.  Every relation on $M_j$-tuples determines a relation
on $M_{j+1}$-tuples by inclusion. Thus every admissible
$M_{j+1}$-tuple for $f|G_{j+1}$ \lq restricts\rq\ to an admissible
$M_j$-tuple for $f|G_j$.  The only \ias\ of $f|G_{j+1}$ that can
contain more than one \ias\ for $f|G_j$ is the one that contains
$H_{j+1}$.  Amalgamating \ias s of $f|G_j$ into a single \ias\ for
$f|G_{j+1}$ can be viewed as a finite set of new relations.  All other
new relations involve the \ias\ containing $H_{j+1}$.

To consolidate ideas and as a warm up we prove a simple estimate on
the $R_j$'s.  A path $\sigma$ is {\em pre-Nielsen} if it is not a
Nielsen path but $f_\#^k(\sigma)$ is a Nielsen path for some $k \ge
1$.

In the following lemmas $G$ is as in Propositions 8.6 and
8.7.

\begin{lemma} \label{simple rank calculation}  If $G$ is as in Propositions 8.6 or
8.7 then the following hold for $r < s$ and $\Delta R = R_s-R_r$.
\begin{enumerate}
\item  If $H_s$ is a fixed edge and $r=s-1$ then $\Delta R = 0$.
\item If $H_{s}$ is a linear edge and $r=s-1$ then $\Delta R =1$.
\item If $H_{s}$ is a non-linear \noneg\ edge and $r=s-1$ then $\Delta
R \le 0$.
\item If $H_r$ is irreducible, $H_s$ is EG and if all strata between
$H_r$ and $H_s$ are zero strata then $\Delta R \le 1$ with equality if
and only if for each edge $E$ of $H_{s}$ the terms in the \ds\ of
$f(E)$ are either edges in $H_{s}$, pre-Nielsen connecting paths or
Nielsen paths in $G_r$.
\end{enumerate}
\end{lemma}

\proof (1) is obvious since there are no new \ias s and no new
relations.  For the remaining items it is clear that $\Delta R \le 1$
since there is at most one new \ias. $\Delta R = 1$ when
$G_{s}\setminus G_r$ is an entire \ias\ for $f|G_{s}$ and when that
\ias\ is not part of any relation.  Equivalently, for each edge $E$ of
$H_{s}$ the terms in the \ds\ of $f(E)$ are either edges in $H_{s}$,
pre-Nielsen connecting paths or Nielsen paths in $G_r$.  (2), (3) and
(4) follow immediately.  \endproof

We now come to the main step in the proofs of
Proposition~\ref{max abel} and Proposition~\ref{IA max abel}.

\begin{lemma} \label{euler} Assume that  $G$ is as in Propositions 8.6 or
8.7,   that $r \le u < s$ and that the following conditions are
satisfied.
\begin{enumerate}
\item [(1)] $G_r$ and $G_s$ have no valence one vertices.
\item [(2)] For $r < l \le u$, \ $H_l$ is a single edge $E_l$ whose
terminal vertex is in $G_r$ and whose initial vertex has valence one
in $G_{u}$.
\item [(3)] For $ u< l < s$, \ $H_l$ is a zero stratum.
\item [(4)] $H_s$ is an EG stratum.
\end{enumerate}  
Let $H_{rs} = \cup_{l=r+1}^sH_l$, let $\Delta R = R_s-R_r$ and let
$\Delta\chi = \chi(G_r) - \chi(G_s)$.  If there is a vertex $v \in
G_r$ and a fixed direction at $v$ determined by an edge of $H_{rs}$
let $\delta=1$; otherwise $\delta= 0$.

Then
$$ \Delta R \le 2 \Delta \chi - \delta.
$$ Moreover, if the inequality is an equality then one of the
following holds:
\begin{description}

\item [(a)]  $H_{rs}$ is  an \fpss\ and $\delta = 0$.
\item [(b)]  $H_{rs}$ is  a partial \fpss\ and $\delta = 1$.
\end{description}
\end{lemma}

\proof If $H_s$ is disjoint from $G_u$ then $s = r+1$ and $H_s$ is a
component of $G_s$.  In this case, $\delta=0$, $\Delta R = 1$, $\Delta
\chi \ge 1$ and the lemma is clear.  We assume for the remainder of
the proof that $H_s$, and hence each component of $H_s$, has non-empty
intersection with $G_u$.

Item (4) and Corollary~3.2.2 of \cite{bfh:tits1}  imply
that
$$
\Delta \chi \ge 2.   
$$
Let $H_{us} = \cup_{l=u+1}^{s}H_l$.  Thus $G_s =
G_{u} \cup H_{us}$.  Denote $G_{u} \cap H_{us}$ by $\V$, the
cardinality of $\V$ by $V$ and the number of components in $ H_{us}$
by $C_{us}$.  Then
$$
\Delta\chi \ge V - C_{us}
$$
with equality if and only if each component of $H_{us}$ is contractible
and 
$$
\Delta\chi \ge C_{us}
$$ with equality if and only if each component of $H_{us}$ is
  topologically either an arc whose interior is disjoint from $G_u$ or
  a loop that intersects $G_u$ in a single point.  Thus
$$ 2\Delta\chi \ge V
$$ with equality if and only if each component of $H_{us}$ is
topologically an arc whose interior is disjoint from $G_u$.  On the
other hand, if each component of $H_{us}$ is topologically an arc
whose interior is disjoint from $G_u$ then there are no illegal turns
in $H_s$.  This would imply the existence of $m > 0$ so that for any
loop $\gamma \subset G_s$ that intersects $H_s$ non-trivially, the
number of edges of $H_s$ in $f^m_\#(\gamma)$ would be strictly larger
than the number of edges of $H_s$ in $\gamma$.  This can not be true
as one easily sees by considering the loops $\gamma_{-m}$ satisfying
$f_\#^m(\gamma_{-m}) = \gamma$.  We conclude that
$$ 2\Delta\chi > V.
$$

For $r < l \le u$, the stratum $H_l$ is a single edge $E_l$.  We write
$E_l \in \E_L$ if $E_l$ is linear.  The initial endpoints $\V_L$ of
the edges in $\E_L$ have valence one in $G_{u}$.  We denote the
cardinality of $\V_L$ by $V_L$.  Lemma~\ref{simple rank calculation}
implies that
$$
\Delta R \le V_L + 1.
$$ Note also that if $V = V_L$ then $\delta = 0$.  Thus $V -\delta \ge
V_L$ and
$$
\Delta R \le V_L + 1 \le V +1 -\delta.
$$ If $C_{us} = 1$ then $\Delta\chi \ge V -C_{us} = V-1$ and
$\Delta\chi \ge 2$  imply that
$$ 
2 \Delta\chi \ge  V+1
$$
with equality only if $\Delta\chi = 2$ and $V=3$.  Thus
$$
2 \Delta\chi - \Delta R -\delta  \ge  (V+1)- (V+1-\delta)-\delta
     = 0
$$ with equality only if $\Delta\chi = 2$, $V=3$ and $\Delta R +
\delta = 4$.  To complete the proof in the $C_{us}=1$ case, assume
that equality holds.  In the $\delta=0$ case, $3 = V \ge V_L \ge
\Delta R-1 = 3$ so $V= V_L =3$ and $\Delta R = 4$.  Items (1) and (2)
of Notation~\ref{n:fpss} follow from the fact that $V = V_L =3$.
Since $\Delta\chi = V-C_{us}$,\ $H_{us}$ is contractible; being
attached to $G_{u}$ in three places, it is topologically either a
triad or a pair of arcs joined at a point. In the former case the
unique valence three vertex of $H_{us}$ must be the base point of both
a legal turn in $H_s$ and an illegal turn in $H_s$ and so is disjoint
from $G_{s-1}$.  The elements of $\V_L$ are fixed points and so are
not contained in any zero strata.  The remaining vertices in $H_s$, if
any, have valence two and have their links entirely contained in
$H_s$.  There is no loss in erasing these vertices.  Once this done,
there are no zero strata in $H_{us}$ so $s =u+1 =r+4$. Items (3a) and
(3b) are immediate and (3c) follows from Lemma~\ref{simple rank
calculation}(4) and Theorem~\ref{comp sp exists}-(N).  We have now
verified that equality in the $\delta=0$ case corresponds to case (a)
of the lemma.  In the $\delta=1$ case, $3 = V > V_L \ge \Delta R-1 =
2$.  The same argument shows that this corresponds to case (b).  This
completes the proof when $C_{us} =1$.

For $C_{us} > 1$ we need another estimate on $\Delta R$, namely
$$
C_{us}  > 1 \Rightarrow \Delta R +\delta \le V + 2 - C_{us} .
$$ From this, the lemma is completed by
\begin{eqnarray*}
2 \Delta\chi - \Delta R -\delta & > & V- (V+2-C_{us} )\\
     &=& C_{us} -2\\
     &\ge&  0.
\end{eqnarray*}
It remains therefore to prove the estimate and for this we show that
there are enough relations between the \ias\ that contains $H_{us}$
and the \ias s determined by the linear edges corresponding to $\V_L$.

Let $X_1,\ldots,X_d$ be the components of $H_{us}$ that are disjoint
from $G_r$ and intersect $\V$ in a subset of $\V_L$.  Define a graph
$Z$ with one vertex $z_p$ for each $X_p$ and one additional vertex $z$
representing $H_s \cup G_r$.  There is at most one edge connecting any
pair of vertices.  The edges of $Z$ are defined as follows.  Suppose
that $\mu$ is an edge of $H_s$ or a connecting path in $H_{us}$ and
that there is a term $\nu \subset G_{u}$ in the \ds\ of $f(\mu)$ that
has exactly one endpoint in $X_p$. After reversing the orientation on
$\nu$ if necessary, the initial edge of $\nu$ is a linear edge $E'_p
\subset X_p$ and either $\nu = E'_p$ or $\nu$ is a \qe\ subpath.  If
$\nu = E'_p$ then $E'_p$ belongs to the same \ias\ as $H_s$. If $\nu$
is \qe\ and the terminal edge of $\nu$ is contained in $G_r$ then
there is a linear relation between the coefficients associated to the
\ias\ containing $E'_p$, the \ias\ containing $H_s$ and the \ias\
containing a stratum of $G_r$.  In both of these cases, $Z$ has
exactly one edge connecting $z_p$ to $z$.  Otherwise, $\nu$ is \qe\
and the terminal edge of $\nu$ is some $\bar E'_q \subset X_q$.  In
this case, there is a linear relation between the coefficients
associated to the\ias\ containing $E'_p$, the \ias\ containing $H_s$
and the \ias\ containing $E'_q$; $Z$ has exactly one edge connecting
$z_p$ to $z_q$.

Let $\va$ be an admissible $M_s$-tuple.  If there is an edge
connecting $z_p$ to $z$ then the coordinate of $\va$ corresponding to
the \ias\ containing $E^*_p$ is determined by the coordinates of $\va$
corresponding to the \ias\ containing $H_s$ and to the \ias s
containing the strata of $G_r$.  Thus one does not need to count both
$E'_p$ and $H_s$ when estimating $\Delta_iR$.  Similarly, if $z_p$ and
$z_q$ belong to the same component of $Z$ then the coordinate of $\va$
corresponding to the \ias\ containing $E'_p$ is determined by the
coordinates of $\va$ corresponding to the \ias\ containing $H_s$ and
to the \ias\ containing $E'_q$.  Thus one does not need to count both
$E'_p$ and $E'_q$ when estimating $\Delta R$.  In both cases each edge
of $Z$ allows us to improve our estimate $\Delta R \le V_L + 1$ by
lowering the right hand side by one.

Let $ |Z| $ be the number of edges in $Z$ and note that $V-V_L \ge
C_{us} -d$.  Thus
\begin{eqnarray*}
\Delta R +\delta & \le & V_L + 1 -|Z|  + \delta \\
     &=& V- (V-V_L) -|Z|  +1 + \delta\\
     &\le&  V - C_{us} +d  - |Z| + 1 + \delta
\end{eqnarray*}
and it suffices to show that $d + \delta -|Z| \le 1$.

Choose an edge $E$ and $m \ge 1$ so that $f^m_\#(E)$ crosses every
edge of $H_{us}$.  Terms in the \ds\ of $f^m_\#(E)$ that are not
single edges in $H_s$ or connecting paths in $H_{us}$ are subpaths of
$G_{u}$.  Consider those terms $\tau \subset G_{u}$ that intersect
some $X_p$ in exactly one endpoint and in particular are not loops.
After reversing orientation if necessary, the initial edge of $\tau$
is a linear edge $E'_p$ in $X_p$ and $\tau$ is either a single edge or
a \qe\ path.  In the former case $\tau$ itself, and in the latter
case, some element of the quasi-exceptional family that contains
$\tau$, is a term in the complete splitting of $f_\#(\mu)$ where $\mu$
is either an edge of $H_s$ or a connecting path in $H_{us}$.  In
either case $\tau$ determines an edge in $Z$.  Since $f^m_\#(E)$
crosses every edge in $H_{us}$, all of the vertices $z_p$
corresponding to components of $H_{us}$ are contained in the same
component of $Z$.  If $z$ is also in this component then $|Z| \ge d $
and we are done.  If $z$ is not in this component then every $\tau$ as
above is a \qe\ path with initial and terminal endpoints in components
of $H_{us}$ represented by vertices of $Z$.  It follows that these are
all the components of $H_{us}$ and hence that $\delta=0$.  Thus $|Z|
-\delta = |Z| \ge d-1 $ and we are done. \endproof

\noindent{\bf Proof of Theorem~\ref{max abel}} After reordering the
strata of the filtration we may assume that there exists $1 \le k < K$ and $k =l_0 <l_1 <\dots<l_K = N$ 
such that  
\begin{itemize}
\item $H_j$ is a non-contractible component of $G_j$ if and only if $j
\le k$.  In particular, $G_k$ has no valence one vertices.
\end{itemize}
and such that the following hold for all $1 \le i \le K$
\begin{itemize} 
\item    $G_{l_i}$  has no valence one vertices. In
particular, $G'_{l_{i}}$ does not deformation retract to
$G'_{l_{i-1}}$.  
\item If $l_{i-1} \le j < l_{i}$ and $H_j$ is irreducible then $G_j$
deformation retracts to $G_{l_{i-1}}$.
\end{itemize}

For $1\le i \le K$,  let $\Delta_iR= R_{l_{i}}-R_{l_{i-1}}$, let
$\Delta_i \chi = \chi(G_{l_{i-1}}) - \chi(G_{l_{i}})$, let $\hat H_i =
\cup_{j=l_{i-1}+1}^{l_{i}} H_j$ and let $\delta_i =1$ if there is a
vertex $v \in G_{l_{i-1}}$ and a fixed direction at $v$ determined by an
edge of $\hat H_i$ and $\delta_i= 0$ otherwise.

The following sublemma is an easy extension of Lemma~\ref{euler}.  We
separate it out of the proof for easy reference.

\begin{sublem} \label{four cases}
Assume notation as above.  For $1 \le i \le K$,  
$$ \Delta_iR \le 2 \Delta_i\chi - \delta_i
$$
with equality only if one of  the following holds.  
\begin{description}
\item [(a)]  $\hat H_{i}$ is  an \fpss\ and $\delta_i = 0$.
\item [(b)]  $\hat H_{i}$ is  a partial \fpss\ and $\delta_i = 1$.
\item [(c)] $\hat H_i$ is a single linear edge and $\delta_i = 1$.
\item [(d)] $\hat H_{i}$ is a pair of linear edges with a common
initial vertex and $\delta_i = 0$.
\end{description}
\end{sublem}

\proof If $\hat H_i$ contains an \eg\ stratum $H_j$ then $j=l_{i}$
because $G_j$ does not deformation retract to $G_{l_{i-1}}$.  In this
case, the sublemma follows from Lemma~\ref{euler}.  If no stratum
$H_j$ of $\hat H_i$ is \eg\ then each $H_j$ is a single edge $E_j$ and
$l_{i}$ equals $l_{i-1}+1$ or $l_{i-1}+2$.  In both cases $\Delta_i
\chi = 1$.  If $l_{i} = l_{i-1}+1$ then $\delta_i = 1$; if $l_{i} =
l_{i-1}+2$ then $\delta_i=0$.  Lemma~\ref{simple rank calculation}
implies that if $l_{i} = l_{i-1}+1$ then $\Delta_iR \le 1$ with
equality corresponding to (c) and that if $l_{i} = l_{i-1}+2$ then
$\Delta_iR \le 2$ with equality corresponding to (d).  \endproof

The sublemma implies that 
$$ R_K- R_k =\sum_{i=1}^K\Delta_iR \le \sum_{i=1}^K (2\Delta_i\chi
 -\delta_i) = 2\chi_k -2\chi_K -\sum_{i=1}^K\delta_i.
$$

Denote $|\chi(G_k)|$ by $c$.  Since each component of $G_k$ is a
single stratum and the restriction of $f$ to a rank one component of
$G_k$ is the identity,
$$
R_k \le c
$$ with equality if and only each component of $G_k$ has rank one or
 two.  Thus
\begin{eqnarray*}
2n-3- R_k & \le & 2(n-1-c)  -\sum_{i=1}^K\delta_i\\
     &=& 2n-2-2c -\sum_{i=1}^K\delta_i
\end{eqnarray*}
 implies that 
$$
c \le 2c-R_k \le 1 -\sum_{i=1}^K\delta_i.
$$

If $c= 1$ then each $\delta_i = 0$, which implies by Theorem~\ref{comp
sp exists}-(PER) that no component of $G_k$ has rank one.  It follows
that $G_k = G_1$ has rank two. Moreover, each inequality in the above
displayed equations are equalities. The proposition in this case now
follows from the sublemma.

If $c=0$ then there is at most one non-zero $\delta_i$.  Together
these imply that $G_k = G_1$ has rank one and that there is exactly
one non-zero $\delta_i$.  If $\delta_1 = 1$ then the sublemma
completes the proof; see Remark~\ref{first fpss}.  Suppose then that 
$\delta_i =1$ for $i > 2$. We will modify $\fG$, arranging that $\delta_1 =1$ for the new homotopy equivalence.

Let $v$ and $E_1$ be the unique vertex and edge in $G_1$ and let $E_i$
be the edge in $H_i$ that determines a fixed direction pointing out of
$G_1$.  Since $\delta_1 = 0$, $G_2$ is a single edge $E_2$ satisfying
$f(E_2) =E_2 E_1^{d_2}$ for some $d_2 \ne 0$.  The link $L(G,v)$ of
$v$ in $G$ consists of both ends of $E_1$, the initial end of $E_i$
and the terminal ends of some linear edges $E_j$, including $E_2$.
For each such $E_j$ there is a non-trivial closed path $u_j \subset
G_{j-1}$ such that $f(E_j) = E_j u_j$.

Create a new graph $G'$ by replacing $v$ with a pair of vertices $v_1$
and $v_2$, adding a new edge $E_0$ with one endpoint at $v_1$ and the
other at $v_2$ and partitioning $L(G,v)$ into $L(G',v_1) \cup
L(G',v_2)$ as follows.  Both ends of $E_1$ belong to $L(G',v_1)$.  The
initial endpoint of $E_i$ is in $L(G',v_2)$.  If the terminal end of
$E_j$ is contained in $L(G,v)$ then it is assigned to $L(G',v_1)$ if
$u_j = E_1$ and to $L(G',v_2)$ otherwise.

The map $f':G' \to G'$ induced by $f:G \to G$ satisfies the
conclusions of Theorem~\ref{comp sp exists} but has an invariant
forest, namely the single edge $E_0$.  Orient $E_0$ so that $v_1$ is
its terminal edge and define $u_0$ to be the trivial path at $v_1$.
With the exception of $E_1$, each element of $L(G,v_1)$ is the
terminal endpoint of an edge $E_k$ satisfying $f(E_k) = E_kE_1^{d_k}$.
Define a new homotopy equivalence $g :G' \to G'$ by replacing $d_k$
with $d_k-d_2$.  Note that $f'$ and $g$ are freely homotopic and so
represent the same element of $\Out(F_n)$.  We have changed the
invariant forest from $E_0$ to $E_2$.  Finally, modify $g$ by
collapsing $E_2$ to a point.  We are now back to the case that
$\delta_2=1$.  \qed

\vspace{.1in}

\noindent{\bf Proof of Theorem~\ref{IA max abel}} The proof is a
variation on that of Theorem~\ref{max abel}.  No changes are required
in the proof up through the sublemma so we do not repeat that here.
The rest of the proof follows.

The sublemma implies that 
$$
 R_K- R_k  =\sum_{i=1}^K\Delta_iR \le \sum_{i=1}^K (2\Delta_i\chi -\delta_i)
=  2\chi_k -2\chi_K -\sum_{i=1}^K\delta_i.
$$

Denote $|\chi(G_k)|$ by $c$.  Since each component of $G_k$ 
is a single stratum  and the restriction of $f$ to a rank one component of $G_k$ 
is the identity,  
$$
R_k \le c
$$
 with equality if and only each component of $G_k$ has rank one or two.  
Thus 
\begin{eqnarray*}
2n-4- R_k & \le & 2(n-1-c)  -\sum_{i=1}^K\delta_i\\
     &=& 2n-2-2c -\sum_{i=1}^K\delta_i
\end{eqnarray*}
which implies that 
$$
c \le 2c-R_k \le 2 -\sum_{i=1}^K\delta_i.
$$

If some component of $G_k$ has rank three then $R_k < c$ and the last
displayed inequality is strict $c < 2 -\sum_{i=1}^K\delta_i \le 2$
which is impossible.  Thus each component of $G_k$ has rank one or
two.  Since \ IA$_2$ is trivial \cite{nielsen:out(f2)}, no component of $G_k$ can have
rank two.  (Rank two fixed subgraphs exist but they are not composed
of a single stratum.)

We may therefore assume that $G_k$ is a union of fixed circles. Each
of these circles represents a non-trivial homology class and so can
not be the axis associated to any linear edge.  We claim that $f$
restricts to the identity on the component $C$ of $G_{l_{1}}$ that
has rank greater than one.  If not, then $C$ is obtained from a fixed
circle by adding one EG stratum and perhaps some zero strata.  In this
case, $R_{l_{1}} \le 1$ while $C$ has rank at least three by
Lemma~3.22 of \cite{bfh:tits1}.  It follows that $2 \chi(G_{l_{1}}) -
R_{l_{1}} \ge 3$.  The sublemma then implies that $2\chi(G_K) - R_K
\ge 3$ which contradicts the fact that $\chi(G_K) = n-1$ and $R_K =
2n-4$.  This verifies the claim.  This same argument proves that $C$
has rank two and that $k =1$.  The sublemma completes the proof.
\endproof

\section{Two Families of Abelian Subgroups}  \label{s:for commensurator}
We now return to the simplest examples of maximal rank abelian
subgroups, those that are rotationless, that have linear growth and
that have only one axis.  We prove that these subgroups and their
standard generators can be characterized using only algebraic (as
opposed to dynamical systems) properties.  These results are needed in
the calculation \cite{FarbH:commensurator} of the commensurator of
$\Out(F_n)$.

We begin by relating the rank of $A(\psi)$ to the dynamical properties
of $\psi$ in a special case.

\begin{lemma} \label{only one stratum}
Suppose that $A$ is a maximal rank rotationless abelian subgroup of
$\Out(F_n) $ or $\IA$, that $\psi \in A$ and that $\A(\psi)$ has rank
one.  Then either $\L(\psi)$ has exactly one element and $\psi$
has no axes or $\L(\psi) = \emptyset$ and $\psi$ has exactly one axis
and that axis has multiplicity one.
\end{lemma}  

\proof It suffices to show that $\psi$ has a representative $g : G'
\to G'$ with exactly one non-fixed stratum.

By Lemma~\ref{synthetic is everything} there exists a rotationless
$\oone$ so that $A \subset \A(\oone)$.  Choose $\fG$ and $\filt$ that represent $\oone$ as in
Theorem~\ref{comp sp exists}.  Theorem~\ref{D has finite
index} implies that $D_R(\oone) \cap D_R(\oone^{-1})$ has finite index
in $\A(\oone) = \A(\oone^{-1})$.  After replacing $\psi$ with an
iterate, we may assume that $\psi \in D_R(\oone) \cap
D_R(\oone^{-1})$.

Proposition~\ref{max abel} and Proposition~\ref{IA max abel} imply
that if $\Lambda \in \L(\oone)$ corresponds to an EG stratum $H_r$
then both ends of every leaf of $\Lambda$ intersect $H_r$ infinitely
often.  By Lemma~3.1.15 of \cite{bfh:tits1} each leaf of $\Lambda$ is dense in
$\Lambda$.  In other words $\Lambda$ is minimal.  The symmetric
argument applied to $\oone^{-1}$ shows that every element of
$\L(\oone) \cup \L(\oone^{-1})$, and hence (Lemma~\ref{is generic})
  every element of $\L(\psi)$, is minimal.

We next prove that there is no proper free factor system $\F$ that
carries each element of $\L(\phi) \cup \L(\phi^{-1})$ and each
$\phi$-invariant conjugacy class by assuming that there is such an
$\F$ and arguing to a contradiction.  By Proposition~\ref{comp sp
exists} there exists $f':G' \to G'$ and $\filt$ representing $\phi$
such that $\F = \F(G_r)$ for some $G_r$.  Let $p:G \to G'$ be the
quotient map that collapses each invariant tree to a point.  The
homotopy equivalence $f':G' \to G'$ induced by $\fG$ satisfies the
conclusions of Theorem~\ref{comp sp exists} and has no invariant
forests.  Proposition~\ref{max abel} and Proposition~\ref{IA max abel}
imply that $p(G_r) = G'$. Equivalently, $f|(G\setminus G_r)$ is the
identity. But this implies that $D_R(\phi|\F)$ has the same rank as
$D_R(\phi)$, which contradicts the fact that the maximal rank of an
abelian subgroup of $\Out(F_n)$ [resp. $\IA_n$] is strictly larger
than of a proper free factor system of $\Out(F_n)$ [resp. $\IA_n$].

Since $\L(\psi) \subset \L(\phi) \cup \L(\phi^{-1})$ and since each
$\oone$-invariant conjugacy class is also $\psi$-invariant
(Lemma~\ref{still Nielsen}), no proper free factor system carries each
element of $\L(\psi) $ and each $\psi$-invariant conjugacy class.
These two facts, the minimality of elements of $\L(\psi)$ and the
absence of a free factor system as above, imply that each non-fixed
stratum in a representative $g:G' \to G'$ of $\psi$ is either linear
or EG and that the rank of $D_R(\psi)$ is equal to the number of
non-fixed strata of $g:G' \to G'$.  As this also equals the rank of
$\A(\psi)$ the lemma follows.  \endproof

Let $G$ be the the rose with $n-2$ of its edges subdivided into
two edges.  Thus there are edges $E_1,\ldots, E_{2n-2}$ and vertices
$v_1,\ldots,v_{n-1}$ with $v_1$ the terminal vertex of all edges and
the initial edges of $E_1$ and $E_2$ and with $v_k$ the initial vertex
of $E_{2k-1}$ and $E_{2k}$ for $2 \le k \le n-1$.

For $1 \le i \le 2n-3$, define $f_i : G \to G$ by $E_{i+1} \mapsto
E_{i+1} E_1$.  Choose a basis $x_1,\ldots,x_n$ for $F_n$ and a marking
on $G$ that identifies $x_j$ with the $j^{th}$ loop of $G$.  The
elements $\eta_i \in \Out(F_n)$ determined by $f_i$ are a basis for an
abelian subgroup $A_1$ of rank $2n-3$.  If $i = {2k-2}$ for $k \ge 2$
then $\hat \eta_i$ is defined by $x_k \mapsto x_k x_1$.  If $i =
{2k-1}$ for $k \ge 2 $ then $\hat \eta_i$ is defined by $x_k \mapsto
\bar x_1x_k$.  The remaining element $\hat \eta_1$ is defined by $x_2
\mapsto x_2 x_1$.  Borrowing notation from \cite{FarbH:commensurator}
we say that $A_1$ is the {\em type E subgroup associated to the basis
$x_1,\ldots,x_n$} and that $\eta_1,\ldots, \eta_{2n-3}$ are its {\em
standard generators}.

\begin{com} 
It is not hard to check (see for example Lemma~2.13 of
\cite{FarbH:commensurator}) that $\eta_1$ is conjugate to each
$\eta_j$ and to $\eta_j\eta_l$ if $\{ j,l\} \ne \{2k,2k+1\}$.
Corollary~\ref{in the center} and Lemma~4.5 \ of
\cite{FarbH:commensurator} imply that $\A(\eta_1)$ has rank one.  This
explains the hypothesis in the next lemma.
\end{com}   

\begin{lemma} \label{elementary uniqueness}
Suppose that $\phi_1,\ldots,\phi_{2n-3}$ are a basis for an abelian
subgroup $A$ of $\Out(F_n)$, $n \ge 3$, that each $\A(\phi_j)$
 has rank one and that $\A(\phi_j\phi_l)$ has
rank one if $\{j,l\} \ne \{2k,2k+1\}$.  Then there is a basis
$x_1,\ldots,x_n$ for $F_n$, standard generators $\eta_j$ of the type E
subgroup associated to this basis, and $s,t > 0$ so that $\phi_j^s =
\eta_j^t$.
\end{lemma}

\proof After replacing each $\phi_j$ with $\phi_j^s$ for some fixed $s
\ge 1$ we may assume by Lemma~\ref{synthetic is everything} that there
exists a rotationless $\theta\in A$ so that each $\phi_j \in
\A(\theta)$.  Choose $\fG$ representing $\theta$ as in
Theorem~\ref{comp sp exists}.  The coordinates of $\Omega^{\theta} :
\A(\theta) \to \Z^{2n-3}$ (Definition~\ref{better Omega}) are in one
to one correspondence with the non-fixed irreducible strata of $\fG$
representing $\theta$ and so correspond to linear edges and EG strata
as described in Proposition~\ref{max abel}.

For each $\phi_j$ there exists $\phi_{j'} \ne \phi_j$ such that
$\A(\phi_j\phi_{j'})$ has rank one.

 Suppose that $\psi \in \A(\theta)$, that $\omega_i$ is a coordinate
of $\Omega^{\theta}$ and that $\omega_i(\psi) \ne 0$.  If $\omega_i =
\PF_{\Lambda}$ then $\Lambda \in \L(\theta) \cup \L(\theta^{-1})$ by
Corollary~3.3.1 of \cite{bfh:tits1}.  If $\omega_i$ corresponds to a linear edge with associated axis
$[c]_u$ then $[c]_u$ is an axis for $\psi$; if $\omega_r$ also
corresponds to a linear edge with associated axis $[c]_u$ and if
$\omega_i(\psi) \ne \omega_r(\psi)$ then $[c]_u$ is an axis for $\psi$
with multiplicity greater than one.  Lemma~\ref{only one stratum}
therefore implies that for each $\phi_j$ the coordinates of
$\Omega^{\theta}(\phi_j)$ takes on a single non-zero value and that if
more than one coordinate takes this value then all such coordinates
come from linear edges associated to the same axis.  The same holds
true for the coordinates of $\Omega^{\theta}(\phi_j\phi_{j'})$.

Suppose that $\omega_i = \PF_{\Lambda}$ and that $\omega_i(\phi_j) \ne
0$.  At least one of $\omega_i(\phi_{j'})$ or
$\omega(\phi_{j'}\phi_j)$ is non-zero , say $\omega_i(\phi_{j'})$.
Then $\Omega^{\theta}(\phi_j)$ and $\Omega^{\theta}(\phi_{j'})$ are
contained in a cyclic subgroup of $\Z^{2n-3}$ in contradiction to the
fact that $\phi_j$ and $\phi_{j'}$ generate a rank two subgroup and
the injectivity of $\Omega^{\theta}$.  We conclude that each
coordinate of $\Omega^{\theta}$ corresponds to a linear edge of $\fG$.

Minor variations on this argument show that all linear edges
correspond to the same axis, that only one coordinate of
$\Omega^{\theta}(\phi_j)$ can be non-zero and that the non-zero value
$t$ that is taken is independent of $j$.  The details are left to the
reader. The lemma now follows from the explicit description of $\fG$
given by Proposition~\ref{max abel} and the definition of
$\Omega^{\theta}$.  \endproof
  
There is an analogous result for $\IA$.  For the model subgroup, we
use the same marked graph $G$ as in the definition of type E
subgroups.  Choose a closed path in $G_2$ based at $v_1$ that forms a
circuit and determines a trivial element of homology.  For $1 \le i
\le 2n-4$ define $f_i : G \to G$ by $E_{i+2} \mapsto E_{i+2}w$.  The
elements $\mu_{i,w} \in \Out(F_n)$ determined by $f_i$ are a basis for
an abelian subgroup $A_w$ of $\IA$ with rank $2n-4$.  We think of $w$
as both a path in $G_2$ and an element of the free factor $\langle
x_1, x_2\rangle$.  If $i = {2k-5}$ then $\hat \mu_{i,w}$ is defined by
$x_k \mapsto x_k w$ and if $i = {2k-4}$ then $\hat \eta_i$ is defined
by $x_k \mapsto \bar w x_k$.  Borrowing notation from
\cite{FarbH:commensurator} we say that $A_w$ is the {\em type C
subgroup associated to $w$ and to the basis $x_1,\ldots,x_n$} and that
$\mu_{1,w},\ldots, \mu_{2n-4,w}$ are its {\em standard generators}.

\begin{lemma} \label{IA uniqueness}
Suppose that $\phi_1,\ldots,\phi_{2n-4}$ are a basis for an abelian
subgroup of $\IA$, $n \ge 4$, that each $\A(\phi_j)$  has rank one and that $\A(\phi_j\phi_l)$ has rank one if
$\{j,l\} \ne \{2k,2k+1\}$.  Then there exists a basis $x_1,\ldots,x_n$
for $F_n$, a homologically trivial element $w \in \langle x_1,
x_2\rangle$ and standard generators $\eta_j$ of the type C subgroup
associated to $w$ and this basis, and $s,t > 0$ so that $\phi_j^s =
\eta_j^t$
\end{lemma}

\proof We have assumed that $n \ge 4$ so that for all $\phi_j$ there
exists $\phi_{j'}$ such that $\A(\phi_j\phi_{j'})$ has rank one.
Otherwise the proof of Lemma~\ref{elementary uniqueness} carries over
to this context without modification, $w$ representing the unique axis
of the elements $\phi_j$.  \endproof

\bibliographystyle{amsalpha}

\bibliography{mh}

\end{document}